\documentclass[twoside,11pt]{article}

\usepackage{jmlr2e}
\RequirePackage[OT1]{fontenc}
\usepackage[utf8]{inputenc}
\let\proof\relax

\RequirePackage{amsthm,amsmath,natbib,graphicx,enumerate}
\usepackage{amssymb,tabularx,multicol,multirow,booktabs}
\usepackage{mhequ}
\usepackage{relsize}
\usepackage[top=1in, bottom=1in, left=1in, right=1in]{geometry}
\usepackage[toc,page]{appendix}
\usepackage[english]{babel}
\RequirePackage{xr}
\usepackage{mathtools}
\usepackage{bbm}
\allowdisplaybreaks

\def \cov{\mathit{cov}}
\def \Cov{\mathit{Cov}}

\def \BP{\mathbb{P}}
\def \TE{\mathrm{E}}
\def \bx{\mathbf{x}}
\def \bX{\mathbf{X}}

\def \be{\begin{equs}}
\def \ee{\end{equs}}

\def \E{\mathbb{E}_0}
\def \P{\mathbb{P}_0}

\def \bhat{\hat{b}}

\def \sumn{\sum_{i=1}^n}

\def \ihat {\hat{I}}

\def \k0 {k_0}

\def \k1 {k_1}
\def \k2 {k_2}
\def \k3 {k_3}
\def \jhat {\hat{j}}
\def \Pf {\mathbb{P}_f}
\def \Ef {\mathbb{E}_f}
\def \Pp {\mathbb{P}_P}

\def \lprime {l^{'}}
\def \fmax {\|f\|_{\infty}}

\def \n2 {\frac{1}{n(n-1)}}

\def \xij {(X_i,X_j)}

\def \In {\mathbf{I}_n}

\def \khat {\hat{k}}

\def\aa{\check a}
\def\bb{\check b}

\let\tilde\widetilde

\def \da {\delta a}
\def \db {\delta b}

\def \phihat {\hat{\phi}}
\def \lhat {\hat{l}}
\def \Ltilde{\tilde{L}}
\def \lprime {l^{'}}

 \def\Par{{\mathbf{\mathcal{P}} }}

 \def \gmin {\gamma_{\min}}
 \def \gmax {\gamma_{\max}}

 \def \lmin {l_{\min}}
 \def \lmax {l_{\max}} 
 \def \gtilde {\tilde{g}}
  \def \ftilde {\tilde{f}}
 \def \lhat {\hat{l}}
 \def \lstar {l^*}
 \def \khat{\hat{k}}
 \def \kustar{k^*}
 \def \ahat {\hat{a}}
 \def \bav {\frac{\alpha+\beta}{2}}
 \def \dela {\delta a}
 \def \delb {\delta b}
 \def \bpsi {\boldsymbol{\psi}}

 \def \Uhat{\hat{U}}
 \def \IF{\textit{IF}}
 
\theoremstyle{plain}
\newtheorem{thm}{Theorem}[section]
\numberwithin{equation}{section}
\newtheorem{prop}[thm]{Proposition}

\theoremstyle{remark}
\newtheorem{rem}[thm]{Remark}
\newtheorem{cor}[thm]{Corollary}
\newtheorem{lem}[thm]{Lemma}

\def \bx{\mathbf{x}}
\def \bX{\mathbf{X}}

\def \be{\begin{equs}}
	\def \ee{\end{equs}}

\def \E{\mathbb{E}}
\def \P{\mathbb{P}}

\def \fhat{\hat{f}}

\def \d{\{1,\ldots,d\}}

\def \sumn{\sum_{i=1}^n}

\def \fracn {\frac{1}{n}}

\def \xij {x_{ij}}

\def \g2 {\|g\|_2}

\def \bw {\mathbf{w}}
\def \bW {\mathbf{W}}

\def \sumn{\sum_{i=1}^n}

\def \ihat {\hat{I}}

\def \k0 {k_0}

\def \k1 {k_1}
\def \k2 {k_2}
\def \k3 {k_3}
\def \jhat {\hat{j}}
\def \Pf {\mathbb{P}_f}
\def \Ef {\mathbb{E}_f}

\def \bpsi {\boldsymbol{\psi}}

\def \bpsilk {\boldsymbol{\psi}_{l,k}}

\def \lprime {l^{'}}
\def \fmax {\|f\|_{\infty}}

\def \n2 {\frac{1}{n(n-1)}}

\def \xij {(X_i,X_j)}

\def \In {\mathbf{I}_n}

\def \khat {\hat{k}}

\def\aa{\check a}
\def\bb{\check b}

\let\tilde\widetilde

\def \da {\delta a}
\def \db {\delta b}

\def \pnbar {\overline{P}_n}
\def \qnbar {\overline{Q}_n}
\def \plambda {p_{\lambda}}
\def \qlambda {q_{\lambda}}
\def \pnlambda {P_{\lambda}^n}
\def \qnlambda {Q_{\lambda}^n}
\def \constant {C(B,C_U,J_0)}
\def \bW {\mathbf{O}}
\def \bw {\mathbf{o}}
\def \sumk {\sum\limits_{k \in \Z_j}}
\def \sumv {\sum\limits_{v\in \{0,1\}^d}}
\def \psijkv {\psi_{jk}^v}
\def \constant {C(B,B_U,J_0)}
\def \bW {\mathbf{O}}
\def \bw {\mathbf{o}}
\def \sumk {\sum\limits_{k \in \Z_j}}
\def \sumv {\sum\limits_{v\in \{0,1\}^d}}

\def \jmin {j_{\mathrm{min}}}
\def \jmax {j_{\mathrm{max}}}

\def \ghat {\hat{g}}
\def \jstar {j^{*}}

 \def\a{{\alpha}}
 \def\b{{\beta}}                         
                                        
 \def\d{{\delta}}

 \def\g{{\gamma}}                       \def\GG{\mathbb{G}}                           
                           
\def\I{{\cal I}}                                         
                                                      
 \def\k{{\kappa}}

                 \def\n{{\nu}}                                                        
 
                           \def\PP{\mathbb{P}}
                        
\def\R{{\mathbb{R}}} 				
                 \def\s{{\sigma}}

 				 \def\bX{{\bf X}}
					 
\def\Z{{\cal Z}}                         \def\ZZ{\mathbb{Z}}

\def\half{{\frac12}}

\def\ra{\rightarrow}
\def\da{\downarrow}

\def\dotfil{\leaders\hbox to 1em{\hss.\hss}\hfill}

\def\hexnumber#1{\ifcase#1 0\or 1\or 2\or 3\or 4\or 5\or 6\or 7\or 8\or 9\or A\or B\or C\or D\or E\or F\fi}

\def\BL{{\rm B\kern -0.5pt L}}
\def\UC{{\rm U\kern-0.5pt C}}

\def\lin{\mathop{\rm lin\,}\nolimits}

\def\closedhull{{\overline{\hull}}}

\def\hull{\mathop{\rm conv\ }\nolimits}

\def\ball{\mathop{\rm ball}\nolimits}

\def\Prodint{\mathop{\mathchoice{\vcenter{\epsfysize=1.5em\epsfbox{%
  prodint.eps}}}{\vcenter{\epsfysize=1em\epsfbox{prodint.eps}}}{}{}}}

\def\tdot#1{\kern1.2pt\dot{\vphantom{#1}}\kern-1.2pt
        \dot#1\kern0.8pt\dot{\vphantom{#1}}\kern-0.8pt}

\def\E{\mathord{\rm E}}

\def\Var{\mathop{\rm Var}\nolimits}

\def\cov{\mathop{\rm cov}\nolimits}
\def\Cov{\mathop{\rm Cov}\nolimits}

\def\almosteverywhere{{\rm a.e.}}

\mathchardef\given="626A
\mathcode`:="603A

\def\generalweak#1{\ {\mathchoice{\buildrel #1\over \rightsquigarrow}%
{\raise-2pt\hbox{$\buildrel #1\over \rightsquigarrow$}}{}{}}\ }

\def\generalprob#1{\ {\mathchoice{\raise-1.5pt\hbox{$\buildrel#1\over\ra$}}
{\raise-2pt\hbox{$\buildrel#1\over\ra$}}{}{}}\ }

\def\Pgg{\ {\mathchoice{\raise-1.5pt\hbox{$\buildrel {\small P}\over\gg$}}
{\raise-2pt\hbox{$\buildrel {\small P}\over\gg$}}{}{}}\ }

\def\endskip{}
\def\boxit#1{\vbox{\hrule\hbox{\vrule\kern1pt\vbox{\kern1pt #1\kern1pt}\kern1pt\vrule}\hrule}}
\def\eet#1{}

\usepackage{lastpage}
\jmlrheading{22}{2021}{1-\pageref{LastPage}}{10/19; Revised 4/21}{5/21}{19-892}{Lin Liu, Rajarshi Mukherjee, James M. Robins, Eric Tchetgen Tchetgen}
\ShortHeadings{Adaptive Estimation of Nonparametric Functionals}{Liu, Mukherjee, Robins, Tchetgen Tchetgen}

\firstpageno{1}

\begin{document}

\title{Adaptive Estimation of Nonparametric Functionals}

\author{\name Lin Liu \email linliu@sjtu.edu.cn \\
       \addr Institute of Natural Sciences, MOE-LSC, School of Mathematical Sciences, SJTU-Yale Center for Biostatistics and Data Science \\
       Shanghai Jiao Tong University \\
       Shanghai, 200240, China
       \AND
       \name Rajarshi Mukherjee \email ram521@mail.harvard.edu \\
       \addr Department of Biostatistics \\
       Harvard T. H. Chan School of Public Health \\
       Boston, MA 02115, USA
       \AND
       \name James M. Robins \email robins@hsph.harvard.edu \\
       \addr Department of Epidemiology and Department of Biostatistics \\
       Harvard T. H. Chan School of Public Health \\
       Boston, MA 02115, USA
       \AND
       \name Eric Tchetgen Tchetgen \email ett@wharton.upenn.edu \\
       \addr Wharton Statistics Department \\
       The University of Pennsylvania \\
       Philadelphia, PA 19104, USA}

\editor{G\'{a}bor Lugosi}

\maketitle

\begin{abstract}
We provide general adaptive upper bounds for estimating nonparametric functionals based on second-order U-statistics arising from finite-dimensional approximation of the infinite-dimensional models. We then provide examples of functionals for which the theory produces rate optimally matching adaptive upper and lower bounds. Our results are automatically adaptive in both parametric and nonparametric regimes of estimation and are automatically adaptive and semiparametric efficient in the regime of parametric convergence rate.
\end{abstract}


\section{Introduction} \label{sec:intro}
Estimation of functionals of data generating distribution has always been of central interest in statistics. In nonparametric statistics and machine learning problems, where data generating distributions are parametrized by functions in infinite-dimensional spaces, there exists a comprehensive literature addressing such questions. In particular, a large body of research has been devoted to explore minimax estimation of linear, quadratic functionals, and entropy type functionals in density and/or white noise models. We do not attempt to survey this extensive literature in this area. However, the interested reader can find a comprehensive snapshot of the literature in \cite{ hall1987estimation}, \cite{bickel1988estimating}, \cite{donoho1990minimax1}, \cite{donoho1990minimax2}, \cite{fan1991estimation}, \cite{kerkyacharian1996estimating}, \cite{laurent1996efficient}, \cite{johnstone2001chi}, \cite{cai2003note}, \cite{cai2004minimax}, \cite{cai2005nonquadratic}, \cite{kandasamy2014influence}, \cite{berrett2019efficient}, \cite{han2020optimal, han2020estimation} and other references therein. Although the question of more general nonparametric functionals has received relatively less attention, some fundamental insights regarding estimation of non-linear integral functionals in density and white noise models can be found in \cite{ibragimov2013statistical}, \cite{kerkyacharian1996estimating}, \cite{nemirovski2000topics}, and references therein. 

A general feature of the results obtained while estimating ``smooth" nonparametric functionals is an elbow effect in the rate of estimation based on the regularity of the underlying function classes. For example while estimating quadratic functionals in a $d$ dimensional density model, $\sqrt{n}$-efficient estimation can be achieved as soon as H\"{o}lder exponent $\beta$ of the underlying density exceeds $\frac{d}{4}$, whereas the optimal rate of estimation is $n^{-\frac{4\beta}{4\beta+d}}$ (in root mean squared error sense) for $\beta \leq \frac{d}{4}$. A similar elbow in the rate of estimation exists for estimation of non-linear integral functionals as well. For density model this was demonstrated by \cite{birge1995estimation}, \cite{kerkyacharian1996estimating} and \cite{tchetgen2008minimax}. For signal or white noise model, the problem of general integrated non-linear functionals was studied by \cite{nemirovski2000topics}, but mostly in the $\sqrt{n}$-regime. However, for more complex nonparametric models, the approach for constructing minimax optimal procedures for general non-linear functionals  in non-$\sqrt{n}$-regimes has been rather case specific. Motivated by this, in recent years, \cite{robins2008higher, robins2017minimax} and \cite{mukherjee2017semiparametric} have developed a theory of inference for nonlinear functionals in parametric, semi-parametric, and non-parametric models based on higher order influence functions. 

Most minimax rate optimal estimators proposed in the literature, however, depend explicitly on the knowledge of the smoothness indices. Thus, it becomes of interest to understand the question of adaptive estimation i.e. the construction and analysis of estimators without prior knowledge of the smoothness. The question of adaptation of linear and quadratic functionals has been studied in detail in the context of density, white noise, and nonparametric additive Gaussian noise regression models (\cite{low1992renormalization}, \cite{efromovich1994adaptive},  \cite{efromovich1996optimal}, \cite{tribouley2000adaptive}, \cite{johnstone2001chi}, \cite{efromovich2000adaptive}, \cite{klemela2001sharp}, \cite{laurent2000adaptive}, \cite{cai2005adaptive1}, \cite{cai2006optimal},  \cite{gine2008simple}, \cite{breunig2021adaptive}). Although our theoretical results are very much related to and motivated by the papers on optimal adaptive estimation of quadratic functionals, including \cite{fan1991estimation, johnstone2001chi, efromovich2000adaptive, laurent2000adaptive, cai2006optimal, gine2008simple}, nonetheless we consider a much wider and more complex class of functionals whose first order influence function depends on multiple nuisance functions. Our paper is the first to consider adaptive estimation in this wider class.

In particular, we observe i.i.d copies of a random vector $O = (\mathbf{W}; \bX) \in \mathbb{R}^{m + d}$ with unknown distribution $P$ on each of $n$ study subjects. The variable $\bX$ represents a random vector of baseline covariates such as age, height, weight, etc. Throughout $\bX$ is assumed to have compact support and a density with respect to (w.r.t.) Lebesgue measure in $\mathbb{R}^d$. The variable $\mathbf{W}\in \mathbb{R}^m$ can be thought of as a combination of outcome and treatment variables in some of our examples. In the above setup, we are interested in estimating certain ``smooth" functionals $\phi(P)$ in the sense that under finite-dimensional parametric submodels, they admit first-order derivatives which can be represented as inner products of first order influence functions with score functions \citep{bickel1993efficient}. For some classical examples of these functionals, we provide matching upper and lower bounds on the rate of adaptive minimax estimation over a varying class of smoothness of the underlying functions, provided that the unknown marginal design density of $\bX$ is sufficiently regular.

The contributions of this paper are as follows. Building on the theory of adaptive estimation of linear and quadratic nonparametric functionals in density and Gaussian white noise models, we explore adaptation theory for non-linear functionals in more complex nonparametric models in both the $\sqrt{n}$ (where our adaptive estimators are semi-parametric efficient) and non-$\sqrt{n}$-regimes. The crux of our arguments relies on the observation that when the non-adaptive minimax estimators can be written as a sum of empirical mean type statistics and $2^{\mathrm{nd}}$-order U-statistics, one can provide a unified theory of selecting the ``best" data driven estimator using Lepski type arguments \citep{lepskii1991problem,lepskii1992asymptotically}. Indeed, under certain assumptions on the data generating mechanism $P$, the non-adaptive minimax estimators have the desired structure for a large class of problems \citep{robins2008higher}. This enables us to produce a class of examples where a single method helps provide a desired answer. Although the basic scheme of the approach is straightforward, the proof is somewhat involved because of the dependence of the U-statistics kernels in question on estimated functions from the sample. This prohibits performing a Hoeffding decomposition w.r.t. the whole sample and we can only proceed by a conditional Hoeffding decomposition restricted on carefully chosen good events. In order to prove a lower bound for the rate of adaptation under the non-$\sqrt{n}$-regime for all the examples covered below (and show that a sharp poly-logarithmic penalty is indeed necessary for adaptation), we adapt the results in \cite{birge1995estimation, robins2009semiparametric} for the Hellinger distance between two mixtures of suitable product measures, by producing similar results in chi-square divergence (see Appendix \ref{app:equiv}). Our results apply to several common examples in nonparametric analyses of observational studies. These include, but are not limited to, average treatment effect estimation in causal inference, mean estimation in missing data studies, and error variance estimation in general non-parametric problems. To the best of our knowledge, the optimal adaptive results (both in sense of upper and lower bound) which simultaneously describe both $\sqrt{n}$ and non-$\sqrt{n}$ regimes of estimation, are among the first results in this direction. Moreover, all these examples involve estimating nonparametric nuisance functions (such as a density or regression). Consequently, the proofs of our results shed light on what properties of machine learning algorithms one needs to explore (see e.g. the results in Theorem \ref{thm_adaptive_density_regression} proved for wavelet based estimation) so that a principled use of them in semiparametric theory is possible. 


The rest of the paper is organized as follows. In Section \ref{section_main_results} we provide the main results of the paper in a general form. Section \ref{section_examples} is devoted for applications of the main results in specific examples. Discussions on potential issues and future work are provided in Section \ref{section_discussion}. In Section \ref{section_wavelets and function spaces} we provide a brief discussion on some basic wavelet and function space theory and notations, which we use extensively. Finally Section \ref{section_appendix_proof_of_theorems}, Appendices \ref{section_appendix_hellinger}, \ref{section:proof_of_remaining_theorems}, and \ref{section_appendix_technical_lemmas} are devoted for the proofs of the theorems and collecting useful technical lemmas.

\subsection{Notation}\label{section:notation} 
For data arising from underlying probability distribution $P$ we denote by $\P_P$ and $\E_P$ the probability of an event and expectation under $P$ receptively. In the sequel, we often split the whole sample $\mathcal{D}$ into $M$ disjoint subsamples $\{\mathcal{D}_1, \mathcal{D}_2, \ldots, \mathcal{D}_M\}$. For some subset $\mathcal{J}_m = \{ j_1, j_2, \ldots, j_m \} \subseteq \{1, 2, \ldots, M\}$, we denote by $\P_{P, \{ j_1, j_2, \ldots, j_m \} }$ and $\E_{P, \{ j_1, j_2, \ldots, j_m \} }$ the probability and expectation under $P$ over the sample space of $\cup_{j \in \mathcal{J}_m} \mathcal{D}_j$, while treating the subsamples $\cup_{j \not\in \mathcal{J}_m} \mathcal{D}_j$ as fixed. For example, when we divide the whole sample into three disjoint subsamples $\{ \mathcal{D}_1, \mathcal{D}_2, \mathcal{D}_3 \}$, then $\P_{P, 2, 3} \equiv \P_{P, \{2, 3\}}$ denotes the conditional probability of an event, while treating the subsample $\mathcal{D}_1$ as fixed.

For a bivariate function $h(O_1,O_2)$ let $$S(h(O_1,O_2))=\frac{1}{2}\left[h(O_1,O_2)+h(O_2,O_1)\right]$$ be the symmetrization of $h$. The results in this paper are mostly asymptotic (in $n$) in nature and thus require some standard asymptotic notations.  If $a_n$ and $b_n$ are two sequences of real numbers then $a_n \gg b_n$ (and $a_n \ll b_n$) implies that ${a_n}/{b_n} \rightarrow \infty$ (and ${a_n}/{b_n} \rightarrow 0$) as $n \rightarrow \infty$, respectively. Similarly $a_n \gtrsim b_n$ (and $a_n \lesssim b_n$) implies that $\liminf_{n \rightarrow \infty} {{a_n}/{b_n}} = C$ for some $C \in (0,\infty]$ (and $\limsup_{n \rightarrow \infty} {{a_n}/{b_n}} =C$ for some $C \in [0,\infty)$). Alternatively, $a_n = o(b_n)$ will also imply $a_n \ll b_n$ and $a_n=O(b_n)$ will imply that $\limsup_{n \rightarrow \infty} \ a_n / b_n = C$ for some $C \in [0,\infty)$).

Finally we comment briefly on the various constants appearing throughout the text and proofs. Given that our primary results concern convergence rates of various estimators, we will not emphasize the role of constants throughout and rely on fairly generic notation for such constants. In particular, for any fixed  tuple $v$ of real numbers, $C(v)$ will denote a positive real number which depends on elements of $v$ only. Finally for any linear subspace $L\subseteq L_2 ([0,1]^d)$, let $\Pi\left(h|L\right)$ denote the orthogonal projection of $h$ onto $L$ under the Lebesgue measure. Also, for a function defined on $[0,1]^d$, for $1\leq q <\infty $ we let $\|h\|_q \coloneqq (\int_{[0,1]^d} |h(\bx)|^q d\bx)^{1/q}$ denote the $L_q$ semi-norm of $h$, $\|h\|_{\infty} \coloneqq \sup_{\bx \in [0,1]^d}|h(\bx)|$ the $L_{\infty}$ semi-norm of $h$. We say $h \in L_q([0,1]^d)$ for $q\in [1,\infty]$ if $\|h\|_q<\infty$. Typical functions arising in this paper will be considered to have memberships in certain H\"{o}lder balls $H(\beta,M)$ (see Section \ref{section_wavelets and function spaces} for details). This will imply that the functions are uniformly bounded by a number depending on $M$. However, to make the dependence of our results on the uniform upper bound of functions more clear, we will typically assume a bound $B_U$ over the function classes, and for the sake of compactness will avoid the notational dependence of $B_U$ on $M$. The necessary notation appeared in Lepski-type adaptation scheme will be introduced in Section \ref{section_lepski} when we first define the procedure.

\section{Main Results}\label{section_main_results} 
We divide the main results of the paper into three main parts. First we discuss a general recipe for producing a data-adaptive ``best" estimator from a sequence of estimators based on second-order U-statistics -- which in turn are constructed from compactly supported wavelet based projection kernels (defined in Section \ref{section_wavelets and function spaces}). Next we show that the bound on Hellinger distance obtained in \cite{robins2009semiparametric} directly implies the desired bound on chi-square divergence between mixtures of product measures under certain bounded density assumptions. In subsequent sections, this then serves as a basis of using a version of constrained risk inequality \citep{cai2011testing} for producing matching adaptive lower bounds in the context of estimating non-linear functionals considered in this paper. Finally we provide control over estimators of the design density as well as the regression functions in $L_{\infty}$ norm which not only adapt over H\"{o}lder type smoothness classes (defined in Section \ref{section_wavelets and function spaces}) but also belong to desired regularity classes with probability converging to $1$ sufficiently fast.

\subsection{\textbf{Upper Bound}}\hspace*{\fill}\\ \label{section_main_upper_bound}
Consider a sample of i.i.d data $\mathcal{D} \coloneqq \{ O_i = (\mathbf{W}_i, \bX_i) \sim P, \mathbf{W}_i \in \mathbb{R}^m, \bX_i \in [0, 1]^d, i = 1, \ldots, n \}$ and a real valued functional of interest $\phi(P)$. In the following, we also split $\mathcal{D}$ into two disjoint subsamples $\mathcal{D}_1 = \{ O_1, \ldots, O_{n_1} \}$ and $\mathcal{D}_2 = \{ O_{n_1 + 1}, \ldots, O_{n_1 + n_2} \}$, where $n_1 + n_2 = n$. Given this sample of size $n \geq 2$, consider further, a sequence of estimators $\{ \phihat_{n, k} (L): n \in \mathbb{N}, k = 2^{j d} \text{ for some } j \in \mathbb{N} \}$ of $\phi(P)$ (indexed by function tuple $\mathbf{L}(O)=(L_1(O),L_{2l}(O),L_{2r}(O))$) defined as follows:
\be
\phihat_{n, k} (L) = \frac{1}{n_1} \sum_{i = 1}^{n_1} L_{1}(O_i) - \frac{1}{n_1 (n_1 - 1)} \sum_{1 \leq i_1 \neq i_2 \leq n_1} S \left( L_{2l} (O_{i_1}) K_{V_k} \left(\bX_{i_1}, \bX_{i_2}\right) L_{2r} (O_{i_2}) \right)
\ee
where $K_{V_k} \left(\bX_{i_1}, \bX_{i_2}\right) \equiv K_{V_j} \left(\bX_{i_1}, \bX_{i_2}\right)$ is a resolution $j = \frac{\log_2{k}}{d}$ wavelet projection kernel defined in Section \ref{section_wavelets and function spaces}. Depending on the context, we sometimes denote $V_j$ by $V_k$ to avoid repeated translation between the resolution $j$ and the number of wavelet bases $k$ being used in defining the projection kernels. In all our examples, the non-random functions $L_1, L_{2l}, L_{2r}: \Omega \rightarrow \mathbb{R}$ depend on the true underlying distribution $P$ (e.g. $L_1(\cdot) = ( Y - \E_P (Y | \bX) ) ( A - \E_P (A | \bX) )$ in the example of average treatment effect estimation in Section \ref{sec:ate}, but we silence the dependence on $P$ in the notation for brevity) so we cannot evaluate $\phihat_{n, k} (L)$ purely from data. The (random) functions $\Ltilde_1, \Ltilde_{2l}$ and $\Ltilde_{2r}$ serve as the corresponding data-adaptive estimators of $L_1, L_{2l}$ and $L_{2r}$ constructed from the (training) subsample $\mathcal{D}_2$, without knowledge on the data generating mechanism. Then
\be
\phihat_{n, k} \equiv \phihat_{n, k} (\Ltilde) = \frac{1}{n_1} \sum_{i = 1}^{n_1} \Ltilde_{1}(O_i) - \frac{1}{n_1 (n_1 - 1)} \sum_{1 \leq i_1 \neq i_2 \leq n_1} S \left( \Ltilde_{2l} (O_{i_1}) K_{V_k} \left(\bX_{i_1}, \bX_{i_2}\right) \Ltilde_{2r} (O_{i_2}) \right),
\ee 
which is now a measurable function w.r.t. the $\sigma$-field generated by the observed data sample $\mathcal{D}$. We further assume that $$\max \{ \vert \Ltilde_1(O) \vert, \vert \Ltilde_{2l}(O) \vert, \vert \Ltilde_{2r}(O) \vert \} \leq B, \text{ $P$-almost surely}$$ for a known constant $B$. In addition, assume that $\vert g(\bx) \vert \leq B_U$ $\forall \ \bx$, $g$ being the marginal density of $\bX$ w.r.t. Lebesgue measure. Then the sequence of estimators $\{ \phihat_{n, k}: n \in \mathbb{N}, k = 2^{j d} \text{ for some } j \in \mathbb{N} \}$ of $\phi(P)$ can be thought of as a bias corrected version of a usual first order estimator arising from standard first order influence function theory for ``smooth" functionals $\phi(P)$ \citep{bickel1993efficient, van2000asymptotic}. In particular, the linear empirical mean type term $\frac{1}{n_1} \sum_{i = 1}^{n_1} \Ltilde_1(O_i)$ typically derives from the classical influence function of $\phi(P)$ and the U-statistics quadratic term 
\begin{align}\label{def:u}
\Uhat_{n, k} \coloneqq \frac{1}{n_1 (n_1 - 1)} \sum_{1 \le i_1 \neq i_2 \le n_1} S \left( \Ltilde_{2l}(O_{i_1}) K_{V_k} \left(\bX_{i_1}, \bX_{i_2}\right)\Ltilde_{2r}(O_{i_2}) \right)
\end{align} 
corrects for higher order bias terms. While, specific examples in Section \ref{section_examples} will make the structure of the sequence of estimators more clear, the interested reader will be able to find more detailed theory in \cite{robins2008higher, robins2009quadratic, li2011higher, robins2017minimax}.

The quality of such sequence of estimators will be judged against models for the data generating mechanism $P$. To this end, assume $P \in \Par_{\theta}$ where $\Par_{\theta}$ is a class of data generating distributions indexed by $\theta$ which in turn can vary over an index set $\Theta$. The choices of such a $\Theta$ will be clear from our specific examples in Section \ref{section_examples}, and will typically correspond to smoothness indices of various infinite-dimensional functions parametrizing the data generating mechanism. Further assume that there exist positive real-valued functions $f_1$ and $f_2$ defined on $\Theta$ such that the sequence of estimators $\{ \phihat_{n, k} \}_{k \geq 1}$ satisfies the following bounds for all $\theta \in \Theta$ with known constants $C_A> 0,C'_{A}$, and $C_B > 0$.

\begin{enumerate}[Property (A):]
\item Conditional Bias Bound within ``Good'' Events of $\mathcal{D}_2$: there exists a ``good'' event $\mathcal{I}_2 (n_2)$, a measurable subset w.r.t. the $\sigma$-algebra generated by the sample in $\mathcal{D}_2$, such that, for some absolute constants $C_A, C_A'$,
\be
\sup_{P \in \Par_{\theta}} \left\vert \E_{P, 1} \left( \phihat_{n, k} - \phi(P) \right) \right\vert \mathbbm{1} \left\{ \mathcal{I}_2 (n_2) \right\} \leq C_A \left( k^{-2 \frac{f_1(\theta)}{d}} + n_2^{-f_2(\theta)} \right) P\text{-}\mathrm{a.s.}
\ee
and for some $\eta > 2$,
\be
\sup_{P \in \Par_{\theta}} \P_{P, 2} \left( \overline{\mathcal{I}_2 (n_2)} \right) \leq \frac{C_A' \log{n_2}}{n_2^\eta} \quad P\text{-}\mathrm{a.s.}
\ee
where $\overline{\mathcal{I}}$ stands for the complement of $\mathcal{I}$ for any event $\mathcal{I}$.
\item Conditional Variance Bound: there exists an absolute constant $C_B > 0$, such that
\be
 \sup_{P \in \Par_{\theta}} \E_{P, 1} \left( \phihat_{n, k} - \E_{P, 1} \left( \phihat_{n, k} \right) \right)^2 \leq C_B \left( \frac{k}{n_1^2} \vee \frac{1}{n_1} \right) \quad P\text{-}\mathrm{a.s.}
\ee
\end{enumerate} 

\begin{rem}
Property (A) makes assumptions on certain nuisance functions involved in the problem. In the examples provided in Section \ref{section_examples} these will typically correspond to certain regression and density functions, and we will invoke Theorem \ref{thm_adaptive_density_regression} to show that Property (A) is met if one estimates the nuisance functions with (boundary-corrected) Cohen-Daubechies-Vial type of wavelet kernel projections\footnote{Throughout this paper, we assume that we can evaluate wavelets without numerical error and we provide the reason for making such a simplification in Remark \ref{rem:numeric}.}. 
\end{rem}

\subsubsection{\textbf{Lepski-type adaptation scheme}}\hspace*{\fill}\\ \label{section_lepski}
Given properties (A) and (B), we employ a standard Lepski-type adaptation scheme \citep{lepskii1991problem,lepskii1992asymptotically} to choose an ``optimal" estimator from the collection of estimators $\{ \phihat_{n, k}: n \in \mathbb{N}, k = 2^{j d} \text{ for some } j \in \mathbb{N} \}$. Following the notation in \citet{gine2008simple}, for any given $n_1 \in \mathbb{N}$, $n > 1$, and sufficiently small $\delta \in (0, 1)$, we define the following discretization set:
\begin{align*}
\mathcal{K} \coloneqq \left\{
k \in \left[ n_1^{1 - \delta}, \frac{n_1^2}{(\log{n_1})^4} \right]: 
k_0 = n_1^{1 - \delta}, k_1 = \frac{n_1}{\log{n_1}}, k_2 = \frac{n_1}{\ell(n_1)}, k_{j + 1} = \varrho k_j, j = 2, 3, \ldots 
\right\},
\end{align*}
where $\varrho > 1$, $\ell(n_1)$ is any function such that $\ell(n_1) \rightarrow 0$ and $\ell(n_1) \log{n_1} \rightarrow \infty$ as $n_1 \rightarrow \infty$, and $\ell^{-1}(n_1) < \log{n_1}$ for all $n_1$. Let $N$ be the cardinality of $\mathcal{K}$. Then $k_{N - 1} = \varrho^{N - 3} \frac{n_1}{\ell(n_1)} \le \frac{n_1^2}{(\log{n_1})^4}$ implies that $N = O(\log{n_1})$. Next we define $R(k) = \frac{k}{n_1^2}$ and $d(k)$ for all $k \in [n_1^{1 - \delta}, n_1^2 / (\log{n_1^4})^4]$ as
\begin{align*}
d(k) = \left\{ \begin{array}{ll}
\sqrt{ C_{\mathrm{Lepski}}^2 \log \left( \dfrac{k}{k_0} \right) } & k > k_2, \\
\ell(n_1)^{-1 / 2} & k_0 \le k \le k_2.
\end{array} \right.
\end{align*}
where $C_{\mathrm{Lepski}}$ will be specified later in the proof of Theorem \ref{theorem_general_lepski} and chosen only depending on the known parameters of the data generating mechanism. For example, when $k = \varrho^l n_1 / \ell(n_1)$ for some positive integer $l$, 
$$d(k) = \sqrt{C_{\mathrm{Lepski}}^2 (\delta \log{n_1} + l \log{\varrho} - \log{(\ell(n_1))})} = O((\log{n_1})^{1/2}).$$ 
Now following Lepski's strategy, we define the data-adaptive estimator of $k$ as
\begin{equation}
\khat \coloneqq \min \left\{ 
k \in \mathcal{K}:
\left( \hat\phi_{n, k} - \hat\phi_{n, k'} \right)^2 \le R(k') d^2(k'), \forall k' \ge k, k' \in \mathcal{K} 
\right\}.
\end{equation}
Then for each $k \in \mathcal{K}$, we simply find the corresponding $j$ as $\lfloor \frac{\log_2{k}}{d} \rfloor$ to construct estimator $\phihat_{2, k}$ with wavelets at resolution $j$.

With the notations and data-adaptive estimation scheme described above we now have the following theorem which is the main result in the direction of adaptive upper bound in this paper -- the proof of which can be found in Section \ref{section_appendix_proof_of_theorems}.
\begin{thm}\label{theorem_general_lepski}
Given a known interval $(\tau_{\min}, \tau_{\max})$ with $0\leq \tau_{\min}<\tau_{\max}$, for any $\tau \in (\tau_{\min}, \tau_{\max})$, under Properties (A) and (B) with $n_2 = n / \log{n}$ and $n_1 = n - n_2$, there exists an absolute constant $C$ depending on $(\tau_{\max}, d, C_A, C_A', C_B)$, such that
\begin{enumerate}
\item If $0 < \tau < d / 4$, then
\begin{align*}
\sup_{P \in \Par_{\theta}: \atop f_1(\theta) = \tau, f_2(\theta) > \min \left\{ \frac{4 \tau}{4 \tau + d}, \frac{1}{2} \right\}} \E_P \left( \phihat_{n, \khat} - \phi(P) \right)^2 \leq C \left(\frac{n}{\sqrt{\log{n}}}\right)^{- \frac{8\tau}{4\tau+d}}.
\end{align*}
\item If $\tau = d / 4$, then
\begin{align*}
\sup_{P \in \Par_{\theta}: \atop f_1(\theta) = \tau, f_2(\theta) > \min \left\{ \frac{4 \tau}{4 \tau + d}, \frac{1}{2} \right\}} \E_P \left( \phihat_{n, \khat} - \phi(P) \right)^2 \leq C \left( n \ell(n)^2 \right)^{-1}.
\end{align*}
\item If $\tau > d / 4$ and $\sigma^2 = \E_P \left[ \left( L_1(O) - \E_P L_1(O) \right)^2 \right]$, then for every $P \in \{ P \in \mathcal{P}_\theta: f_1(\theta) = \tau, f_2(\theta) > \min \{ 4 \tau / (4 \tau + d), 1 / 2 \} \}$
\begin{align*}
\sqrt{n} \left( \phihat_{n, \khat} - \phi(P) \right) \rightarrow_d Z \sim \mathcal{N} (0, \sigma^2).
\end{align*}
where $\mathcal{N}(\mu, \sigma^2)$ stands for normal distribution with mean $\mu$ and variance $\sigma^2$.
\end{enumerate}
\end{thm}

A few remarks are in order regarding the implications of Theorem \ref{theorem_general_lepski} as well as the subtleties involved in its proof.

\begin{rem}
Provided one has knowledge of a data generating $\theta$ and therefore of $f_1(\theta)$ and $f_2(\theta)$, one can use the bias and variance properties to achieve an optimal trade-off and subsequently obtain optimal mean squared error in estimating $\phi(P)$ over $P \in \mathcal{P}_{\theta}$ which scales as $n^{- \frac{8 \tau}{4 \tau + d}}$ in the low regularity regime when $0 < \tau < d / 4$ or is $\sqrt{n}$-consistent and semiparametric efficient when $\tau > d / 4$. Theorem \ref{theorem_general_lepski} demonstrates a logarithmic price paid by the estimator $\phihat_{n, \khat}$ in terms of estimating $\phi(P)$ over a class of data generating mechanisms $\left\{ \mathcal{P}_{\theta}: f_1(\theta) = \tau, f_2(\theta) > \frac{4 \tau}{4 \tau + d} \right\}$ -- which will in all examples be non-empty and often infinite. As will be demonstrated in Section \ref{section_examples}, the term $f_1(\theta) = \tau$ usually drives the minimax rate of estimation whereas $f_2(\theta) > \frac{4 \tau}{4 \tau + d}$ is a regularity condition which typically will relate to the marginal density of covariates in observational studies. Moreover, in our examples, the range of $\tau < d / 4$ will necessarily imply the non-existence of $\sqrt{n}$-consistent estimators of $\phi(P)$ over $P \in \Par_{\theta}$ in a minimax sense \citep{robins2009semiparametric}.
\end{rem}
\begin{rem}
The proof of Theorem \ref{theorem_general_lepski} is somewhat involved because of the dependence of the U-statistics kernel on estimated nuisance functions (e.g. the functions $\Ltilde_1, \Ltilde_{2l}$ and $\Ltilde_{2r}$ are estimated from the data). This prohibits performing a Hoeffding decomposition w.r.t. the whole sample and we can only proceed by a conditional Hoeffding decomposition -- allowed by our sample splitting mechanism. Subsequently, the conditioning event needs to be sufficiently well behaved and corresponds to the high probability requirement of the good event defined in Property (A). As we shall see in Section \ref{section_examples}, that the good event corresponds to certain estimates of nuisance functions to be well behaved -- both in terms of sup-norm concentration and membership to desired H\"{o}lder regularity classes. It is for this reason, we need somewhat precise results on adaptive function estimation as well -- as provided by Theorem \ref{thm_adaptive_density_regression} in Section \ref{adaptive_nuisance_function_estimation}.
\end{rem}
\begin{rem}
The sample splitting scheme employed in Theorem \ref{theorem_general_lepski} is not the only option. Many other reasonable options exist, e.g. $n_2 = n / \text{polylog}(n)$. One may also choose $n_1 = n_2 = n / 2$, which only affects the constants when $0 < \tau \leq d / 4$. However, complication arises when $\tau > d / 4$: we need to rely on cross-fitting after sample splitting \citep{chernozhukov2018double} to achieve semiparametric efficiency.
\end{rem}

Finally, it therefore remains to be explored whether this logarithmic price payed in Theorem \ref{theorem_general_lepski} is indeed necessary for the regime $\tau < d / 4$. Using a chi-square divergence inequality developed in the next section along with a suitable version of constrained risk inequality (see Section \ref{section_appendix_technical_lemmas}) we shall argue that the logarithmic price of Theorem \ref{theorem_general_lepski} is indeed necessary for a class of examples naturally arising in many observational studies. 

\subsection{\textbf{Lower Bound}}\label{section:lower_bound}\hspace*{\fill}\\
In terms of lower bounds, we shall use the constrained risk inequality of \cite{cai2011testing} (see Lemma \ref{lemma:prop_cai} in Section \ref{section_appendix_technical_lemmas}) -- which in turn requires the control over the chi-square divergence between the mixture of suitable product measures. In this section we provide such a bound which will be used in the examples in Section \ref{section_examples}.

We closely follow the strategy in \cite{birge1995estimation, robins2009semiparametric} who already provide useful control for the Hellinger distance between mixtures of suitable product measures. It turns out that under certain boundedness assumptions on the Radon-Nykodym derivatives of the measures involved in the problem, one can convert the available bounds on the Hellinger distance quite easily to the chi-squared divergence related bounds. More precisely, as in \cite{robins2009semiparametric}, let $O_1, \ldots, O_n$ be a random sample from a density $p$ w.r.t. measure $\mu$ on a sample space $(\chi,\mathcal{A})$. For $k \in \mathbb{N}$, let $\chi = \cup_{j = 1}^k \chi_j$ be a measurable partition of the sample space. Given a vector $\lambda=(\lambda_1,\ldots,\lambda_k)$ in some product measurable space $\Lambda= \Lambda_1\times\cdots\times \Lambda_k$ let $P_{\lambda}$ and $Q_{\lambda}$ be probability measures on $\chi$ such that
\begin{itemize}
	\item $P_{\lambda}(\chi_j)=Q_{\lambda}(\chi_j)=p_j$ for every $\lambda$ and some $(p_1,\ldots,p_k)$ in the $k$-dimensional simplex.
	\item The restrictions $P_{\lambda}$ and $Q_{\lambda}$ to $\chi_j$ depends on the $j^{\mathrm{th}}$ coordinate $\lambda_j$ of $\lambda=(\lambda_1,\ldots,\lambda_k)$ only.
\end{itemize}
For $p_{\lambda}$ and $q_{\lambda}$ densities of the measures $P_{\lambda}$ and $Q_{\lambda}$ respectively that are jointly measurable in the parameters $\lambda$ and the observations, and $\pi$ a probability measure on $\Lambda$, define $p=\int p_{\lambda}d{\pi(\lambda)}$ and $q =\int q_{\lambda} d(\pi(\lambda))$, and set
$$a=\max_j\sup_{\lambda}\int_{\chi_j}\frac{(p_{\lambda}-p)^2}{p_{\lambda}}\frac{d\mu}{p_j},$$  
$$b=\max_j\sup_{\lambda}\int_{\chi_j}\frac{(q_{\lambda}-p_{\lambda})^2}{q_{\lambda}}\frac{d\mu}{p_j},$$ 
$$\tilde{c}=\max_j \sup_{\lambda}\int_{\chi_j}\frac{p^2}{p_{\lambda}}\frac{d\mu}{p_j},$$
$$\d=\max_j\sup_{\lambda}\int_{\chi_j}\frac{(q-p)^2}{p_{\lambda}}\frac{d\mu}{p_j}.$$

With the notations and definitions as above we now have the following proposition which is a direct application of \cite{robins2009semiparametric}[Theorem 2.1]. 
\begin{prop}\label{prop_chisquare_affinity}
	Suppose that $n p_j(1 \vee a \vee b\vee \tilde{c}) \leq A$ for all $j$ and for all $\lambda$, $\underline{B}\leq p_{\lambda}\leq \overline{B}$ for positive constants $A,\underline{B},\overline{B}$. Then there exists a $C>0$ that depends only on $A,\underline{B},\overline{B}$, such that, for any product probability measure $\pi=\pi_1 \otimes\cdots\otimes\pi_k$, one has
	$$\chi^2 \left( \int P_{\lambda} d (\pi(\lambda)), \int Q_{\lambda} d (\pi(\lambda)) \right) \leq e^{Cn^2(\max_j p_j) (b^2 + a b)+ C n \d} - 1,$$
	where $\chi^2(\nu_1,\nu_2)=\int \left(\frac{d\nu_2}{d\nu_1}-1\right)^2 d\nu_1$ is the chi-square divergence between two probability measures $\nu_1$ and $\nu_2$ with $\nu_2 \ll \nu_1$.
\end{prop}
\begin{proof}
Based on the above construction, \cite{robins2009semiparametric}[Theorem 2.1] showed that, there exists some absolute constant $C'$ such that
$$
H^2 \left( \int P_{\lambda}d(\pi(\lambda)), \int Q_{\lambda}d(\pi(\lambda)) \right) \le C' n^2 (\max_j p_j) (b^2 + ab) + C' n \delta.
$$
where $H(\nu_1, \nu_2)$ denotes the Hellinger distance between two probability measures $\nu_1$ and $\nu_2$. As $p_\lambda \ge \underline{B}$, so is $\int P_{\lambda} d ( \pi (\lambda) ) \ge \underline{B}$. By setting $C = 4 \underline{B}^{-1} C'$, we have
\begin{align}
& \chi^2 \left( \int P_{\lambda} d (\pi(\lambda)), \int Q_{\lambda} d (\pi(\lambda)) \right) \nonumber\\
& \le \exp \left\{ 4 \underline{B}^{-1} H^2 \left( \int P_{\lambda}d(\pi(\lambda)), \int Q_{\lambda}d(\pi(\lambda)) \right) \right\} - 1 \label{EQ:EQUIV} \\
& \le \exp \left\{ 4 \underline{B}^{-1} (C' n^2 (\max_j p_j) (b^2 + ab) + C' n \delta) \right\} - 1 \nonumber \\
& = e^{ C n^2 (\max_j p_j) (b^2 + ab) + C n \delta } - 1 \nonumber
\end{align}
where the first inequality \eqref{EQ:EQUIV} is proved in Appendix \ref{app:equiv} for the sake of completeness. 
\end{proof}

\subsection{\textbf{Adaptive Estimation of Nuisance Functions}}\label{adaptive_nuisance_function_estimation}\hspace*{\fill}\\
As will be evident from our examples in Section \ref{section_examples} that the applications of Theorem \ref{theorem_general_lepski} require certain estimates of nuisance functions to be well behaved -- both in terms of sup-norm concentration and membership to desired H\"{o}lder regularity classes. In this section we provide such estimators of regression and density functions using Lepski-type arguments \citep{lepskii1991problem,lepskii1992asymptotically}. Although the construction and adaptation proof of such estimators are quite standard, we will need a few additional results on sufficiently high probability inclusion of the adaptive estimators in suitable H\"{o}lder regularity classes. We provide the necessary ingredients below.   

Consider i.i.d. data on $O_i=(W_i, \bX_i) \sim P$ for a scalar variable $W$ such that $|W| \leq B_U$ and $\E_P \left( W | \bX \right) = f(\bX)$ almost surely, and $\bX \in [0, 1]^d$ has density $g$ such that $0 < B_L \leq g(\bx) \leq B_U < \infty$ for all $\bx \in [0, 1]^d$. Although to be precise, we should put subscripts $P$ to $f, g$, we omit this since the context of their use is clear. We assume H\"{o}lder type smoothness (defined in Section \ref{section_wavelets and function spaces}) on $f, g$ and let $\mathcal{P}(s, \gamma)=\{P: (f, g) \in H(s, M) \times H(\g, M), | f(\bx) | \leq B_U, B_L \leq g(\bx) \leq B_U \ \forall \ \bx \in [0, 1]^d\}$ denote classes of data generating mechanisms indexed by the smoothness indices. Then we have the following theorem (the proof of which can be found in Appendix \ref{section:proof_of_remaining_theorems}), which considers adaptive estimation of $f, g$ in $L_{\infty}$ norm over $(s, \gamma) \in (s_{\min}, s_{\max}) \times (\gmin, \gmax)$, for given positive real numbers $s_{\min}, s_{\max}, \gmin, \gmax$. 
 
\begin{thm}\label{thm_adaptive_density_regression}
	If $\g_{\min} > s_{\max}$, then there exist $\fhat$ and $\ghat$ depending on $M$, $B_L$, $B_U$, $s_{\min}$, $s_{\max}$, $\gmin$, $\gmax$, and choice of wavelet bases $\bpsi_{0, 0}^0, \bpsi_{0, 0}^1$ (defined in Section \ref{section_wavelets and function spaces}) such that the following hold for every $(s, \gamma) \in (s_{\min}, s_{\max}) \times (\gmin, \gmax)$ with a large enough $C > 0$ depending possibly on  $M, B_L, B_U$, and $\gmax$.
	\begin{equation}\label{eq:moment_g}
	\sup\limits_{P \in \mathcal{P}(s, \gamma)} \E_P \| \ghat - g \|_{\infty} \leq C^{\frac{d}{2 \gamma + d}} \left(\frac{n}{\log{n}}\right)^{-\frac{\g}{2 \g + d}},
	\end{equation}
	\begin{equation}\label{eq:tail_g}
	\sup\limits_{P \in \mathcal{P}(s, \gamma)} \P_P \left( \Vert \ghat - g \Vert_\infty \geq C^{\frac{d}{d + 2 \gamma}} \left( \frac{n}{\log{n}} \right)^{-\frac{\gamma}{2 \gamma + d}} \right)\leq \frac{\log{n}}{n^3},
	\end{equation}
	\begin{equation}\label{eq:holder_g}
	\sup\limits_{P \in \mathcal{P}(s, \gamma)} \P_P \left( \ghat\notin H(\g, C) \right) \leq \frac{1}{n^2},
	\end{equation}
	\begin{equation}\label{eq:moment_f}
	\sup\limits_{P \in \mathcal{P}(s, \gamma)} \E_P \| \fhat - f \|_{\infty} \leq C^{\frac{d}{2 s + d}} \left(\frac{n}{\log{n}}\right)^{-\frac{s}{2 s + d}},
	\end{equation}
	\begin{equation}\label{eq:tail_f}
	\sup\limits_{P \in \mathcal{P}(s, \gamma)} \P_P \left( \Vert \fhat - f \Vert_\infty \geq C^{\frac{d}{d + 2 s}} \left( \frac{n}{\log{n}} \right)^{-\frac{s}{2 s + d}} \right)\leq \frac{\log{n}}{n^3},
	\end{equation}
	\begin{equation}\label{eq:holder_f}
	\sup\limits_{P \in \mathcal{P}(s, \gamma)} \P_P \left( \fhat\notin H(s, C) \right) \leq \frac{1}{n^2},
	\end{equation}
	\begin{equation}\label{eq:bound}
	|\fhat(\bx)| \leq 2 B_U \ \text{and}\ B_L / 2 \leq \ghat(\bx) \leq 2B_U \,\,\,\, \forall \ \bx \in [0,1]^d.
	\end{equation}
\end{thm}
\begin{rem} 
A close look at the proof of Theorem \ref{thm_adaptive_density_regression} shows that the proof continues to hold for $s_{\min} = 0$. Moreover, although we did not keep track of our constants, the purpose of keeping them in the form above is to show that the multiplicative constants are not arbitrarily big when $\beta$ is large.
\end{rem}

Theorems of the flavor of Theorem \ref{thm_adaptive_density_regression} are not uncommon in literature (see \cite[Chapter 8]{gine2016mathematical} for details) and the proof is somewhat standard for most parts of the result. However, proving \eqref{eq:holder_g} and \eqref{eq:holder_f} for data driven adaptive estimators requires somewhat more care than standard results in the literature. In particular, results of the kind stating that $\ghat \in H(\gamma,C)$ with high probability uniformly over $\Par(\beta,\gamma)$ for a suitably large constant $C$ are often very easy to demonstrate. However, our proof shows that a suitably bounded estimator $\ghat$, which adapts over smoothness and satisfies $\ghat \in H(\gamma,C)$ with probability larger than $1-\frac{1}{n^{\kappa}}$ uniformly over $\Par(\beta,\gamma)$, for any $\kappa>0$ and correspondingly large enough $C$. This result in turn turns out to be crucial for the purpose of controlling suitable bias terms in our functional estimation problems -- as specified by the good event property in Property (A) in Theorem \ref{theorem_general_lepski}. Additionally, the results concerning $\fhat$ are relatively less common in an unknown design density setting. Indeed, adaptive estimation of regression function with random design over Besov type smoothness classes has been obtained by model selection type techniques by \cite{baraud2002model} for the case of Gaussian errors. Our results in contrast hold for any regression model with bounded outcomes and compactly supported covariates having suitable marginal design density. 

\begin{rem}
\label{rem:dnn}
Theorem \ref{thm_adaptive_density_regression} estimates the nuisance functions in our problem using standard Cohen-Daubechies-Vial type of wavelet kernel projections. It is an open problem if estimators based on other more general machine learning algorithms such as random forests or deep neural networks (DNNs) still satisfy Property (A). Indeed, recent results in DNNs using ReLU activation functions \citep{schmidt2020nonparametric, farrell2021deep, chen2019nonparametric, chen2019efficient} are non-adaptive to the underlying regularity of the function classes and little is known about the Fourier coefficients of their outputs -- a property which we will crucially require in the good events defined in Property (A) accompanying Theorem \ref{theorem_general_lepski}. Some recent results \citep{rahaman2019spectral, xu2020frequency, luo2020fourier} connecting the training dynamics \citep{luo2021phase} to the behavior of Fourier coefficients in DNNs might be an interesting direction to investigate for future work.
\end{rem}

\section{Examples}\label{section_examples}
In this section we discuss applications of Theorem \ref{theorem_general_lepski} and Proposition \ref{prop_chisquare_affinity} in producing rate optimal adaptive estimators of certain nonparametric functionals commonly arising in statistical and causal inference literature. These include (i) treatment effect type functional, (ii) mean functional in missing data studies, and (iii) quadratic and variance functional in nonparametric regression. Before proceeding we note that the results of this paper can be extended to include the whole class of doubly robust functionals considered in \cite{robins2008higher}. However we only provide specific examples here to demonstrate the clear necessity to pay a sharp poly-logarithmic penalty for adaptation in low regularity regimes. The proof of the results in this section can be found in Appendix \ref{section:proof_of_remaining_theorems}. 
\subsection{\textbf{Weighted Average Treatment Effect}}\label{sec:ate}
In this subsection, we consider estimating the ``treatment effect" of a treatment on an outcome in presence of multi-dimensional confounding variables \citep{crump2009dealing, robins1992estimating}. To be more specific, we consider a binary treatment $A$ and response $Y$ and $d$-dimensional covariate vector $\bX$, and let $\tau$ be the variance weighted average treatment effect defined as
\be
	\tau \coloneqq \E\left(\frac{\Var(A | \bX) c(\bX)}{\E(\Var(A | \bX))}\right) = \frac{\E(\Cov(Y, A | \bX))}{\E(\Var(A | \bX))},
\ee
where
\be
	\label{treatmenteffectfunction_nonparametric}
	c(\bx) = \E(Y | A = 1, \bX = \mathbf{x}) - \E(Y | A = 0, \bX = \mathbf{x}).
\ee
Under the assumption of no unmeasured confounding, $c(\bx)$ is often referred to as the average treatment effect among subjects with covariate value $\bX = \bx$. The reason of referring to $\tau$ as the average treatment effect can be further understood by considering the semiparametric regression model that assumes

\be
	\label{treatmenteffectfunction}
	c(\bx) = \phi^{*} \ \text{for all} \ \bx,
\ee 
or specifically the partial linear model
\be
	\E(Y | A, \bX) = \phi^{*} A + b(\bX).
\ee

It is clear that under \eqref{treatmenteffectfunction} $\tau$ equals $\phi^{*}$. Moreover, the inference on $\tau$ is closely related to the estimation of $\E(\Cov(Y,A|\bX))$ \citep{robins2008higher}. Specifically, point and interval estimators for $\tau$ can be constructed from point and interval estimators of $\E(\Cov(Y,A|\bX))$. To be more specific, for any fixed $\tau^{*} \in \mathbb{R}$, one can define $Y^{*}(\tau^{*})=Y-\tau^{*} A$ and consider $$\phi(\tau^{*})=\E((Y^{*}(\tau^{*})-\E(Y^{*}(\tau^{*})|\bX))(A-\E(A|\bX)))=\E(\Cov(Y^{*}(\tau^{*}),A|\bX)).$$
It is easy to check that $\tau$ is the unique solution of $\phi(\tau^{*})=0$. Consequently, if we can construct estimator $\hat{\phi}(\tau^*)$ of $\phi(\tau^{*})$, then $\hat{\tau}$ satisfying $\phi(\hat{\tau}) = 0$ is an estimator of $\tau$ with desired properties. Moreover, $(1-\alpha)$ confidence set for $\tau$ can be constructed as the set of values of $\tau^{*}$ for which $(1-\alpha)$ interval estimator of $\phi(\tau^{*})$ contains the value $0$. Finally, since $\E_P ( \Cov_P (Y, A | \bX) ) = \E_P ( \E_P ( Y | \bX ) \E_P ( A | \bX ) ) - \E_P( A Y )$, and $\E_P ( A Y )$ is estimable at a parametric rate, the crucial part of the problem hinges on the estimation of $\E_P ( \E_P ( Y | \bX ) \E_P ( A | \bX ) )$.

For the rest of this section, we assume that we observe $n$ i.i.d. copies of $O = (Y, A, \bX) \sim P$ and we want to estimate $\phi(P)= \E_P \left(\Cov_P(Y,A|\bX)\right)$. We assume that the marginal distribution of $\bX$ has a density w.r.t. Lebesgue measure on $\R^d$ that has a compact support, which we assume to be $[0,1]^d$ and let $g$ be the marginal density of $\bX$ (i.e. $\E_P (h(\bX)) = \int_{[0,1]^d} h(\bx) g(\bx) d\bx$ for all $P$-integrable function $h$), $a(\bX) \coloneqq \E_P(A|\bX)$, $b(\bX) \coloneqq \E_P(Y|\bX)$, and $c(\bX)=\E_P \left( Y | A = 1, \bX \right) - \E_P \left( Y | A = 0, \bX \right)$. Although to be precise, we should put subscripts $P$ to $a, b, g, c$, we omit this since the context of their use is clear. 

To connect this example with the notation employed in Theorem \ref{theorem_general_lepski}, we proceed as follows. Let $\Theta \coloneqq \{ \theta = (\alpha, \beta, \gamma): \frac{\a + \b}{2d} > 0, \gmax \geq \g > \gamma_{\min}\geq 2 ( 1 + \epsilon ) \max \{ \a, \b \} \}$ for some fixed $\epsilon > 0$, and let $\Par_{\theta}$ denote all data generating mechanisms $P$ satisfying the following conditions for known positive constants $M, B_L, B_U$.
\begin{enumerate}[(1)]
	\item $\max \{ |Y|, |A| \} \leq B_U$ a.s. $P$.
	\item $a \in H(\a, M)$, $b \in H(\b, M)$, and $g \in H(\g, M)$. 	 	
	\item $0 < B_L < g(\bx) < B_U$ for all $\bx \in [0,1]^d$.
\end{enumerate}

Note that we do not put any assumptions on the function $c$. Indeed for $Y$ and $A$ binary random variables, the functions $a, b, g, c$ are variation independent. Following our discussion above, we will discuss adaptive estimation of $\phi(P) = \E_P ( \Cov_P (Y, A | \bX ) ) = \E_P ( ( Y - b(\bX) ) ( A - a (\bX) ) )$ over $P \in \Par$. In particular, we summarize our results on upper and lower bounds on the adaptive minimax estimation rate of $\phi(P)$ in the following theorem. 
 \begin{thm}\label{thm_treatment_effect} Assume $(1)-(3)$ and $(\a, \b, \g) \in \Theta$.  Then the following hold for  positive $C, C'$ depending on $M, B_L, B_U, \gmax$ and any sequence $\ell(n)\rightarrow 0$ such that $\ell(n)\log{n}\rightarrow \infty$.
 \begin{enumerate}[(i)]
 	\item (Upper Bound) There exists an estimator $\phihat$, depending only on  $M$, $B_L$, $B_U$, $\gmax$ such that
 	\begin{align*}
	& \sup_{P\in \Par_{(\alpha,\beta,\gamma)}}\E_P\left(\phihat-\phi(P)\right)^2 \\
	& \leq \left\{ \begin{array}{ll}
	C\left(\frac{n}{\sqrt{\log{n}}}\right)^{- \frac{4\alpha+4\beta}{2\alpha+2\beta+d}} & \text{ if } 0 < (\alpha + \beta) / 2 < d / 4, \\
	C\left( n \ell(n)^2 \right)^{-1} & \text{ if } (\alpha + \beta) / 2 = d / 4,
	\end{array} \right.
	\end{align*}
	and if $(\alpha + \beta) / 2 > d / 4$, with $$\sigma^2 = \E_P \left[ \left\{ ( Y - b(\bX) ) ( A - a (\bX) ) - \phi(P) \right\}^2 \right],$$ for every $P \in \{ P \in \Par_{(\alpha,\beta,\gamma)} \}$, we have
	$$
	\sqrt{n} \left( \phihat_{n, \khat} - \phi(P) \right) \rightarrow_d Z \sim \mathcal{N} ( 0, \sigma^2 ).
	$$
 	\item (Lower Bound) Suppose $\{A,Y\}\in \{0,1\}^2$ and $0 < (\alpha + \beta) / 2 < d / 4$. If one has 
 	$$\sup_{P\in \Par_{(\alpha,\beta,\gamma)}}\E_P\left(\phihat-\phi(P)\right)^2 \leq C\left(\frac{n}{\sqrt{\log{n}}}\right)^{-\frac{4\alpha+4\beta}{2\alpha+2\beta+d}},$$
 	for an estimator $\hat{\phi}$ of $\phi(P)$. Then there exists a class of distributions  $\Par_{(\a^{'},\b^{'},\g^{'})}$ such that 
 	$$\sup_{P^{'}\in \Par_{(\a^{'},\b^{'},\g^{'})}}\E_{P'}\left(\phihat-\phi(P^{'})\right)^2 \geq C'\left(\frac{n}{\sqrt{\log{n}}}\right)^{-\frac{4\alpha'+4\beta'}{2\alpha'+2\beta'+d}}.$$
 \end{enumerate}
 \end{thm}
 
Theorem \ref{thm_treatment_effect} describes the adaptive minimax estimation rate of the treatment effect functional in both low regularity regime ($\frac{\alpha+\beta}{2} \leq \frac{d}{4}$) i.e. when $\sqrt{n}$-rate estimation is not possible and regular $\sqrt{n}$-regimes ($\frac{\alpha+\beta}{2} > \frac{d}{4}$) where our result demonstrates adaptive semi-parametric efficiency. Note that unlike the classical adaptation problem regarding quadratic functional estimation in density and white noise models, the adaptation in this problem is w.r.t. to hyperplanes $\{(\alpha,\beta):(\alpha+\beta)/2d=\tau\}$. Finally, we note that the case of $\frac{\alpha+\beta}{2} = \frac{d}{4}$ also incurs an additional penalty (which can be made to grow arbitrarily slow) over usual $\sqrt{n}$-rate of convergence resulted from adaptation. Although we do not state it explicitly, the proof of the lower bound shows that if the rate of convergence $\frac{\alpha+\beta}{2}=\frac{d}{4}$ is $O(n^{-1})$ then one pays a polynomial penalty for $\frac{\alpha+\beta}{2}<\frac{d}{4}$ -- which is indeed suboptimal. 

It is worth noting that if the set of $(\alpha,\beta)$ \textit{only} includes the case $\frac{\alpha+\beta}{2} > \frac{d}{4}$, one can indeed obtain adaptive and even semi-parametric efficient estimation of the functionals studied here without effectively any assumption on $g$. When $\frac{\alpha+\beta}{2} > \frac{d}{4}$, the dependence of the smoothness $\gamma$ of the density $g(\cdot)$ on $\a$ and $\b$ can also be relaxed by considering estimators using higher-order U-statistics. But we do not further pursue such direction in this paper. The interested reader can find the relevant details in \cite{robins2017minimax,mukherjee2017semiparametric}.

\begin{rem}[Further comments on the dependence of $\gamma_{\min}$ on $\alpha \vee \beta$]
In Theorem \ref{thm_treatment_effect} (and also in Theorem \ref{thm_missing_data} and \ref{thm_quadratic_functional_regression} in the next two subsections), we provide a sufficient condition on the dependence of the known smoothness lower bound $\gmin$ of the density of $\bX$ on $\alpha \vee \beta$ - the maximum between the smoothness of $a(\cdot)$ and $b(\cdot)$ - to obtain matching adaptive minimax upper and lower bound. At a first sight, such dependence looks dimension free. However, when $\alpha + \beta \geq d / 4$, $\gmin$ implicitly grows with $d$. As a matter of fact, one can easily obtain sharper dependence of $\gmin$ on $\alpha$, $\beta$, and $d$ by looking at the proof of Theorem \ref{thm_treatment_effect} in the appendix. We decide to keep the constraint on $\gmin$ in the current form for ease of exposition. As mentioned above, such dependence can even further relaxed by considering estimators based on even higher-order U-statistics. We refer to Section \ref{section_discussion} for more discussions on this.
\end{rem}

\subsection{\textbf{Mean Response in Missing Data Models}} \label{sec:mar}
Suppose we have $n$ i.i.d observations on $O = (Y A,A, \bX)\sim P$, for a response variable $Y\in \mathbb{R}$ which is conditionally independent of the missingness indicator variable $A\in \{0,1\}$ given covariate information $\bX$. In literature, this assumption is typically known as the missing at random (MAR) assumption. Under this assumption, our quantity of interest $\phi(P)=\E_P(Y)$ is identifiable as $\E_P[\E_P[\left(Y|A=1,X\right)]]$ from the observed data. This model is a canonical example of a study with missing response variable and to make this assumption reasonable, the covariates must contain the information on possible dependence between response and missingness. We refer the interested readers to \cite{tsiatis2007semiparametric} for the history of statistical analysis of MAR and related models.

To lay down the mathematical formalism for minimax adaptive estimation of $\phi(P)$ in this model, let $f$ be the marginal density of $\bX$, $a(\bX)\coloneqq\E_P(A|\bX) = \BP_{P} (A = 1 | \bX)$ which is often called ``the propensity score'' in the causal inference literature, and $b(\bX)\coloneqq\E_P(Y|A=1,\bX)=\E_P(Y|\bX)$, and $g(\bX)=f(\bX)/a(\bX)$ (with the convention of the value $+\infty$ when dividing by $0$). Although to be precise, we should put subscripts $P$ to $a,b,g$, we omit this since the context of their use is clear. 

To connect this example with the notation employed in Theorem \ref{theorem_general_lepski}, we proceed as follows. Let $\Theta\coloneqq\{\theta=(\alpha,\beta,\gamma): \frac{\a+\b}{2d} > 0, \gmax\geq \g >\gamma_{\min}\geq 2(1+\epsilon)\max\{\a,\b\} \}$ for some fixed $\epsilon>0$, and let $\Par_{\theta}$ denote all data generating mechanisms $P$ satisfying the following conditions for known positive constants $M,B_L,B_U$.
\begin{enumerate}[(1)]
 	    \item $|Y| \leq B_U$.
 	 	\item $a \in H(\a,M)$, $b \in H(\b,M)$, and $g\in H(\g,M)$. 
 	 	\item $B_L<g(\bx), a(\bx)<B_U$ for all $\bx \in [0,1]^d$.
 	 \end{enumerate}
 We then have the following result.
 	 \begin{thm}\label{thm_missing_data} Assume $(1)-(3)$ and $(\a,\b,\g)\in \Theta$. Then the following hold for  positive $C,C'$ depending on $M,B_L,B_U,\gmax$ and any sequence $\ell(n)\rightarrow 0$ such that $\ell(n)\log{n}\rightarrow \infty$.
 	  \begin{enumerate}[(i)]
 	  	\item (Upper Bound) There exists an estimator $\phihat$, depending only on  $M$,$B_L$,$B_U$, $\gmax$ such that
 	  	\begin{align*}
	& \sup_{P\in \Par_{(\alpha,\beta,\gamma)}}\E_P\left(\phihat-\phi(P)\right)^2 \\
	& \leq \left\{ \begin{array}{ll}
	C\left(\frac{n}{\sqrt{\log{n}}}\right)^{-\frac{4\alpha+4\beta}{2\alpha+2\beta+d}} & \text{ if } 0 < (\alpha + \beta) / 2 < d / 4, \\
	C\left( n \ell(n)^2 \right)^{-1} & \text{ if } (\alpha + \beta) / 2 = d / 4,
	\end{array} \right.
	\end{align*}
	and if $(\alpha + \beta) / 2 > d / 4$, with $$\sigma^2 = \E_P \left[ \left\{ A a(\bX) ( Y - b(\bX) ) + b(\bX) - \phi(P) \right\}^2 \right],$$ for every $P \in \{ P \in \Par_{(\alpha,\beta,\gamma)} \}$, we have
	$$
	\sqrt{n} \left( \phihat_{n, \khat} - \phi(P) \right) \rightarrow_d Z \sim \mathcal{N} ( 0, \sigma^2 ).
	$$
 	  	\item (Lower Bound) Suppose $\{A,Y\}\in \{0,1\}^2$ and $0 < (\alpha + \beta) / 2 < d / 4$. If one has 
 	  	$$\sup_{P\in \Par_{(\alpha,\beta,\gamma)}}\E_P\left(\phihat-\phi(P)\right)^2 \leq C\left(\frac{n}{\sqrt{\log{n}}}\right)^{-\frac{4\alpha+4\beta}{2\alpha+2\beta+d}},$$
 	  	for an estimator $\hat{\phi}$ of $\phi(P)$. Then there exists a class of distributions  $\Par_{(\a^{'},\b^{'},\g^{'})}$ such that 
 	  	$$\sup_{P^{'}\in \Par_{(\a^{'},\b^{'},\g^{'})}}\E_{P'}\left(\phihat-\phi(P^{'})\right)^2 \geq C'\left(\frac{n}{\sqrt{\log{n}}}\right)^{-\frac{4\alpha'+4\beta'}{2\alpha'+2\beta'+d}}.$$
 	  \end{enumerate}
 	  \end{thm}

Once again, Theorem \ref{thm_missing_data} describes the adaptive minimax estimation rate of the treatment effect functional in both low regularity regime ($\frac{\alpha+\beta}{2} < \frac{d}{4}$) i.e. when $\sqrt{n}$-rate estimation is not possible and regular $\sqrt{n}$-regimes ($\frac{\alpha+\beta}{2} > \frac{d}{4}$)  where our result demonstrates adaptive semi-parametric efficiency. The case of $\frac{\alpha+\beta}{2} = \frac{d}{4}$ also incurs an additional penalty (which can be made to grow arbitrarily slow) over usual $\sqrt{n}$-rate of convergence resulted from adaptation. As mentioned in the end of Section \ref{sec:ate}, the dependence on the smoothness $\gamma$ of the function $g(\cdot)$ can be relaxed using higher-order U-statistics.

\subsection{\textbf{Quadratic and Variance Functionals in Regression Models}}
Consider a sample of $n$ i.i.d copies of $O = (Y, \bX) \sim P$ and the functional of interest is the expected value of the square of the regression of $Y$ on $\bX$. Specifically suppose we want to estimate $\phi(P)= \E_P\left(\left\{E_P(Y|\bX)\right\}^2\right)$. 
Assume that distribution of $\bX$ has a density w.r.t. Lebesgue measure on $\R^d$ that has a compact support, which we assume to be $[0, 1]^d$ for sake of simplicity. Let $g$ be the marginal density of $\bX$, and $b(\bX)\coloneqq\E_P(Y|\bX)$. 

To connect this example with the notation employed in Theorem \ref{theorem_general_lepski}, we proceed as follows. Let $\Theta\coloneqq\{\Par(\beta,\gamma): \beta / d > 0, \gmax\geq \g > \gamma_{\min}\geq 2(1+\epsilon)\beta \}$ for some fixed $\epsilon>0$, and by $\Par(\beta,\gamma)$ we consider all data generating mechanisms $P$ satisfying the following conditions for known positive constants $M,B_L,B_U$.
 \begin{enumerate}[(1)]
 	    \item $\max\{|Y|\}\leq B_U$.
 	 	\item $b \in H(\b,M)$, and $g\in H(\g,M)$. 
 	 	\item $0<B_L<g(\bx)<B_U$ for all $\bx \in [0,1]^d$.
 	 \end{enumerate}
\begin{thm}\label{thm_quadratic_functional_regression} Assume $(1)-(3)$ and $(\b,\g)\in \Theta$.  Then the following hold for  positive $C,C'$ depending on $M,B_L,B_U,\gmax$ and any sequence $\ell(n)\rightarrow 0$ such that $\ell(n)\log{n}\rightarrow \infty$.
 \begin{enumerate}[(i)]
 	\item (Upper Bound) There exists an estimator $\phihat$, depending only on  $M$,$B_L$,$B_U$, $\gmax$ such that
 	\begin{align*}
	\sup_{P \in \Par_{(\beta,\gamma)}}\E_P\left(\phihat-\phi(P)\right)^2 \leq \left\{ \begin{array}{ll}
	C\left(\frac{n}{\sqrt{\log{n}}}\right)^{-\frac{8\beta}{4\beta+d}} & \text{ if } 0 < \beta < d / 4, \\
	C\left( n \ell(n)^2 \right)^{-1} & \text{ if } \beta = d / 4,
	\end{array} \right.
	\end{align*}
	and if $\beta > d / 4$, with $\sigma^2 = \E_P \left[ \left( (2 Y - b(\bX)) b(\bX) - \phi(P) \right)^2 \right]$, for every $P \in \{ P \in \Par_{(\beta,\gamma)} \}$, we have
	$$
	\sqrt{n} \left( \phihat_{n, \khat} - \phi(P) \right) \rightarrow_d Z \sim \mathcal{N} ( 0, \sigma^2 ).
	$$
 	\item (Lower Bound) Suppose $Y\in \{0,1\}^2$ and $0 < \beta < d / 4$. If one has 
 	$$\sup_{P\in \Par(\beta,\gamma)}\E_P\left(\phihat-\phi(P)\right)^2 \leq C\left(\frac{n}{\sqrt{\log{n}}}\right)^{-\frac{8\beta}{4\beta+d}},$$
 	for an estimator $\hat{\phi}$ of $\phi(P)$. Then there exists a class of distributions  $\Par(\b^{'},\g^{'})$ such that 
 	$$\sup_{P^{'}\in \Par(\b^{'},\g^{'})}\E_{P}\left(\phihat-\phi(P^{'})\right)^2 \geq C'\left(\frac{n}{\sqrt{\log{n}}}\right)^{-\frac{8\beta'}{4\beta'+d}}.$$
 \end{enumerate}
 \end{thm}
\begin{rem}
Although Theorem \ref{thm_quadratic_functional_regression} and the discussion before that are made in the context of estimating a particular quadratic functional in the context of a regression framework, it is worth noting that the result also extends to estimating classical quadratic functionals in density models \citep{efromovich1996optimal,gine2008simple}.
\end{rem}
One can also consider in the same setup, the estimation of functionals related to the conditional variance of $Y$ under such a regression model, which has been studied in detail by \cite{hall1989variance, ruppert1997local, fan1998efficient, brown2007variance,cai2008adaptive}. Whereas, the minimax optimal and adaptive results in \cite{brown2007variance,cai2008adaptive} are in a equi-spaced fixed design setting, one can use an analogue of Theorem \ref{thm_quadratic_functional_regression} to demonstrate a rate adaptive estimator and corresponding matching lower bound, with a mean-squared error of the order of $\left(\frac{n}{\sqrt{\log{n}}}\right)^{-\frac{8\beta}{4\beta+d}}$ for estimating $\E_P ( \Var_P ( Y | \bX ) )$ adaptively over H\"{o}lder balls of regularity $\beta < \frac{d}{4}$. As noted by \cite{robins2008higher}, this rate is higher than the rate of estimating the conditional variance in mean-squared error for the equi-spaced fixed design \citep{cai2008adaptive}. In a similar vein, one can also obtain this type of results for the estimation of conditional variance under the assumption of homoscedasticity i.e. $\sigma^2 \coloneqq \Var \left( Y | \bX = \bx \right)$ for all $\bx \in [0,1]^d$. In particular, there exists an estimator for $\sigma^2$ with mean-squared error of order $\left(\frac{n}{\sqrt{\log{n}}}\right)^{-\frac{8\beta}{4\beta+d}}$ when $\beta<\frac{d}{4}$ and of order $n^{-1}$ when $\beta>\frac{d}{4}$. As demonstrated recently by \cite{shen2020optimal}, the (non-adaptive) minimax lower bound for estimating homoscedastic variance $\sigma^2$ in any dimension is indeed $n^{-\frac{8\beta}{4\beta+d}}$ when $\beta<\frac{d}{4}$ and our proof results in an adaptive estimator of the same rate in general dimensions (which also guarantees efficient estimation whenever $\beta>\frac{d}{4}$) while assuming sufficient smoothness on the marginal density of $\bX$. For instance, a candidate sequence of $\phihat_{n, k}$'s for this purpose can be constructed by taking $A=Y$ in the treatment effect example considered in Section \ref{sec:ate}. The requirement of not needing a smoothness on the marginal density of $\bX$ can be removed for $d=1$ (see e.g. \cite{robins2008higher,shen2020optimal}) but doing so in higher dimensions remains a challenge. Although we do not pursue it here, it is possible to follow the lines of argument in \cite{shen2020optimal} to show the requirement of a logarithmic penalty while trying to adapt in the region $\beta<\frac{d}{4}$.




\section{Discussion}\label{section_discussion}
In this paper, we have extended the results for adaptive estimation of non-linear integral functionals from density estimation and Gaussian white noise models, to move towards adaptive estimation of non-linear functionals in more complex nonparametric models. Below we make a few comments on some of the assumptions and future research directions.
\begin{enumerate}\itemsep3pt

    \item [(i)] The adaptation considered in this paper is w.r.t. the smoothness of certain nuisance functions that naturally arise in a class of nonparametric problems. For the three concrete examples considered above, the smoothness of the regression functions decides the minimax rates of the problem and consequently our results naturally adapt to this regularity. However, since our estimators are based on second-order U-statistics, we can only perform a second-order bias correction (see e.g. \cite{robins2008higher,robins2017minimax}) and this results in the requirement of sufficient large smoothness of the marginal density of the covariates. Using higher order U-statistics this dependence can be substantially lowered while considering non-adaptive minimax rates \citep{robins2008higher,robins2017minimax,mukherjee2017semiparametric}. However, such higher order U-statistics make the adaptation proof substantially more difficult and it is currently beyond the scope of this paper. Moreover, it also remains open to decide the minimum smoothness required from the marginal density of the covariates for even non-adaptive minimax risks to go through in general (fixed) dimensions $d>1$ (for $d=1$ \cite{shen2020optimal} shows that one can do without the requirement on density smoothness in one of the special cases considered here). 
    
    \item [(ii)] The examples considered here involve bounded outcomes. This requirement is mostly used while getting suitable tail bounds for the U-statistics under study. Although we do not pursue it here, we believe that this condition can be weakened by assuming sub-Gaussian or other related light tails (by carefully using results in \cite{houdre2003exponential,gine2000exponential,chakrabortty2018tail}) for the conditional error distribution of the outcome given the covariates. 

    \item [(iii)] The sequence of estimators constructed here depends both on the estimation of multivariate regressions and density functions. The requirement of estimating a multivariate density as well as the requirement of smoothness on the marginal density of the covariates can be entirely removed when making the assumptions such that $\sqrt{n}$-consistent estimation is possible. More precisely, in all of the three examples above, there exists a sequence of estimators which does not require estimation or any smoothness requirement on the marginal density of the covariates and automatically adapts to the $\sqrt{n}$-consistent and semiparametric efficient estimators in the $\sqrt{n}$-estimable range. This sequence of estimators was recently constructed in \cite{mukherjee2017semiparametric} -- but they do not generalize to the non-$\sqrt{n}$-estimable part of the problems.
    
    \item [(iv)] Adaptive estimation of the functionals is only a first step of inference and provides no guidelines for statistical inference. It remains an extremely interesting question to explore honest adaptive confidence set constructions for the class of nonparametric functionals considered above.
    
    \item [(v)] Finally, as mentioned in Remark \ref{rem:dnn}, it is interesting to investigate if the adaptation theory for nonparametric functional estimation developed in this paper can be generalized when the nuisance functions (including the design densities) are estimated by deep neural networks instead of classical nonparametric regression such as wavelets regression.
    
\end{enumerate}

\section{Wavelets, Projections, and H\"{o}lder Spaces}\label{section_wavelets and function spaces}\hspace*{\fill}\\
We work with certain Besov-H\"{o}lder type spaces which we define in terms of moduli of wavelet coefficients of continuous functions. 
For $d>1$, consider expansions of functions $h \in L_2\left([0,1]^d\right)$ on an orthonormal basis of compactly supported bounded wavelets of the form 
\be
h(\bx)&=\sum_{m \in \mathbb{Z}^d}\langle h, \bpsi_{0, m}^0\rangle \bpsi_{0, m}^0(\bx)+ \sum_{l=0}^{\infty}\sum_{m \in \mathbb{Z}^d}\sum\limits_{v \in \{0, 1\}^d \setminus \{0\}^d}\langle h, \bpsi_{l, m}^v\rangle \bpsi_{l, m}^v(\bx) , 
\ee
where the base functions $\bpsi_{l, m}^v$ are orthogonal for different indices $(l, m, v)$ and are scaled and translated versions of the $2^d$ $S$-regular base functions $\bpsi_{0, 0}^v$ with $S > \beta$, i.e., $\bpsi_{l, m}^v(\bx)=2^{ld/2}\bpsi_{0,0}^v(2^l \bx-m)=\prod_{j=1}^{d}2^{\frac{l}{2}}\psi_{0,0}^{v_j}\left(2^lx_j-m_j\right)$ for $m=(m_1,\ldots,m_d) \in \ZZ^d$ and $v=(v_1,\ldots,v_d)\in \{0,1\}^d$ with $\psi_{0,0}^0=\phi$ and $\psi_{0,0}^1=\psi$ being the scaling function and mother wavelet of regularity $S$ respectively as defined in one dimensional case. As our choices of wavelets, we will throughout use compactly supported scaling and wavelet functions of (boundary-corrected) Cohen-Daubechies-Vial type with $S$ first null moments \citep{cohen1993wavelets}. In view of the compact support of the wavelets, for each resolution level $l$ and index $v$, only $O(2^{ld})$ base elements $\psi_{l,m}^v$ are non-zero on $[0,1]^{d}$; let us denote the corresponding set of indices $m$ by $\Z_l$ and obtain the representation,
\be
h(\bx)&=\sum_{m \in \Z_{J_0}}\langle h, \bpsi_{J_0, m}^0\rangle \bpsi_{J_0, m}^0(\bx)+ \sum_{l=J_0}^{\infty}\sum_{m \in \Z_l}\sum\limits_{v \in \{0, 1\}^d \setminus \{0\}^d}\langle h, \bpsi_{l, m}^v\rangle \bpsi_{l, m}^v(\bx),\label{eqn:wavelet_expansion_aad_compact_j0}
\ee
where $J_0=J_0(S)\geq 1$ is such that $2^{J_0}\geq S$ \citep{cohen1993wavelets, gine2016mathematical}. Thereafter, let for any $h \in L_2 ([0,1]^d)$, $\|\langle h,\bpsi_{\lprime,\cdot}\rangle\|_2$ be the vector $L_2$ norm of the vector $$\left(\langle h,\bpsi^v_{\lprime, m'}\rangle: m' \in \mathcal{Z}_{\lprime},v \in \left\{0,1\right\}^d\right).$$

We will be working with projections onto subspaces defined by truncating expansions as above at certain resolution levels. For example letting
\be
V_j\coloneqq\mathrm{span}\left\{\bpsilk^v, J_0 \leq l\leq j, m \in \Z_l, v\in \{0,1\}^d\right\}, j \geq J_0, \label{eqn:defn_vj}
\ee
one immediately has the following orthogonal projection kernel onto $V_j$ as 
\be
K_{V_{j}}\left(\bx_1,\bx_2\right)=\sum_{m \in \Z_{J_0}}\bpsi_{J_0, m}^0(\bx_1)\bpsi_{J_0, m}^0(\bx_2)+\sum_{l=J_0}^{j}\sum_{m \in \Z_l}\sum_{v \in \{0,1\}^d \setminus \{0\}^d}\bpsi_{l, m}^v(\bx_1)\bpsi_{l, m}^v(\bx_2). \label{eqn:projection_vj}
\ee
As mentioned in Section \ref{section_main_upper_bound}, we sometimes also denote $V_j$ by $V_k$ for $k = 2^{j d}$ . Owing to the MRA property of the wavelet bases, it is easy to see that $K_{V_j}$ has the equivalent representation as 
\be
K_{V_{j}}\left(\bx_1,\bx_2\right)=\sum_{m \in \Z_j}\sum_{v \in \{0,1\}^d}\psi_{j, m}^v\left(\bx_1\right)\psi_{j, m}^v\left(\bx_2\right). \label{eqn:projection_vj_alt}
\ee

Thereafter, using $S$-regular scaling and wavelet functions of Cohen-Daubechies-Vial type with $S>\beta$ let
\be
H(\beta,M) \coloneqq\left\{\begin{array}{c}h \in C\left([0,1]^d\right): \\  2^{J_0(\beta+\frac{d}{2})}\|\langle h,\bpsi^{0}_{J_0\cdot}\rangle\|_{\infty}+\sup\limits_{l \geq 0,m \in \mathbb{Z}^d,v \in \{0,1\}^d \setminus \{0\}^d}2^{l(\beta+\frac{d}{2})}|\langle h, \bpsi_{l,m}^v\rangle|\leq M \end{array} \right\}, \label{eqn:besov_multidim}
\ee
with $C\left([0,1]^d\right)$ being the set of all continuous bounded functions on $[0,1]^d$. It is standard result in the theory of wavelets that $H(\beta,M)$ is related to classical H\"{o}lder-Zygmund spaces with equivalent norms (see \citet[Chapter 4]{gine2016mathematical} for details). For $0<\beta<1$, $H(\beta,M)$ consists of all functions in $C\left([0,1]^d\right)$ such that $\|f\|_{\infty}+\sup\limits_{\bx_1,\bx_2\in [0,1]^d}\frac{|f(\bx_1)-f(\bx_2)|}{\|\bx_1-\bx_2\|^{\beta}}\leq C(M)$. For non-integer $\beta>1$, $H(\beta,M)$ consists of all functions in $C\left([0,1]^d\right)$ such that $f^{(\lfloor \beta \rfloor)}\in C\left([0,1]^d\right)$ for any partial $f^{(\lfloor \beta \rfloor)}$ of order $\lfloor \beta \rfloor$ of $f$ and $\|f\|_{\infty}+\sup\limits_{\bx_1,\bx_2\in [0,1]^d}\frac{|f^{(\lfloor \beta \rfloor)}(\bx_1)-f^{(\lfloor \beta \rfloor)}(\bx_2)|}{\|\bx_1-\bx_2\|^{\beta-\lfloor \beta \rfloor}}\leq C(M)$. Therefore, the functions in $H(\beta,M)$ are automatically uniformly bounded by a number depending on the radius $M$. \\

\begin{rem}
\label{rem:numeric}
Since the Cohen-Daubechies-Vial type wavelets in general do not have closed-form expressions, we briefly comment on the numerical issues when in practice we are only able to evaluate wavelets to a certain accuracy. In our paper, it suffices to approximate the scaling and wavelet functions $\phi$ and $\psi$ by $\tilde{\phi}$ and $\tilde{\psi}$ to the accuracy such that: for any $h \in L_{2} ([0, 1]^{d})$
\begin{align*}
\Vert h_{j_{\max}} - \tilde{h}_{j_{\max}} \Vert_{\infty} \ll n^{- 1 / 2}
\end{align*}
because the finest resolution $j_{\max}$ is chosen so that $2^{j_{\max} d} \asymp n^{2}$. Here $h_{j_{\max}}$ is the projection of $h$ onto $V_{j_{\max}}$ and $\tilde{h}_{j_{\max}}$ is the projection of $h$ onto $\tilde{V}_{j_{\max}}$, where
\begin{align*}
\tilde{V}_{j}\coloneqq\mathrm{span}\left\{\tilde{\bpsi}_{lk}^v, J_0 \leq l\leq j, m \in \Z_l, v\in \{0,1\}^d\right\}, j \geq J_0,
\end{align*}
and $\tilde{\bpsi}_{lk}^v$ is defined in the same way as $\bpsi_{lk}^v$ but with $\phi$ and $\psi$ replaced by $\tilde{\phi}$ and $\tilde{\psi}$. As a result, it suffices to have the following:
\begin{align*}
& \ \left\Vert \begin{array}{c}
\sum_{m \in \Z_{J_0}}\langle h, \bpsi_{J_0, m}^0\rangle \bpsi_{J_0, m}^0(\bx)+ \sum_{l=J_0}^{j_{\max}}\sum_{m \in \Z_l}\sum\limits_{v \in \{0, 1\}^d \setminus \{0\}^d}\langle h, \bpsi_{l, m}^v\rangle \bpsi_{l, m}^v(\bx) \\
- \sum_{m \in \Z_{J_0}}\langle h, \tilde{\bpsi}_{J_0, m}^0\rangle \tilde{\bpsi}_{J_0, m}^0(\bx) - \sum_{l=J_0}^{j_{\max}}\sum_{m \in \Z_l}\sum\limits_{v \in \{0, 1\}^d \setminus \{0\}^d}\langle h, \tilde{\bpsi}_{l, m}^v\rangle \tilde{\bpsi}_{l, m}^v(\bx)
\end{array} \right\Vert_{\infty} \ll n^{-1 / 2} \\
\Leftarrow & \ \underbrace{\vert \Z_{J_0} \vert \max_{m \in \Z_{J_0}} \Vert \langle h, \bpsi_{J_0, m}^0\rangle \bpsi_{J_0, m}^0(\bx) - \langle h, \tilde{\bpsi}_{J_0, m}^0\rangle \tilde{\bpsi}_{J_0, m}^0(\bx) \Vert_{\infty}}_{(A)} \ll n^{- 1 / 2} \\
& \text{and } \underbrace{j_{\max} \vert \Z_{j_{\max}} \vert \max_{l = J_{0}, \cdots, j_{\max}} \max_{m \in \Z_l} \max_{v \in \{0, 1\}^d \setminus \{0\}^d} \Vert \langle h, \bpsi_{l, m}^v\rangle \bpsi_{l, m}^v(\bx) - \langle h, \tilde{\bpsi}_{l, m}^v\rangle \tilde{\bpsi}_{l, m}^v(\bx) \Vert_{\infty}}_{(B)} \ll n^{- 1 / 2}.
\end{align*}

Now our goal is to find the appropriate accuracy $\delta$ with $\Vert \tilde{\phi} - \phi \Vert_{\infty} + \Vert \tilde{\psi} - \psi \Vert_{\infty} \leq \delta$ such that $(A) \ll n^{- 1 / 2}$ and $(B) \ll n^{- 1 / 2}$ hold. $\Vert \tilde{\phi} - \phi \Vert_{\infty} + \Vert \tilde{\psi} - \psi \Vert_{\infty} \leq \delta$ leads to the following:
\begin{equation} \label{single}
\Vert \bpsi_{l, m}^v - \tilde{\bpsi}_{l, m}^v \Vert_{\infty} \lesssim 2^{l d / 2} d \delta,
\end{equation}
where we use the trivial identity 
\begin{equation} \label{trivial}
\prod_{i = 1}^{d} a_{i} - \prod_{i = 1}^{d} \tilde{a}_{i} \equiv \sum_{i = 1}^{d} \prod_{j = 1}^{i} a_{j} (a_{i} - \tilde{a}_{i}) \prod_{j' = i + 1}^{d} \tilde{a}_{j'}
\end{equation} 
for any $a_{1}, \cdots, a_{d}; \tilde{a}_{1}, \cdots, \tilde{a}_{d}$.

With such an error bound, for $(A) \ll n^{- 1 / 2}$ and $(B) \ll n^{- 1 / 2}$ to hold, it suffices to have the following, where we again use the identity \eqref{trivial} together with \eqref{single}:
\begin{align*}
& \ j_{\max} \vert \Z_{j_{\max}} \vert 2^{j_{\max} d / 2} d \delta 2^{j_{\max} d / 2} \ll n^{- 1 / 2} \\
\Rightarrow & \ \log (n^{2}) 2^{j_{\max} d} 2^{j_{\max} d} \delta \ll n^{- 1 / 2} \\
\Rightarrow & \ \log (n^{2}) n^{4} \delta \ll n^{- 1 / 2} \\
\Rightarrow & \ \delta \ll n^{- 9 / 2} / \log (n)
\end{align*}
where the second $2^{j_{\max} d / 2}$ in the first line is due to the product of two $\bpsi_{j_{\max}, m}$'s (and the product of two $\tilde{\bpsi}_{j_{\max}, m}$'s) within the $\Vert \cdot \Vert_{\infty}$ in term (B), in the second line we use $\vert \Z_{j_{\max}} \vert = O (2^{j_{\max} d})$ and in the third line we use $2^{j_{\max} d} \asymp n^{2}$. In conclusion, if we make sure that the approximation error of the scaling and wavelet functions satisfies $\delta \ll n^{- 9 / 2} / \log (n)$, the numerical error will not affect the statistical rate of convergence. Therefore, without loss of generality, we simply assume that we can exactly evaluate the scaling and wavelet functions in the rest of the paper.

In the above calculation we did not try to find the largest $\delta$ such that the statistical rate is unaffected so the above requirement on $\delta$ may very well not be the weakest possible. Finally, we remark that unlike wavelets, numerical accuracy will be an integral part of the analysis if one wants to generalize the theory in this paper to the case in which nuisance functions are estimated by deep neural networks trained with (stochastic) gradient descent or its variants.
\end{rem}

\section{Proof of Theorem \ref{theorem_general_lepski}}\label{section_appendix_proof_of_theorems}
\phantomsection
\begin{proof}
In this proof we repeatedly use the fact that for any fixed $m\in \mathbb{N}$ and $a_1,\ldots,a_m$ real numbers, one has by H\"{o}lder's inequality $\vert a_1 + \ldots + a_m \vert^p\leq C(m, p) \left( |a_1|^p + \ldots + |a_m|^p \right)$ for $p > 1$. Fix $\theta$ such that $\{ P \in \mathcal{P}_\theta: f_1(\theta) = \tau, f_2(\theta) > \min \{ 4 \tau / (4 \tau + d), 1 / 2 \} \}$ and our result will hold over all $\tau$ such that $f_1(\theta) = \tau, f_2(\theta) > \min \{ 4 \tau / (4 \tau + d), 1 / 2 \}$. Throughout the proof, we use $C, C', C''$ (with subscripts $1$, $2$, etc.) to denote generic universal constants only depending on the known parameters of the data generating mechanism (and $C_{\mathrm{Lepski}}$, which also only depends on the known parameters). When it is clear from the derivation, the same constant notation may change its value from line to line. We will highlight which known parameters a constant depends on when it is necessary for understanding.
	
	The corresponding oracle choice of $k$ is defined through balancing bias-variance tradeoff:
\begin{equation}
\kustar \coloneqq \left\{ \begin{array}{ll}
\min \left\{ k \in \mathcal{K}: k^{- 2 \tau / d} \leq \frac{1}{C_{\mathrm{Lepski}}^{1 / 2}} R(k)^{1 / 2} d(k), k > k_2 \right\} & \tau / d < 1 / 4 \\
k_2 \equiv \dfrac{n_1}{\ell(n_1)} & \tau / d = 1 / 4 \\
k_1 \equiv \dfrac{n_1}{\log{n_1}} & \tau / d > 1 / 4
\end{array} \right.
\end{equation}
where $\kustar$ is of order $\left( n_1 / \sqrt{\log{n_1}} \right)^{2 / (4 \tau / d + 1)}$ when $\tau / d < 1 / 4$ up to a constant depending on $C_{\mathrm{Lepski}}$. By construction, we have, for $n_1$ large enough
$$
\kustar{}^{- 2 \tau / d} \le \frac{1}{C_{\mathrm{Lepski}}^{1 / 2}} R(\kustar)^{1 / 2} d(\kustar).
$$	
	A key step in bounding the risk of an adaptive estimator through Lepski's adaptation scheme is to show that the probability of selecting a $k$ greater than necessary ($\kustar$) is small. To this end, we first prove the following lemma in order to bound $\P_P \left( \khat > \kustar \right)$, when Properties (A) and (B) are met:
\begin{lem}\label{lem:undersmoothing}
Under Properties (A) and (B), given a large enough universal constant $C_{\mathrm{Lepski}}$ depending on known parameters $(\tau_{\max}, C_A, C_A', C_B, B, B_U)$ of the data generating mechanism, for any $k > \kustar$, we have
\begin{align*}
\Pp \left( \hat{k} = k \right) \le \left\{ \begin{array}{ll}
C \log{n} \exp \{ - C d(k)^2 \} & \text{ if } k > k_2, \\
C \left( \exp \{ - C d(k_2)^2 \} + \log{n} \exp \{ - C d(k_3)^2 \} \right) & \text{ if } k = k_2.
\end{array} \right.
\end{align*}
for some constants $C > 0$ depending on $C_{\mathrm{Lepski}}$.
\end{lem}
\begin{proof}
For any $k$, let $k_-$ be the previous element of $k$ in the discretization set $\mathcal{K}$. Then
\begin{align*}
\Pp \left( \hat{k} = k \right) & = \Pp \left( \exists k' \geq k: \left\vert \hat\phi_{n, k'} - \hat\phi_{n, k_-} \right\vert > R(k')^{1 / 2} d(k') \right) \\
& \leq \sum_{k \le k' \le k_{N - 1}} \Pp \left( \left\vert \hat\phi_{n, k'} - \hat\phi_{n, k_-} \right\vert > R(k')^{1 / 2} d(k') \right).
\end{align*}
Recall the ``good'' event $\mathcal{I}_2 (n_2)$ introduced in Property (A). Then we decompose each summand on the RHS of the above display as follows:
\begin{align*}
& \; \BP_P \left( \left\vert \hat\phi_{n, k'} - \hat\phi_{n, k_-} \right\vert > R(k')^{1 / 2} d(k') \right) \\
= & \; \BP_P \left( \left\vert \hat\phi_{n, k'} - \hat\phi_{n, k_-} \right\vert > R(k')^{1 / 2} d(k'), \mathcal{I}_2 (n_2) \right) \\
& + \BP_P \left( \left\vert \hat\phi_{n, k'} - \hat\phi_{n, k_-} \right\vert > R(k')^{1 / 2} d(k'), \overline{\mathcal{I}_2 (n_2)} \right) \\
\leq & \; \underbrace{\BP_P \left( \left\vert \hat\phi_{n, k'} - \hat\phi_{n, k_-} \right\vert > R(k')^{1 / 2} d(k'), \mathcal{I}_2 (n_2) \right)}_{A} + \underbrace{\BP_P \left( \overline{\mathcal{I}_2 (n_2)} \right)}_{B}
\end{align*}
where $B \leq \frac{C_A' \log{n_2}}{n_2^2}$ by Property (A). We are only left to bound the first term $A$:
\begin{align*}
A = \TE_{P, 2} \left( \P_{P, 1} \left( \left\vert \hat\phi_{n, k'} - \hat\phi_{n, k_-} \right\vert > R(k')^{1 / 2} d(k') \right) \mathbbm{1} \left\{ \mathcal{I}_2 (n_2) \right\} \right).
\end{align*}
We attempt to control $\P_{P, 1} \left( \left\vert \hat\phi_{n, k'} - \hat\phi_{n, k_-} \right\vert > R(k')^{1 / 2} d(k') \right)$ through the following inequality decomposition:
\begin{align*}
& \; \P_{P, 1} \left( \left\vert \hat\phi_{n, k'} - \hat\phi_{n, k_-} \right\vert > R(k')^{1 / 2} d(k') \right) \\
= & \; \P_{P, 1} \left( \left\vert \Uhat_{n, k'} - \Uhat_{n, k_-} \right\vert > R(k')^{1 / 2} d(k') \right) \\
= & \; \underbrace{\P_{P, 1} \left( \left\vert \Uhat_{n, k_-} - \E_{P, 1} \Uhat_{n, k_-} \right\vert > \frac{1}{2} R(k')^{1 / 2} d(k') \right)}_{A_1} \\
& + \underbrace{\P_{P, 1} \left( \left\vert \Uhat_{n, k'} - \E_{P, 1} \Uhat_{n, k'} \right\vert > \frac{1}{2} R(k')^{1 / 2} d(k') - \left\vert \E_{P, 1} \Uhat_{n, k'} - \E_{P, 1} \Uhat_{n, k_-} \right\vert \right)}_{A_2},
\end{align*}
where we recall the definition of $\Uhat_{n, k}$ in equation \eqref{def:u}. Now we need to control terms $\E_{P, 2} [ A_1 \mathbbm{1} \{ \mathcal{I}_2 (n_2) \} ]$ and $\E_{P, 2} [ A_2 \mathbbm{1} \{ \mathcal{I}_2 (n_2) \} ]$ respectively. First we apply the exponential concentration inequality in Lemma \ref{lemma_ustat_tail_use} for U-statistics of order two to obtain the following upper bound on $\E_{P, 2} [ A_1 \mathbbm{1} \{ \mathcal{I}_2 (n_2) \} ]$:
\begin{align*}
\E_{P, 2} [ A_1 \mathbbm{1} \{ \mathcal{I}_2 (n_2) \} ] \le C_1 \exp \{ - C_1' d(k')^2 \},
\end{align*}
where $C_1$ and $C_1'$ are chosen sufficiently large depending on $C_A$ defined in Property (A), $C_B$ defined in Property (B), $(B, B_U)$ defined in Section \ref{section_main_upper_bound}, $J_0$ determined by $\tau_{\max}$, a known upper bound on the smoothness indices adapted over, and $C_{\mathrm{Lepski}}$ to be specified later.

Next we analyze $\TE_{P, 2} [ A_2 \mathbbm{1} \{ \mathcal{I}_2 (n_2) \} ]$. To do this, we need to first control $\left\vert \TE_{P, 1} \hat\phi_{n, k'} - \TE_{P, 1} \hat\phi_{n, k_-} \right\vert$ within the ``good'' event $\mathcal{I}_2 (n_2)$.

\begin{align*}
& \; \left\vert \TE_{P, 1} \Uhat_{n, k'} - \TE_{P, 1} \Uhat_{n, k_-} \right\vert \mathbbm{1} \{ \mathcal{I}_{2} (n_2) \} \\
= & \; \left\vert \TE_{P, 1} \hat\phi_{n, k'} - \TE_{P, 1} \hat\phi_{n, k_-} \right\vert \mathbbm{1} \{ \mathcal{I}_{2} (n_2) \} \\
\leq & \; 2 C_A \left( n_2^{ - f_2(\theta)} + k_-^{- \frac{2 f_1(\theta)}{d}} \right)
\end{align*}
where the last inequality is simply a consequence of Property (A). Since $k > \kustar$, $k_- = k - 1 \ge \kustar$, and thus for some absolute constant $C$ depending on $C_A$,
\begin{align*}
\left\vert \TE_{P, 1} \Uhat_{n, k'} - \TE_{P, 1} \Uhat_{n, k_-} \right\vert \mathbbm{1} \{ \mathcal{I}_{2} (n_2) \} & \leq \left\{ 
\begin{array}{ll}
C {\kustar}^{- \frac{2 f_1(\theta)}{d}} & \text{ if } \tau / d \leq 1 / 4 \\
C n^{- f_2(\theta)} & \text{ if } \tau / d > 1 / 4 \\
\end{array} \right. \\
& \leq \frac{1}{4} R(\kustar)^{1 / 2} d(\kustar) \le \frac{1}{4} R(k')^{1 / 2} d(k')
\end{align*}
where we use $n_2 = n / \log{n}$, the definition of $\kustar$ and $f_2(\theta) > 1 / 2$ if $\tau / d > 1 / 4$ and we choose $C_{\mathrm{Lepski}}$ large enough depending on $C_A$.

Then again Lemma \ref{lemma_ustat_tail_use} implies:
\begin{align*}
\TE_{P, 2} [ A_2 \mathbbm{1} \{ \mathcal{I}_{2} (n_2) \} ] \leq C_2 \exp \{ - C_2' d(k')^2 \}
\end{align*}
where similarly $C_2$ and $C_2'$ are chosen sufficiently large depending on $C_A$ defined in Property (A) (which in turn determines $C_{\mathrm{Lepski}}$), $C_B$ defined in Property (B), $(B, B_U)$ defined in Section \ref{section_main_upper_bound}, and $J_0$ determined by $\tau_{\max}$, a known upper bound of the smoothness index adapted over.

Combining the bounds on $\TE_{P, 2} [ A_1 \mathbbm{1} \{ \mathcal{I}_{2} (n_2) \} ]$ and $\TE_{P, 2} [ A_2 \mathbbm{1} \{ \mathcal{I}_{2} (n_2) \} ]$, we have
\begin{align*}
& \BP_P \left( \left\vert \phihat_{n, k'} - \phihat_{n, k_-} \right\vert > R(k')^{1 / 2} d(k') \right) \\
& \leq C \exp \{ - C' d(k')^2 \} + C'' \frac{\log{n}}{n^2} \\
& \leq \left\{ \begin{array}{ll}
\dfrac{C \ell(n_1)}{n_1^\delta} \leq \dfrac{C \ell(n)}{n^\delta} & k' > k_2, \\
C \exp\{ - \ell^{-1}(n_1) \} \leq C \exp \{ - \ell^{-1}(n) \} & k' = k_2
\end{array} \right.
\end{align*}
for some absolute constant $C$ depending on $(C_A, C_A', C_B, B, B_U, \tau_{\max}, C_{\mathrm{Lepski}})$.

Finally, recall that in the beginning of the proof we bound $\P_P \left( \khat = k \right)$ by
\begin{align*}
\Pp \left( \hat{k} = k \right) & = \Pp \left( \exists k' \geq k: \left\vert \hat\phi_{n, k'} - \hat\phi_{n, k_-} \right\vert > R(k')^{1 / 2} d(k') \right) \\
& \leq \sum_{k \le k' \le k_{N - 1}} \Pp \left( \left\vert \hat\phi_{n, k'} - \hat\phi_{n, k_-} \right\vert > R(k')^{1 / 2} d(k') \right),
\end{align*}
and there are at most $O(\log{n_1}) = O(\log{n})$ summands in total, therefore when $k > k_2$
$$
\P_P \left( \hat{k} = k \right) \leq C \log{n} \max_{k \le k' \le k_{N - 1}} \frac{\ell(n_1)}{n_1^\delta} \leq C \frac{\log{(n)} \ell(n)}{n^{\delta}},
$$
whereas when $k = k_2$,
\begin{align*}
\P_P \left( \hat{k} = k_2 \right) & \leq \Pp \left( \left\vert \hat\phi_{n, k_2} - \hat\phi_{n, k_1} \right\vert > R(k_2)^{1 / 2} d(k_2) \right) \\
& + \sum_{k_3 \le k' \le k_{N - 1}} \Pp \left( \left\vert \hat\phi_{n, k'} - \hat\phi_{n, k_1} \right\vert > R(k')^{1 / 2} d(k') \right) \\
& \leq C \left( \exp \{ - \ell^{-1}(n) \} + \frac{\log{(n)} \ell(n)}{n^{\delta}} \right).
\end{align*}
\end{proof}

	We further introduce the notation for the mean zero random variable $\IF \coloneqq L_1(O) - \phi(P)$, which encodes the first-order influence function \citep{van2000asymptotic} of $\phi(P)$ evaluated at the true nuisance functions. Returning to the proof of Theorem \ref{theorem_general_lepski}, our goal is to obtain an upper bound on the following risk of our adaptive estimator $\phihat_{n, \khat}$:
	\begin{align*}
	& \; \sup_{P \in \Par_{\theta}: \atop f_1(\theta) = \tau, f_2 (\theta) > \min \left\{ \frac{4 \tau}{4 \tau + d}, \frac{1}{2} \right\}} \E_{P} \left[ \left( \phihat_{n, \khat} - \phi(P) - \frac{1}{n_1} \sum_{i = 1}^{n_1} \IF_i \right)^2 \right] \\
	= & \; \sup_{P \in \Par_{\theta}: \atop f_1(\theta) = \tau, f_2 (\theta) > \min \left\{ \frac{4 \tau}{4 \tau + d}, \frac{1}{2} \right\}} \E_{P, 2} \left[ \E_{P, 1} \left[ \left( \phihat_{n, \khat} - \phi(P) - \frac{1}{n_1} \sum_{i = 1}^{n_1} \IF_i \right)^2 \right] \right].
	\end{align*}
	We start with controlling the following conditional risk bound:
	\begin{align*}
	 & \; \E_{P, 1} \left[ \left( \phihat_{n, \khat} - \phi(P) - \frac{1}{n_1} \sum_{i = 1}^{n_1} \IF_i \right)^2 \right] \\
	 \leq & \; \underbrace{\E_{P, 1} \left[ \left( \phihat_{n, \khat} - \phi(P) - \frac{1}{n_1} \sum_{i = 1}^{n_1} \IF_i \right)^2 \mathbbm{1} \left\{ \khat \leq \kustar \right\} \right]}_{T_1} \\
	 & + \underbrace{\E_{P, 1} \left[ \left( \phihat_{n, \khat} - \phi(P) - \frac{1}{n_1} \sum_{i = 1}^{n_1} \IF_i \right)^2 \mathbbm{1} \left\{ \khat > \kustar \right\} \right]}_{T_2}.
	\end{align*}
	Below we control the terms $T_1$ and $T_2$ separately.
\subsection*{Control of $T_1$}
We decompose $T_1$ as follows.
\be
T_1 & = \E_{P, 1} \left[ \mathbbm{1} \left\{ \khat \leq \kustar \right\} \left( \phihat_{n, \khat} - \phi(P) - \frac{1}{n_1} \sum_{i = 1}^{n_1} \IF_i \right)^2 \right] \\
& \leq 4 \E_{P, 1} \left[ \mathbbm{1} \left\{ \khat \leq \kustar \right\} \left\{ \begin{array}{c}
\underbrace{\left( \phihat_{n, \khat} - \phihat_{n, \kustar} \right)^2}_{T_{11}} + \underbrace{\left( \E_{P, 1} \phihat_{n, \kustar} - \phi(P) \right)^2}_{T_{12}} \\
+ \left( \phihat_{n, \kustar} - \E_{P, 1} \phihat_{n, \kustar} - \frac{1}{n_1} \sum_{i = 1}^{n_1} \IF_i \right)^2
\end{array} \right\} \right] \\
& \leq 4 \left\{ \E_{P, 1} \left[ \mathbbm{1} \left\{ \khat \leq \kustar \right\} T_{11} \right] + T_{12} + \underbrace{\E_{P, 1} \left[ \left( \phihat_{n, \kustar} - \E_{P, 1} \phihat_{n, \kustar} - \frac{1}{n_1} \sum_{i = 1}^{n_1} \IF_i \right)^2 \right]}_{T_{13}} \right\}.
\ee
Now we control $T_{11}, T_{12}$ and $T_{13}$ separately. First, by definition of $\khat$, for some absolute constant $C$ depending on $C_{\mathrm{Lepski}}$,
\begin{align*}
T_{11} & \le R(\kustar) d(\kustar)^2 = \frac{\kustar}{n_1^2} d(\kustar)^2 \\
& \le \left\{ \begin{array}{ll}
C \left( \dfrac{n_1}{\sqrt{\log{n_1}}} \right)^{- \frac{8 \tau / d}{4 \tau / d + 1}} \asymp C \left( \dfrac{n}{\sqrt{\log{n}}} \right)^{- \frac{8 \tau / d}{4 \tau / d + 1}} & \tau / d < 1 / 4, \\
\dfrac{\ell^{-2}(n_1)}{n_1} \asymp \dfrac{\ell^{-2}(n)}{n} & \tau / d = 1 / 4, \\
\dfrac{\ell^{-1}(n_1)}{n_1 \log{n_1}} \asymp \dfrac{\ell^{-1}(n)}{n \log{n}} \ll \dfrac{1}{n} & \tau / d > 1 / 4.
\end{array} \right.
\end{align*} 
The above upper bound also applies to $\E_{P, 1} T_{11}$ and $\E_P T_{11}$. In terms of $T_{13}$, by Lemma \ref{lemma_ustat_tail_use} (ii), we have
\begin{align*}
\E_{P, 2} T_{13} & \le 2 \E_{P, 2} \E_{P, 1} \left[ \begin{array}{c}
\left( \Uhat_{n, \kustar} - \E_{P, 1} \Uhat_{n, \kustar} \right)^2 \\
+ \left( \frac{1}{n_1} \sum_{i = 1}^{n_1} \Ltilde_1(O_i) - \E_{P, 1} \Ltilde_1(O_i) - \IF_i \right)^2 
\end{array} \right] \\
& \le C \left( \frac{\kustar}{n_1^2} + o \left( \frac{1}{n_1} \right) \right)
\end{align*}
for some absolute constant $C$ depending on $(C_B, B, B_U, \tau_{\max})$, the known parameters of the problem. The fact that $$\E_{P, 2} \E_{P, 1} \left[ \left( \frac{1}{n_1} \sum_{i = 1}^{n_1} \Ltilde_1(O_i) - \E_{P, 1} \Ltilde_1(O_i) - \IF_i \right)^2 \right]$$ scales as $o(n_1^{-1})$ can be seen from the following argument. Within the ``good'' event $\mathcal{I}_{2} (n_2)$, this term is $n_1^{-1} \E_{P, 1} [ ( \Ltilde_1(O) - L_1(O) - \E_{P, 1} ( \Ltilde_1(O) - L_1(O) ) )^2 ] = o(n_1^{-1})$ because $\E_{P, 1} [ ( \Ltilde_1(O) - L_1(O) )^2 ] = o(1)$ in $\mathcal{I}_{2} (n_2)$. Outside the ``good'' event $\mathcal{I}_{2} (n_2)$, $\E_{P, 1} \left[ \left( \frac{1}{n_1} \sum_{i = 1}^{n_1} \Ltilde_1(O_i) - \E_{P, 1} \Ltilde_1(O_i) - \IF_i \right)^2 \right]$ scales as $$n_1^{-1} \E_{P, 1} [ ( \Ltilde_1(O) - L_1(O) - \E_{P, 1} ( \Ltilde_1(O) - L_1(O) ) )^2 ] = O(n_1^{-1})$$ because now we only know $\E_{P, 1} [ ( \Ltilde_1(O) - L_1(O) - \E_{P, 1} ( \Ltilde_1(O) - L_1(O) ) )^2 ]$ is bounded almost surely. However, since the probability outside the ``good'' event is $o(1)$ by Property (A), $$\E_{P, 2} \E_{P, 1} \left[ \left( \frac{1}{n_1} \sum_{i = 1}^{n_1} \Ltilde_1(O_i) - \E_{P, 1} \Ltilde_1(O_i) - \IF_i \right)^2 \mathbbm{1} \left\{ \overline{\mathcal{I}_2 (n_2)} \right\} \right] = o(n_1^{-1}).$$ Then combining the above two points, we see that $\E_{P, 2} \E_{P, 1} [ ( \frac{1}{n_1} \sum_{i = 1}^{n_1} \Ltilde_1(O_i) - \E_{P, 1} \Ltilde_1(O_i) - \IF_i )^2 ]$ scales as $o(n_1^{-1})$.

As a result, there exists some universal positive constant $C$ depending only on known parameters such that
\begin{align*}
& \sup_{P \in \Par_{\theta}: \atop f_1(\theta) = \tau, f_2 (\theta) > \min \left\{ \frac{4 \tau}{4 \tau + d}, \frac{1}{2} \right\}} \E_{P, 2} T_{13} \le C \left( \frac{\kustar}{n_1^2} + o \left( \frac{1}{n_1} \right) \right) \\
& \le \left\{ \begin{array}{ll}
\dfrac{C \kustar}{n_1^2} \asymp \dfrac{C \kustar}{n^2} & \tau / d \le 1 / 4, \\
o(n_1^{-1}) = o(n^{-1}) & \tau / d > 1 / 4.
\end{array} \right.
\end{align*}

In terms of $T_{12}$, we again need to analyze its behavior under two different scenarios: 
\begin{enumerate}
\item Within the ``good'' event $\mathcal{I}_{2} (n_2)$, for some absolute constant $C$ depending on $C_A$:
\begin{align*}
& \; \E_{P, 2} T_{12} \mathbbm{1} \left\{ \mathcal{I}_{2} (n_2) \right\} \\
\le & \; C_A \left( n_2^{- 2 f_2(\theta)} + \kustar{}^{- \frac{4 f_1(\theta)}{d}} \right) \\
\le & \; \left\{ \begin{array}{ll}
C \left( \dfrac{n_1}{\sqrt{\log{n_1}}} \right)^{- \frac{8 \tau / d}{4 \tau / d + 1}} \asymp C \left( \dfrac{n}{\sqrt{\log{n}}} \right)^{- \frac{8 \tau / d}{4 \tau / d + 1}} & \tau / d < 1 / 4, \\
\dfrac{C}{n} & \tau / d = 1 / 4, \\
o(n^{-1}) & \tau / d > 1 / 4.
\end{array} \right.
\end{align*}
\item Outside the ``good'' event $\mathcal{I}_{2} (n_2)$: we still have $T_{12} < B_T$ for some constant $B_T > 0$ depending on $B$ and $B_U$ almost surely.
\end{enumerate}
Summarizing the above arguments, we have
\begin{align*}
\E_{P, 2} T_{12} & = \E_{P, 2} T_{12} \mathbbm{1} \left\{ \overline{\mathcal{I}_{2} (n_2)} \right\} + \E_{P, 2} T_{12} \mathbbm{1} \left\{ \mathcal{I}_{2} (n_2) \right\} \\
& \leq B_T \P_{P, 2} \left( \overline{\mathcal{I}_{2} (n_2)} \right) + C_A \left( n_2^{- 2 f_2(\theta)} + \kustar{}^{- \frac{4 f_1(\theta)}{d}} \right) \\
& \leq \left\{ \begin{array}{ll} 
C \left( \dfrac{n}{\sqrt{\log{n}}} \right)^{- \frac{8 \tau / d}{4 \tau / d + 1}} & \tau / d < 1 / 4, \\
\dfrac{C}{n} & \tau / d = 1 / 4, \\
o(n^{-1}) & \tau / d > 1 / 4.
\end{array} \right.
\end{align*}
Combining the above bounds for $\E_{P, 2} \left[ \mathbbm{1} \left\{ \khat \leq \kustar \right\} T_{11} \right], \E_{P, 2} T_{12}$ and $\E_{P, 2} T_{13}$, we get, for some universal positive constant $C > 0$,
\allowdisplaybreaks
\begin{align*}
T_1 \le \left\{ \begin{array}{ll}
C \left( \dfrac{n}{\sqrt{\log{n}}} \right)^{- \frac{8 \tau / d}{4 \tau / d + 1}} & \tau / d < 1 / 4, \\
\dfrac{C \ell^{-2}(n)}{n} & \tau / d = 1 / 4, \\
o(n^{-1}) & \tau / d > 1 / 4.
\end{array} \right.
\end{align*}

\subsection*{Control of $T_2$}
For some $p, q > 1$ such that $\frac{1}{p} + \frac{1}{q} = 1$, we have
\begin{align*}
T_2 & \le \sum_{\kustar < k \le k_{N - 1}} \E_{P, 1} \left[ \left( \phihat_{n, \khat} - \phi(P) - \frac{1}{n_1} \sum_{i = 1}^{n_1} \IF_i \right)^2 \mathbbm{1} \left\{ \khat = k \right\} \right] \\
& = \sum_{\kustar < k \le k_{N - 1}} \E_{P, 1} \left[ \left( \phihat_{n, k} - \phi(P) - \frac{1}{n_1} \sum_{i = 1}^{n_1} \IF_i \right)^2 \mathbbm{1} \left\{ \khat = k \right\} \right] \\
& \le 2 \left( \E_{P, 1} \phihat_{n, \kustar} - \phi(P) \right)^2 \\
& + 2 \sum_{\kustar < k \le k_{N - 1}} \E_{P, 1} \left[ \left( \phihat_{n, k} - \E_{P, 1} \phihat_{n, k} - \frac{1}{n_1} \sum_{i = 1}^{n_1} \IF_i \right)^2 \mathbbm{1} \left\{ \khat = k \right\} \right] \\
& \le 2 \underbrace{\left( \E_{P, 1} \phihat_{n, \kustar} - \phi(P) \right)^2}_{T_{21}} \\
& + 2 \underbrace{\sum_{\kustar < k \le k_{N - 1}} \E_{P, 1}^{\frac{1}{q}} \left[ \left( \phihat_{n, k} - \E_{P, 1} \phihat_{n, k} - \frac{1}{n_1} \sum_{i = 1}^{n_1} \IF_i \right)^{2q} \right] \P_{P, 1}^{\frac{1}{p}} \left( \khat = k \right)}_{T_{22}}.
\end{align*}
We next control $T_{21}$ and $T_{22}$ separately. First we observe that $T_{21}$ is exactly $T_{12}$. Thus applying the same bound, we have
\begin{align*}
\E_{P, 2} T_{21} \le \left\{ \begin{array}{ll}
C \left( \dfrac{n}{\sqrt{\log{n}}} \right)^{- \frac{8 \tau / d}{4 \tau / d + 1}} & \tau / d < 1 / 4, \\
\dfrac{C}{n} & \tau / d = 1 / 4, \\
o(n^{-1}) & \tau / d > 1 / 4.
\end{array} \right.
\end{align*}

In the control of $T_{22}$, we utilize the result in Lemma \ref{lem:undersmoothing}. 
\begin{align*}
T_{22} \le & \; C(q) \sum_{\kustar < k \le k_{N - 1}} \E_{P, 1}^{\frac{1}{q}} \left[ \begin{array}{c}
\left( \Uhat_{n, k} - \E_{P, 1} \Uhat_{n, k} \right)^{2q} \\
+ \left( \frac{1}{n_1} \sum_{i = 1}^{n_1} \Ltilde_1(O_i) - \E_{P, 1} \Ltilde_1(O) - \IF_i \right)^{2q}
\end{array} \right] \P_{P, 1}^{\frac{1}{p}} \left( \khat = k \right) \\
\le & \; C'(q) \sum_{\kustar < k \le k_{N - 1}} \frac{k}{n^2} \P_{P, 1}^{\frac{1}{p}} \left( \khat = k \right)
\end{align*}
where the first inequality is due to H\"{o}lder inequality, the second inequality follows from Lemma \ref{lemma_ustat_tail_use} (ii) with $k > \kustar$ implying that $k > n_1 = n(1 - o(1))$, and the third inequality is a consequence of Lemma \ref{lem:undersmoothing}. Here $C'(q)$ absorbs the constants (depending on $(B, B_U, \tau_{\max}, C_B)$) in front of the rate (which is $k / n^2$) of $$\E_{P, 1}^{\frac{1}{q}} \left[ \left( \Uhat_{n, k} - \E_{P, 1} \Uhat_{n, k} \right)^{2q} + \left( \frac{1}{n_1} \sum_{i = 1}^{n_1} \Ltilde_1(O_i) - \E_{P, 1} \Ltilde_1(O) - \IF_i \right)^{2q} \right]$$ into $C(q)$.

Next we divide our analyses into three different scenarios:
\begin{enumerate}
\item $\tau / d < 1 / 4$: $\kustar > k_2$. Then following Lemma \ref{lem:undersmoothing}, we have
\allowdisplaybreaks
\begin{align*}
T_{22} & \le \frac{C'(q)}{n^2} \log^{\frac{1}{p}}{(n)} \sum_{\kustar < k \le k_{N - 1}} k \exp \{ - C d(k)^2 / p \} \\
& \le \frac{C'(q) n^{\frac{C C_{\mathrm{Lepski}}^2 (1 - \delta)}{p}}}{n^2} \log^{\frac{1}{p}}{(n)} \sum_{\kustar < k \le k_{N - 1}} k \cdot k^{- C C_{\mathrm{Lepski}}^2 / p} \\
& = \frac{C'(q) n^{2 (1 - \delta)}}{n^2} \log^{\frac{1}{p}}{(n)} \sum_{\kustar < k \le k_{N - 1}} k \cdot k^{-2} \\
& \le C'(q) n^{-2\delta} \log^{\frac{1}{p}}{(n)} C' \log{n} \kustar{}^{-1} \\
& \le C'(q) C' n^{-2\delta} \log^{\frac{1 + p}{p}}{(n)} \left( \frac{n}{\sqrt{\log{n}}} \right)^{- \frac{2}{4 \tau / d + 1}} \\
& = o \left( \left( \frac{n}{\sqrt{\log{n}}} \right)^{- \frac{8 \tau / d}{4 \tau / d + 1}} \right)
\end{align*}
where the second line follows by plugging in the definition of $d(k)$ when $k > \kustar > k_2$, the third line holds by choosing the constant $p$ such that $C C_{\mathrm{Lepski}}^2 / p = 2$, the fourth line utilizes the cardinality of $\mathcal{K}$ being bounded by $\log{n}$ up to some constant $C' > 0$ and $k > \kustar$, the fifth line follows from the definition of $\kustar$ when $\tau / d < 1 / 4$ and in the sixth line we simply use $8 \tau / d < 2$. Here we want to emphasize that the choice of $p$ depends only on the constants $C$ in Lemma \ref{lem:undersmoothing} and $C_{\mathrm{Lepski}}$, not on the sample size $n$ or $k$.

\item $\tau / d = 1 / 4$: $\kustar = k_2 = n_1 / \ell(n_1)$. Then similar to the calculations above by setting $p$ such that $C_1' C_{\mathrm{Lepski}}^2 / p = 2$ and choosing $0.5 < \delta < 1$, we have
\allowdisplaybreaks
\begin{align*}
T_{22} & \le \frac{C'(q) C' n^{2(1 - \delta)}}{n^2} \log^{\frac{1 + p}{p}}{(n)} k_2^{-1} \\
& \le C'(q) C' n^{-2\delta} \log^{\frac{1 + p}{p}}{(n)} \ell(n) \\
& = o(n^{-1}).
\end{align*}

\item $\tau / d > 1 / 4$: $\kustar = k_1 = n_1 / \log{n_1}$. Then compared to $T_{22}$ in $\tau / d \le 1 / 4$, we need to add one extra term $$\E_{P, 1}^{\frac{1}{q}} \left[ \left( \phihat_{n, k_2} - \E_{P, 1} \phihat_{n, k_2} - \frac{1}{n_1} \sum_{i = 1}^{n_1} \IF_i \right)^{2q} \right] \P_{P, 1}^{\frac{1}{p}} \left( \khat = k_2 \right):$$
\begin{align*}
T_{22} \le & \; \frac{C'(q)}{n^2} \log^{\frac{1}{p}}{(n)} \sum_{\kustar < k \le k_{N - 1}} k \exp \{ - C d(k)^2 / p \} \\
& + \E_{P, 1}^{\frac{1}{q}} \left[ \left( \phihat_{n, k_2} - \E_{P, 1} \phihat_{n, k_2} - \frac{1}{n_1} \sum_{i = 1}^{n_1} \IF_i \right)^{2q} \right] \P_{P, 1}^{\frac{1}{p}} \left( \khat = k_2 \right) \\
\le & \; o(n^{-1}) + \E_{P, 1}^{\frac{1}{q}} \left[ \left( \phihat_{n, k_2} - \E_{P, 1} \phihat_{n, k_2} - \frac{1}{n_1} \sum_{i = 1}^{n_1} \IF_i \right)^{2q} \right] \P_{P, 1}^{\frac{1}{p}} \left( \khat = k_2 \right) \\
\le & \; o(n^{-1}) + C (n \ell(n))^{-1} \left( \exp\{ - C d(k_2)^2 \} + \log{n} \exp\{ - C d(k_3)^2 \} \right)^{\frac{1}{p}} \\
= & \; o(n^{-1})
\end{align*}
where the last line follows from the calculation as in the case $\tau / d = 1 / 4$.
\end{enumerate}

Hence combining the above findings on the control of $T_{22}$, we have
\begin{align*}
\E_{P, 2} T_{22} = \left\{ \begin{array}{ll}
o \left( \frac{n}{\sqrt{\log{n}}} \right)^{- \frac{8 \tau / d}{4 \tau / d + 1}} & \tau / d < 1 / 4, \\
o(n^{-1}) & \tau / d \ge 1 / 4.
\end{array} \right.
\end{align*}
Finally, summarizing the above results on $T_1$ and $T_2$, we have
\begin{align*}
& \; \E_{P} \left[ \left( \phihat_{n, \khat} - \phi(P) - \frac{1}{n_1} \sum_{i = 1}^{n_1} \IF_i \right)^2 \right] \le \E_{P, 2} ( T_1 + T_2 ) \\
\le & \; \left\{ \begin{array}{ll}
O \left( \left( \dfrac{n}{\sqrt{\log{n}}} \right)^{- \frac{8 \tau / d}{4 \tau / d + 1}} \right) & \tau / d < 1 / 4, \\
O \left( \dfrac{1}{n \ell^2(n)} \right) & \tau / d = 1 / 4, \\
o \left( \dfrac{1}{n} \right) & \tau / d > 1 / 4.
\end{array} \right.
\end{align*}
In particular, when $\tau / d > 1 / 4$, the above result implies that
\begin{align*}
\E_{P} \left[ \left\vert \phihat_{n, \khat} - \phi(P) - \frac{1}{n_1} \sum_{i = 1}^{n_1} \IF_i \right\vert \right] = o(n^{- 1 / 2})
\end{align*}
which further implies that
\begin{align*}
\sqrt{n} \left( \phihat_{n, \khat} - \phi(P) \right) = \frac{1}{\sqrt{n_1}} \sum_{i = 1}^{n_1} \IF_i + o_P(1) \rightarrow_{d} \mathcal{N}(0, \sigma^2).
\end{align*}
\end{proof}

\section*{Acknowledgement}
We would like to thank two anonymous referees, the action editor G\'{a}bor Lugosi, and Zhi-Qin John Xu (Shanghai Jiao Tong University) for helpful comments on this paper.

\section*{Funding Information}
Lin Liu and James M. Robins were supported by the U.S. Office of Naval Research (ONR) grant N000141912446, and National Institutes of Health (NIH) awards R01AG057869 and R01AI127271. Lin Liu was also partially sponsored by Shanghai Pujiang Program Research grant 20PJ1408900, Shanghai Municipal Science and Technology Major Project 2021SHZDZX0102, and Major Program of National Natural Science Foundation of China 12090024. Rajarshi Mukherjee’s research was partially supported by NSF Grant EAGER-1941419. Eric Tchetgen Tchetgen was supported by NIH awards R01GM139926, R01AG065276, R01CA222147 and R01AI27271.


\bibliography{biblio_adaptation}
\renewcommand{\theHsection}{A\arabic{section}}
\appendix

\section{Proof of (\ref{EQ:EQUIV})} \label{app:equiv}\label{section_appendix_hellinger}
	\begin{proof}
		For any two probability measures $\nu_1$ and $\nu_2$, denote their total-variation distance and Hellinger distance as $\text{TV} \left( \nu_1, \nu_2 \right) \coloneqq \frac{1}{2} \int \vert d \nu_1 - d \nu_2 \vert $ and $H ( \nu_1, \nu_2 ) \coloneqq \left( \int ( \sqrt{d \nu_1} - \sqrt{d \nu_2} )^2 \right)^{1 / 2}$ respectively. Since we assume $d \nu_2$ is strictly bounded from below by $\underline{B}$, we have
		\begin{align*}
		\chi^2 \left( \nu_1, \nu_2 \right) = \int \frac{(d \nu_1 - d \nu_2)^2}{d \nu_1} \le \underline{B}^{-1} \int (d \nu_1 - d \nu_2)^2 \le 4 \underline{B}^{-1} \text{TV}^2 \left( \nu_1, \nu_2 \right) \le 4 \underline{B}^{-1} H^2 \left( \nu_1, \nu_2 \right)
		\end{align*}
		where in the last inequality we apply the well-known Le Cam's inequality that the Hellinger distance upper bounds the total variation distance \citep[Lemma 2.3]{tysbakov2009book}. Furthermore, using a simple Taylor expansion argument, we have
		\begin{align*}
		\chi^2 ( \nu_1, \nu_2 ) \leq 4 \underline{B}^{-1} H^2 \left( \nu_1, \nu_2 \right) \le \exp \{ 4 \underline{B}^{-1} H^2 \left( \nu_1, \nu_2 \right) \} - 1.
		\end{align*}
	\end{proof}
	
	\section{Proof of Remaining Theorems}\label{section:proof_of_remaining_theorems}
	\subsection*{Proof of Theorem \ref{thm_adaptive_density_regression}} 
	\begin{proof}
		Let
		$$2^{\jmin d}=\lfloor \left(\frac{n}{\log{n}}\right)^{\frac{1}{2 s_{\max} / d + 1}} \rfloor, \quad 2^{\jmax d}=\lfloor \left(\frac{n}{\log{n}}\right)^{\frac{1}{2 s_{\min} / d + 1}} \rfloor,$$ $$2^{\lmin d}=\lfloor \left(\frac{n}{\log{n}}\right)^{\frac{1}{2\gmax/d+1}} \rfloor, \quad 2^{\lmax d}=\lfloor \left(\frac{n}{\log{n}}\right)^{\frac{1}{2\gmin/d+1}} \rfloor.$$
		Without loss of generality assume that we have data $\{ \mathbf{x}_i, y_i \}_{i =1}^{2n}$. We split it into two disjoint and equal-sized parts and use the second part to construct the estimator $\hat{g}$ of the design density $g$  and use the resulting $\ghat$ to construct the adaptive estimates of the regression functions from the first half of the sample. 
		Throughout we choose the regularity of our wavelet bases to be larger than $\gamma_{\max}$ for the desired approximation and moment properties to hold. As a result our constants depend on $\gamma_{\max}$.
		
		Define $\mathcal{T}_1=[\jmin,\jmax]\cap \mathbb{N}$ and $\mathcal{T}_2=[\lmin,\lmax]\cap \mathbb{N}$. For $l \in \mathcal{T}_2$, let $\ghat_l(\bx)=\frac{1}{n}\sum_{i=n+1}^{2n}K_{V_l}\left(\bX_i,\bx\right)$. Now, let
		\be
		\lhat=\min\left\{j\in \mathcal{T}_2: \ \|\ghat_j-\ghat_{l}\|_{\infty}\leq C^*\sqrt{\frac{2^{ld}ld}{n}}, \ \forall l \in \mathcal{T}_2 \ \text{s.t.} \ l \geq j \right\}.
		\ee
		where $C^*$ is a constant  ({depending on $\gmax, B_U$}) that can be determined from the proof hereafter. Thereafter, consider the estimator $\gtilde\coloneqq\ghat_{\lhat}$. In the following, we first prove inequality \eqref{eq:moment_g}.
		
		Fix a $P \coloneqq (f, g) \in \mathcal{P}(s,\gamma)$. To analyze the estimator $\gtilde$, we begin with standard bias variance type  analysis for the candidate estimators $\ghat_l$. First note that for any $\bx \in [0,1]^d$, using standard facts about compactly supported wavelet basis having regularity larger than $\gamma_{\max}$ \citep{hardle1998wavelets}, one has, for a constant $C_1$ depending only on $q$ and the wavelet basis used, that
		\be
		|\E_P\left(\ghat_l(\bx)\right)-g(\bx)|&=|\Pi\left(g|V_l\right)(\bx)-g(\bx)|\leq C_1 M 2^{-ld\frac{\g}{d}}. \label{eqn:bias_lepski_density}
		\ee
		Above we have used the fact that 
		\be
		\sup\limits_{h \in H(\gamma,M)}\|h-\Pi(h|V_l)\|_{\infty}\leq C_1 M 2^{-l\gamma}. \label{eqn:holder_approx}
		\ee
		Also, by standard arguments about compactly supported wavelet basis having regularity larger than $\gamma_{\max}$ \citep{gine2016mathematical}, one has, for a constant $C_2 \coloneqq C(B_U,\bpsi_{0,0}^0, \bpsi_{0,0}^1,\gmax)$, that
		\be 
		\E_P\left(\|\ghat_l(\bx)-\E_P\left(\ghat_l(\bx)\right)\|_{\infty}\right)& \leq C_2\sqrt{\frac{2^{ld} ld}{n}}. \label{eqn:variance_lepski_density}
		\ee
		Therefore, by inequalities \eqref{eqn:bias_lepski_density} and \eqref{eqn:variance_lepski_density}, and triangle inequality,
		\be 
		\E_{P,2}\|\ghat_l-g\|_{\infty} & \leq C_1 M 2^{-ld\frac{\g}{d}}+C_2\sqrt{\frac{2^{ld} ld}{n}}.
		\ee 
		Define,
		\be
		\lstar & \coloneqq \min \left\{l \in \mathcal{T}_2: C_1 M 2^{-ld\frac{\g}{d}}\leq C_2\sqrt{\frac{2^{ld} ld}{n}}\right\}.
		\ee
		The definition of $\lstar$ implies that for $n$ sufficiently large,
		\be 
		2^d \left(\frac{C_1}{C_2}M\right)^{\frac{2d}{2\gamma+d}}\left(\frac{n}{\log{n}}\right)^{\frac{d}{2\gamma+d}}\leq 2^{\lstar d}\leq 2^{d+1} \left(\frac{C_1}{C_2}M\right)^{\frac{2d}{2\gamma+d}}\left(\frac{n}{\log{n}}\right)^{\frac{d}{2\gamma+d}}. \\ \label{eqn:lstar_size}
		\ee
		The error analysis of $\tilde{g}$ can now be carried out as follows:
		\be 
		\E_{P,2}\|\tilde{g}-g\|_{\infty}&=\E_{P,2}\|\tilde{g}-g\|_{\infty}\mathbbm{1}\left(\lhat\leq \lstar\right)+\E_{P,2}\|\tilde{g}-g\|_{\infty}\mathbbm{1}\left(\lhat>\lstar\right)\\& \coloneqq I + II. \label{eqn:lepski_decomposition_density}
		\ee
		We first control term $I$ as follows:
		\be
		I&=\E_{P,2}\|\tilde{g}-g\|_{\infty}\I\left(\lhat\leq \lstar\right)\\&\leq \E_{P,2}\|\ghat_{\lhat}-\ghat_{\lstar}\|_{\infty}\mathbbm{1}\left(\lhat\leq \lstar\right)+\E_{P,2}\|\ghat_{\lstar}-g\|_{\infty}\mathbbm{1}\left(\lhat\leq \lstar\right)\\
		& \leq C^{*}\sqrt{\frac{2^{\lstar d}\lstar d}{n}}+C_1 M 2^{-\lstar d\frac{\g}{d}}+C_2\sqrt{\frac{2^{\lstar d} \lstar d}{n}}\\
		& \leq (C^*+2C_2)\sqrt{\frac{2^{\lstar d}\lstar d}{n}}\leq 2^{d+1} \left(\frac{C_1}{C_2}M\right)^{\frac{2d}{2\gamma+d}} \left(\frac{n}{\log{n}}\right)^{-\frac{\gamma}{2\gamma+d}}.\\ \label{eqn:control_I_density}
		\ee
		The control of term $II$ is easier if one has suitable bounds on $\|\ghat_l-g\|_{\infty}$. To this end note that, for any fixed $\bx \in [0, 1]^d$, there exists a constant $C_3 \coloneqq C(\bpsi_{0,0}^0, \bpsi_{0,0}^1, \gmax)$ such that
		\be 
		|\ghat_{l}(\bx)|&\leq \frac{1}{n}\sum_{i=1}^n \sum_{m \in \Z_l}\sum_{v \in \{0,1\}^d}|\psi_{l, m}^v(\bX_i)||\psi_{l, m}^v(\bx)|\leq C_3 2^{ld}.
		\ee
		This along with the fact that $\|g\|_{\infty}\leq B_U$, implies that for $n$ sufficiently large,
		\be 
		\|\ghat_l-g\|_{\infty}\leq C_3 2^{ld}+B_U\leq 2C_3 2^{ld}.
		\ee
		In the above display the last inequality follows since $l \geq  \lmin \geq   \left(\frac{n}{\log{n}}\right)^{\frac{1}{2\gmax/d+1}}$. Therefore,
		\be 
		II & \leq C_3 \sum_{l = \lstar + 1}^{\lmax} 2^{ld} \P_{P, 2} \left( \lhat = l \right). \label{eqn:control_II_density_step1}
		\ee
		We now complete the control over II by suitably bounding $\P\left(\lhat=l\right)$. To this end, note that for any $l >\lstar$,
		\be
		\ & \mathbb{P}_{P,2}\left(\lhat=l\right) \\& \leq \sum_{l>\lstar}\mathbb{P}_{P,2}\left(\|\ghat_l-\ghat_{\lstar}\|_{\infty}> C^*\sqrt{\frac{2^{ld}ld}{n}}\right)\\
		&\leq \sum_{l>\lstar}\left\{\begin{array}{c}\mathbb{P}_{P,2}\left(\|\ghat_{\lstar}-\E\left(\ghat_{\lstar}\right)\|_{\infty}> \frac{C^*}{2}\sqrt{\frac{2^{ld}ld}{n}}-\|\E_{P,2}\left(\ghat_{\lstar}\right)-\E_{P,2}\left(\ghat_{l}\right)\|_{\infty}\right)\\+\mathbb{P}_{P,2}\left(\|\ghat_l-\E_{P,2} \left(\ghat_{l}\right)\|_{\infty}> \frac{C^*}{2}\sqrt{\frac{2^{ld}ld}{n}}\right)\end{array}\right\}\\
		& \leq \sum_{l>\lstar}\left\{ \begin{array}{c}\mathbb{P}_{P,2}\left(\|\ghat_{\lstar}-\E\left(\ghat_{\lstar}\right)\|_{\infty}> \frac{C^*}{2}\sqrt{\frac{2^{ld}ld}{n}}-\|\Pi\left(g|V_{\lstar}\right)-\Pi\left(g|V_{l}\right)\|_{\infty}\right)\\+\mathbb{P}\left(\|\ghat_l-\E_{P,2} \left(\ghat_{l}\right)\|_{\infty}> \frac{C^*}{2}\sqrt{\frac{2^{ld}ld}{n}}\right)\end{array}\right\}\\
		& \leq \sum_{l>\lstar}\left\{ \begin{array}{c}\mathbb{P}_{P,2}\left(\|\ghat_{\lstar}-\E\left(\ghat_{\lstar}\right)\|_{\infty}> \frac{C^*}{2}\sqrt{\frac{2^{ld}ld}{n}}-2C_2\sqrt{\frac{2^{l^*d}\lstar d}{n}}\right)\\+\mathbb{P}\left(\|\ghat_l-\E_{P,2} \left(\ghat_{l}\right)\|_{\infty}> \frac{C^*}{2}\sqrt{\frac{2^{ld}ld}{n}}\right)\end{array}\right\}\\
		& \leq \sum_{l>\lstar}\left\{ \begin{array}{c}\mathbb{P}_{P,2}\left(\|\ghat_{\lstar}-\E\left(\ghat_{\lstar}\right)\|_{\infty}> (\frac{C^*}{2}-2C_2)\sqrt{\frac{2^{ld}ld}{n}}\right)\\+\mathbb{P}\left(\|\ghat_l-\E_{P,2} \left(\ghat_{l}\right)\|_{\infty}> \frac{C^*}{2}\sqrt{\frac{2^{ld}ld}{n}}\right)\end{array}\right\}\\
		& \leq \sum_{l>l^*} 2\exp{\left(-Cld\right)}. \label{eqn:control_II_density_step2}
		\ee
		In the fourth and fifth lines of the above series of inequalities, we have used inequality \eqref{eqn:holder_approx} and the definition of $\lstar$ respectively. The last inequality in the above display holds for an absolute constant $C>0$ depending only on $B_U, \bpsi_{0,0}^0, \bpsi_{0,0}^1$ and the inequality follows from Lemma \ref{lemma_linear_projection_tail_bound} provided we choose  $C^*$ large enough depending on $M, B_U,\bpsi_{0,0}^0, \bpsi_{0,0}^1,\gmax$. In particular, this implies that, choosing $C^*$ large enough will guarantee that there exists a $\eta>3$  such that for large enough $n$, one has for any $l>\lstar$
		\be 
		\P(\lhat=l)\leq n^{-\eta}. \label{eqn:wrongchoicelepski_density}
		\ee
		This along with inequality \eqref{eqn:control_II_density_step1} and the choice of $\lmax$ implies that
		\be 
		II&\leq C_3\sum_{l>\lstar}2^{ld}n^{-\eta}=C_3\sum_{l>\lstar}\frac{2^{ld}}{n}n^{-\eta+1}\leq \frac{\lmax}{n^{\eta-1}}\leq \frac{\log{n}}{n}. \label{eqn:control_II_density_step3}
		\ee
		
		Finally combining equations \eqref{eqn:control_I_density} and  \eqref{eqn:control_II_density_step3},  we have the existence of an estimator $\gtilde$ depending on $M,B_U$, and $\gmax$ (once we have fixed our choice of the scaling and wavelet functions; see Section \ref{section_wavelets and function spaces} in the main text), such that for every $(s,\gamma) \in [s_{\min},s_{\max}]\times [\gmin,\gmax]$,
		$$\sup\limits_{P \in \mathcal{P}(s,\gamma)} \E_P\|\gtilde-g\|_{\infty} \leq  \left(C\right)^{\frac{d}{2\gamma+d}}\left(\frac{n}{\log{n}}\right)^{-\frac{\g}{2\g+d}},$$ 
		with a large enough positive $C$ depending on  $M, B_U$, and $\gmax$.
		
		We next show that uniformly over $P \in \mathcal{P}(s, \gamma) $, $\gtilde$ belongs to $H(\gamma,C)$ with probability at least $1-1/n^2$, for a large enough constant $C$ depending on  $M, B_U$, and $\gmax$. Towards this end, note that, for any $C>0$, $\lprime >0$, and $h \in L_2 ([0,1]^d)$, let $\|\langle h,\bpsi_{\lprime,\cdot}\rangle\|_2$ be the vector $L_2$ norm of the vector $\left(\langle h,\bpsi^v_{\lprime, m'}\rangle: m' \in \mathcal{Z}_{\lprime},v \in \left\{0, 1\right\}^d \setminus \{0\}^d\right)$. We have, 
		\be 
		\ & \P_{P,2}\left(2^{\lprime(\gamma+\frac{d}{2})}\|\langle \gtilde, \bpsi_{\lprime,\cdot}\rangle \|_{\infty}>C\right)
		\\&= \sum_{l=\lmin}^{\lmax}\P_{P,2}\left(2^{\lprime(\gamma+\frac{d}{2})}\|\langle \ghat_l, \bpsi_{\lprime,\cdot}\rangle \|_{\infty}>C, \lhat=l\right) \mathbbm{1} \left(\lprime \leq l\right)\\
		&= \sum_{l=\lmin}^{\lstar}\P_{P,2}\left(2^{\lprime(\gamma+\frac{d}{2})}\|\langle \ghat_l, \bpsi_{\lprime,\cdot}\rangle \|_{\infty}>C, \lhat=l\right)\mathbbm{1} \left(\lprime \leq l\right)\\&+\sum_{l=\lstar+1}^{\lmax}\P_{P,2}\left(2^{\lprime(\gamma+\frac{d}{2})}\|\langle \ghat_l, \bpsi_{\lprime,\cdot}\rangle \|_{\infty}>C, \lhat=l\right)\mathbbm{1} \left(\lprime \leq l\right)\\
		& \leq \sum_{l=\lmin}^{\lstar}\P_{P,2}\left(2^{\lprime(\gamma+\frac{d}{2})}\|\langle \ghat_l, \bpsi_{\lprime,\cdot}\rangle \|_{\infty}>C\right)\mathbbm{1} \left(\lprime \leq l\right)\\&+\sum_{l=\lstar+1}^{\lmax}\P_{P,2}\left( \lhat=l\right)\mathbbm{1} \left(\lprime \leq l\right)\\
		& \leq \sum_{l=\lmin}^{\lstar}\P_{P,2}\left(2^{\lprime(\gamma+\frac{d}{2})}\|\langle \ghat_l, \bpsi_{\lprime,\cdot}\rangle \|_{\infty}>C\right)\mathbbm{1} \left(\lprime \leq l\right)+\sum_{l>l^*}n^{-\eta},\label{eqn:ghat_fourier_lprime}
		\ee
		where the last inequality follows from \eqref{eqn:wrongchoicelepski_density} for some $\eta>3$ provided $C^*$ is chosen large enough as before. Now,
		\be 
		\ & \P_{P,2}\left(2^{\lprime(\gamma+\frac{d}{2})}\|\langle \ghat_l, \bpsi_{\lprime,\cdot}\rangle \|_{\infty}>C\right)\\
		&\leq \P_{P,2}\left(2^{\lprime(\gamma+\frac{d}{2})}\|\langle \ghat_l, \bpsi_{\lprime,\cdot}\rangle-\E_{P,2}\left(\langle \ghat_l, \bpsi_{\lprime,\cdot}\rangle\right) \|_{\infty}>C/2\right)\\&+ \mathbbm{1} \left(2^{\lprime(\gamma+\frac{d}{2})}\|\E_{P,2}\left(\langle \ghat_l, \bpsi_{\lprime,\cdot}\rangle\right) \|_{\infty}>C/2\right)\\
		&=\P_{P,2}\left(2^{\lprime(\gamma+\frac{d}{2})}\|\langle \ghat_l, \bpsi_{\lprime,\cdot}\rangle-\E_{P,2}\left(\langle \ghat_l, \bpsi_{\lprime,\cdot}\rangle\right) \|_{\infty}>C/2\right)
		\ee
		if $C > 2M$ (by definition in equation \eqref{eqn:besov_multidim}). Therefore, from \eqref{eqn:ghat_fourier_lprime}, one has for any $C > 2M$, 
		\be 
		\ &\P_{P,2}\left(2^{\lprime(\gamma+\frac{d}{2})}\|\langle \ghat, \bpsi_{\lprime,\cdot}\rangle \|_{\infty}>C\right)\\
		&\leq \sum_{l=\lmin}^{\lstar}\P_{P,2}\left(2^{\lprime(\gamma+\frac{d}{2})}\|\langle \ghat_l, \bpsi_{\lprime,\cdot}\rangle-\E_{P, 2}\left(\langle \ghat_l, \bpsi_{\lprime,\cdot}\rangle\right) \|_{\infty}>C/2\right)\mathbbm{1} \left(\lprime \leq l\right)\\
		&+\sum_{l=l^*}^{\lmax} n^{-3} \mathbbm{1} \left( \lprime \leq l \right). \label{eqn:ghat_fourier_lprime_second}
		\ee
		Considering the first term of the last summand of the above display, we have
		\be 
		\ & \sum_{l=\lmin}^{\lstar}\P_{P,2}\left(2^{\lprime(\gamma+\frac{d}{2})}\|\langle \ghat_l, \bpsi_{\lprime,\cdot}\rangle-\E_{P,2}\left(\langle \ghat_l, \bpsi_{\lprime,\cdot}\rangle\right) \|_{\infty}>C/2\right)\mathbbm{1} \left(\lprime \leq l\right)\\
		& \leq \sum_{l=\lmin}^{\lstar} \sum_{m \in \Z_{\lprime}}\sum_{v\in \{0,1\}^d}\P_{P,2}\left(\left|\frac{1}{n}\sum_{i=n+1}^{2n}\left(\psi_{\lprime, m}^v(\bX_i)-\E_{P,2}\left(\psi_{\lprime, m}^v(\bX_i)\right)\right)\right|>\frac{C/2}{2^{\lprime(\gamma+\frac{d}{2})}}\right)\mathbbm{1} \left(\lprime \leq l\right).
		\ee
		By Bernstein's inequality, for any $\lambda > 0$,
		\be 
		\ &\P_{P, 2} \left( \left| \frac{1}{n} \sum_{i = n + 1}^{2 n} \left( \psi_{\lprime, m}^v (\bX_i) - \E_{P, 2} \left( \psi_{\lprime, m}^v (\bX_i) \right) \right) \right| > \lambda \right) \\
		& \leq 2 \exp{\left(-\frac{n\lambda^2}{2\left(\sigma^2+\|\psi_{\lprime, m}^v\|_{\infty}\lambda/3\right)}\right)},
		\ee
		where $\sigma^2 = \E_{P, 2} \left( \psi_{\lprime, m}^v (\bX_i) - \E_{P, 2} \left( \psi_{\lprime, m}^v (\bX_i) \right) \right)^2$. Indeed, there exists a constant $C_4$ depending on $\bpsi_{0, 0}^0, \bpsi_{0, 0}^1, \gmax$ such that $\sigma^2 \leq C_4$ and $\|\psi_{\lprime, m}^v \|_{\infty} \leq C_4 2^{\frac{\lprime d}{2}}$. Therefore, 
		\be 
		\ & \sum_{l=\lmin}^{\lstar} \sum_{m \in \Z_{\lprime}}\sum_{v\in \{0,1\}^d}\P_{P,2}\left(\left|\frac{1}{n}\sum_{i=n+1}^{2n}\left(\psi_{\lprime, m}^v(\bX_i)-\E_{P,2}\left(\psi_{\lprime, m}^v(\bX_i)\right)\right)\right|>\frac{C/2}{2^{\lprime(\gamma+\frac{d}{2})}}\right)\mathbbm{1} \left(\lprime \leq l\right)\\
		&\leq 2\sum_{l=\lmin}^{\lstar} \sum_{m \in \Z_{\lprime}}\sum_{v\in \{0,1\}^d}\exp{\left(-\frac{C^2}{8C_4}\frac{n 2^{-2\lprime(\gamma+\frac{d}{2})}}{1+\frac{C}{2}2^{\frac{\lprime d}{2}}2^{-\lprime(\gamma+\frac{d}{2})}}\right)}\mathbbm{1} \left(\lprime \leq l\right)\\
		&=2 \sum_{l=\lmin}^{\lstar} \sum_{m \in \Z_{\lprime}}\sum_{v\in \{0,1\}^d}\exp{\left(-\frac{C^2}{8C_4}\frac{n 2^{-2 \lprime\gamma}}{2^{\lprime d}+\frac{C}{2}2^{\lprime(d-\gamma)}}\right)}\mathbbm{1} \left(\lprime \leq l\right)\\
		& \leq 2\sum_{l=\lmin}^{\lstar} \sum_{m \in \Z_{\lprime}}\sum_{v\in \{0,1\}^d}\exp{\left(-\frac{C^2}{8(1+\frac{C}{2})C_4}\frac{n 2^{-2\lprime\gamma}}{2^{\lprime d}}\right)}\mathbbm{1} \left(\lprime \leq l\right)\\
		&=2\sum_{l=\lmin}^{\lstar} \sum_{m \in \Z_{\lprime}}\sum_{v\in \{0,1\}^d}\exp{\left(-\frac{C^2}{8(1+\frac{C}{2})C_4}\frac{n 2^{-2\lstar\gamma}}{2^{\lstar d}\lstar d}2^{(\lstar-\lprime)(d+2\gamma)}\lstar d\right)}\mathbbm{1} \left(\lprime \leq l\right)\\
		& \leq 2\sum_{l=\lmin}^{\lstar} \sum_{m \in \Z_{\lprime}}\sum_{v\in \{0,1\}^d}\exp{\left(-\frac{ C^2C_2}{2^{d+3}C_1(1+\frac{C}{2})C_4}2^{(\lstar-\lprime)(d+2\gamma)}\lstar d\right)}\mathbbm{1} \left(\lprime \leq l\right) \label{eqn:lstar_use}\\
		& \leq 2\sum_{l=\lmin}^{\lstar} \sum_{m \in \Z_{\lprime}}\sum_{v\in \{0,1\}^d}\exp{\left(-\frac{ C^2C_2}{2^{d+3}C_1(1+\frac{C}{2})C_4}l d\right)}\mathbbm{1} \left(\lprime \leq l\right) \\
		& \leq 2\sum_{l=\lmin}^{\lstar} C(\psi_{0,0}^0,\psi_{0,0}^1)2^{\lprime d}\exp{\left(-\frac{ C^2C_2}{2^{d+3}C_1(1+\frac{C}{2})C_4}l d\right)}\mathbbm{1} \left(\lprime \leq l\right)\\
		& \leq 2\sum_{l=\lmin}^{\lstar} C(\psi_{0,0}^0,\psi_{0,0}^1)\exp{\left(-\left(\frac{ C^2C_2}{2^{d+3}C_1(1+\frac{C}{2})C_4}-1\right)l d\right)}\mathbbm{1} \left(\lprime \leq l\right)\\
		& \leq 2\lmax C(\psi_{0,0}^0,\psi_{0,0}^1)\exp{\left(-\left(\frac{ C^2C_2}{2^{d+3}C_1(1+\frac{C}{2})C_4}-1\right)\lmin d\right)}\mathbbm{1} \left(\lprime \leq \lmax\right),
		\ee
		if $\frac{ C^2C_2}{2^{d+3}C_1(1+\frac{C}{2})C_4}\geq 1$. The above inequality \eqref{eqn:lstar_use} uses the definition of $\lstar$. Indeed choosing $C$ large enough, one can guarantee,  $\left(\frac{ C^2C_2}{2^{d+3}C_1(1+\frac{C}{2})C_4}-1\right)\lmin d\geq 4\log{n}$. Such a choice of $C$ implies that, 
		\be 
		\ & \sum_{l=\lmin}^{\lstar} \sum_{m \in \Z_{\lprime}}\sum_{v\in \{0,1\}^d}\P_{P,2}\left(\left|\frac{1}{n}\sum_{i=n+1}^{2n}\left(\psi_{\lprime, m}^v(\bX_i)-\E_{P,2}\left(\psi_{\lprime, m}^v(\bX_i)\right)\right)\right|>\frac{C/2}{2^{\lprime(\gamma+\frac{d}{2})}}\right)\mathbbm{1} \left(\lprime \leq l\right)\\&\leq \frac{C(\psi_{0,0}^0,\psi_{0,0}^1)}{n^{3}}\mathbbm{1}(\lprime \leq \lmax),
		\ee
		which in turn implies that, for $C$ sufficiently large (depending on $M,\psi_{0,0}^0,\psi_{0,0}^1$) one has
		\be 
		\P_{P,2}\left(2^{\lprime(\gamma+\frac{d}{2})}\|\langle \ghat, \bpsi_{\lprime,\cdot}\rangle \|_2>C\right) &\leq  \frac{C(\psi_{0,0}^0,\psi_{0,0}^1)+1}{n^{3}} \mathbbm{1}(\lprime \leq \lmax).
		\ee
		This along with the logarithmic in $n$ size of $\lmax$ implies that for sufficiently large $n$, uniformly over $P \in \mathcal{P}(s,\gamma) $, $\gtilde$ belongs to $H(\gamma,C)$ with probability at least $1-1/n^2$, for a large enough constant $C$ depending on  $M,B_U$, and $\gmax$ (the choice of $\psi_{0,0}^0,\psi_{0,0}^1$ being fixed by specifying a regularity $S>\gmax$).
		
		However this $\gtilde$ does not satisfy the desired point-wise bounds. To achieve this let $\phi$ be a $C^{\infty}$ function such that $\psi(x) |_{[B_L, B_U]} \equiv x$ while $\frac{B_L}{2} \leq \psi(x) \leq 2B_U$ for all $x$. Finally, consider the estimator $\hat{g}(\bx) = \psi(\tilde{g}(\bx))$. We note that $|g(\bx) - \hat{g}(\bx)| \leq |g(\bx) - \tilde{g}(\bx)|$--- thus $\hat{g}$ is adaptive to the smoothness of the design density. The boundedness of the constructed estimator follows from the construction. Finally, we wish to show that almost surely, the constructed estimator belongs to the H\"{o}lder space with the same smoothness, possibly of a different radius. This is captured by the next lemma, the proof of which can be completed by following arguments similar to the proof of Lemma 3.1 in \cite{mukherjee2018optimal}. In particular, recall the definition of $H(s, M)$ in Section \ref{section_wavelets and function spaces}, 
		
		\begin{lem}
			\label{lemma:function_truncation}
			For all $h \in H(s, M)$, $\psi(h) \in H(s, C(M, s))$, where $C(M, s)$ is a universal constant dependent only on $M, s$ and independent of $h \in H(s, M)$. 
		\end{lem}
		
		The probabilistic tail bound \eqref{eq:tail_g} follows from essentially the same argument as the moment bound \eqref{eq:moment_g}. For the sake of completeness, we present the proof below. We first control the tail probability of $\Vert \tilde{g} - g \Vert_\infty$.
		\begin{align*}
		& \; \P_{P, 2} \left( \Vert \tilde{g} - g \Vert_\infty \ge (\tilde{C})^{\frac{d}{2 \gamma + d}} \left( \frac{n}{\log{n}} \right)^{- \frac{\gamma}{2 \gamma + d}} \right) \\
		= & \; \P_{P, 2} \left( \Vert \tilde{g} - g \Vert_\infty \ge (\tilde{C})^{\frac{d}{2 \gamma + d}} \left( \frac{n}{\log{n}} \right)^{- \frac{\gamma}{2 \gamma + d}}, \lhat \le \lstar \right) + \P_{P, 2} \left( \Vert \tilde{g} - g \Vert_\infty \ge (\tilde{C})^{\frac{d}{2 \gamma + d}} \left( \frac{n}{\log{n}} \right)^{- \frac{\gamma}{2 \gamma + d}}, \lhat > \lstar \right) \\
		\le & \; \P_{P, 2} \left( \Vert \tilde{g} - g \Vert_\infty \ge (\tilde{C})^{\frac{d}{2 \gamma + d}} \left( \frac{n}{\log{n}} \right)^{- \frac{\gamma}{2 \gamma + d}}, \lhat \le \lstar \right) + \P_{P, 2} \left( \lhat > \lstar \right).
		\end{align*}
		Recall that we have shown that, for some $\eta > 3$, $$\P_{P, 2} \left( \lhat > \lstar \right) = \sum_{l = \lstar + 1}^{\lmax} \P_{P, 2} \left( \lhat = l \right) \le \sum_{l = \lstar + 1}^{\lmax} n^{-\eta} \le \frac{\lmax}{n^\eta} \lesssim \frac{\log{n}}{n^3}.$$
		Thus we are left to control the following term: 
		\allowdisplaybreaks
		\begin{align*}
		& \; \P_{P, 2} \left( \Vert \tilde{g} - g \Vert_\infty \ge (\tilde{C})^{\frac{d}{2 \gamma + d}} \left( \frac{n}{\log{n}} \right)^{- \frac{\gamma}{2 \gamma + d}}, \lhat \le \lstar \right) \\
		\le & \; \P_{P, 2} \left( \Vert \tilde{g} - \ghat_{\lstar} \Vert_\infty + \Vert \ghat_{\lstar} - g \Vert_\infty \ge (\tilde{C})^{\frac{d}{2 \gamma + d}} \left( \frac{n}{\log{n}} \right)^{- \frac{\gamma}{2 \gamma + d}}, \lhat \le \lstar \right) \\
		\le & \; \P_{P, 2} \left( \Vert \ghat_{\lstar} - g \Vert_\infty \ge (\tilde{C})^{\frac{d}{2 \gamma + d}} \left( \frac{n}{\log{n}} \right)^{- \frac{\gamma}{2 \gamma + d}} - C^{*} \sqrt{\frac{2^{\lstar d} \lstar d}{n}} \right) \\
		\le & \; \P_{P, 2} \left( \Vert \ghat_{\lstar} - \Pi_{P, 2} [ g \vert V_{\lstar} ] \Vert_\infty + \Vert \Pi_{P, 2} [ g \vert V_{\lstar} ] - g \Vert_\infty \ge (\tilde{C})^{\frac{d}{2 \gamma + d}} \left( \frac{n}{\log{n}} \right)^{- \frac{\gamma}{2 \gamma + d}} - C^{*} \sqrt{\frac{2^{\lstar d} \lstar d}{n}} \right) \\
		\le & \; \P_{P, 2} \left( \Vert \ghat_{\lstar} - \Pi_{P, 2} [ g \vert V_{\lstar} ] \Vert_\infty \ge (\tilde{C})^{\frac{d}{2 \gamma + d}} \left( \frac{n}{\log{n}} \right)^{- \frac{\gamma}{2 \gamma + d}} - C^{*} \sqrt{\frac{2^{\lstar d} \lstar d}{n}} - C_1 M 2^{ - \lstar \gamma } \right) \\
		\le & \; \P_{P, 2} \left( \Vert \ghat_{\lstar} - \Pi_{P, 2} [ g \vert V_{\lstar} ] \Vert_\infty \ge (\tilde{C})^{\frac{d}{2 \gamma + d}} \left( \frac{n}{\log{n}} \right)^{- \frac{\gamma}{2 \gamma + d}} - (C^{*} + C_2) \sqrt{\frac{2^{\lstar d} \lstar d}{n}} \right) \\
		\le & \; \P_{P, 2} \left( \Vert \ghat_{\lstar} - \Pi_{P, 2} [ g \vert V_{\lstar} ] \Vert_\infty \ge \left\{ \left( C \right)^{\frac{d}{2 \gamma + d}} - C^{*} - C_2 \right\} \sqrt{\frac{2^{\lstar d} \lstar d}{n}} \right) \\
		\le & \; \exp \{ - C' \lstar d \} \le n^{- \eta}
		\end{align*}
		for some $\eta > 3$, where the second and fifth inequalities follow from the definitions of $\lhat$ and $\lstar$ respectively, the fourth inequality follows from \eqref{eqn:bias_lepski_density}, the sixth inequality is a consequence of \eqref{eqn:lstar_size} and the seventh inequality is implied by Lemma \ref{lemma_linear_projection_tail_bound}. 
		
		For verification purposes, \eqref{eqn:lstar_size} implies
		\begin{align*}
		& \; \left\{ 2^{-d} \left( \frac{C_1 M}{C_2} \right)^{- \frac{2d}{2 \gamma + d}} \left( \frac{n}{\log{n}} \right)^{- \frac{d}{2 \gamma + d}} \right\}^{\frac{\gamma}{d}} \ge \left( 2^{- \lstar d} \right)^{\frac{\gamma}{d}} \\
		\Rightarrow & \; C_1 M \left( \frac{C_1 M}{C_2} \right)^{- \frac{2 \gamma}{2 \gamma + d}} \left( \frac{n}{\log{n}} \right)^{- \frac{\gamma}{2 \gamma + d}} \ge C_1 M 2^{- (\lstar - 1) \gamma} \ge C_2 \sqrt{\frac{2^{(\lstar - 1) d} (\lstar - 1) d}{n}} \\
		\Rightarrow & \; (\tilde{C})^{\frac{d}{2 \gamma + d}} \left( \frac{n}{\log{n}} \right)^{- \frac{\gamma}{2 \gamma + d}} \ge \left( \frac{\tilde{C} C_2}{C_1 M} \right)^{\frac{d}{2 \gamma + d}} 2^{- \frac{d}{2}} \sqrt{\frac{\lstar - 1}{\lstar}} \sqrt{\frac{2^{\lstar d} \lstar d}{n}} \ge C^{\frac{d}{2 \gamma + d}} \sqrt{\frac{2^{\lstar d} \lstar d}{n}}.
		\end{align*} 
		Thus combining the above calculations, we have
		\begin{align*}
		\P_{P, 2} \left( \Vert \gtilde - g \Vert_\infty \ge (\tilde{C})^{\frac{d}{2 \gamma + d}} \left( \frac{n}{\log{n}} \right)^{- \frac{\gamma}{2 \gamma + d}} \right) \lesssim \frac{\log{n}}{n^3} + \frac{1}{n^3} \lesssim \frac{\log{n}}{n^3}.
		\end{align*}
		Since $\vert g(\bx) - \ghat(\bx) \vert \le \vert g(\bx) - \gtilde(\bx) \vert$, the above display immediately implies
		\begin{align*}
		\P_{P, 2} \left( \Vert \ghat - g \Vert_\infty \ge (\tilde{C})^{\frac{d}{2 \gamma + d}} \left( \frac{n}{\log{n}} \right)^{- \frac{\gamma}{2 \gamma + d}} \right) \lesssim \frac{\log{n}}{n^3}.
		\end{align*}
		
		Now, the construction of $\hat{f}$ satisfying the desired properties of Theorem \ref{thm_adaptive_density_regression} can be done following ideas from the proof of Theorem 1.1 of \cite{mukherjee2018optimal}. In particular, construct the estimator $\hat{g}$ of the design density $g$ as above from the second part of the sample and let for $j \in \mathcal{T}_1$, $\hat{f}_j(\bx)=\frac{1}{n}\sum\limits_{i=1}^n \frac{W_i}{\ghat(\bX_i)}K_{V_j}\left(\bX_i, \bx\right)$. 
		Now, let
		\be
		\jhat=\min\left\{j\in \mathcal{T}_1: \ \|\fhat_j-\fhat_{j'}\|_{\infty}\leq C^{**}\sqrt{\frac{2^{j'd}j'd}{n}}, \ \forall j' \in \mathcal{T}_1 \ \text{s.t.} \ j' \geq j \right\}.
		\ee
		where $C^{**}$ depends only on the known parameters of the problem and can be determined from the proof hereafter. Thereafter, consider the estimator $\ftilde \coloneqq \fhat_{\jhat}$.
		Now define
		\be
		\jstar & \coloneqq \min \left\{j \in \mathcal{T}_1: 2^{-j d\frac{s}{d}}\leq \sqrt{\frac{2^{jd} jd}{n}}\right\},
		\ee
		Therefore
		\be 
		\E_{P}\|\ftilde-f\|_{\infty}&\leq\E_{P}\|\fhat_{\jhat}-f\|_{\infty}\mathbbm{1}(\jhat\leq \jstar)+\E_{P}\|\fhat_{\jhat}-f\|_{\infty}\mathbbm{1}(\jhat> \jstar).\label{eqn:fhat_basic_inequality}
		\ee
		Thereafter using Lemma \ref{lemma_linear_projection_tail_bound} and equation \eqref{eqn:holder_approx} we have
		\be 
		\ & \E_{P}\|\fhat_{\jhat}-f\|_{\infty}\mathbbm{1}(\jhat\leq \jstar)\\
		&\leq \E_{P}\|\fhat_{\jhat}-\fhat_{\jstar}\|_{\infty}\mathbbm{1}(\jhat\leq \jstar)+\E_{P}\|\fhat_{\jstar}-f\|_{\infty}\\
		&\leq C^{**}\sqrt{\frac{2^{\jstar d}\jstar d}{n}}+\E_{P,2}\|\fhat_{\jstar}-\E_{P,1}(\fhat_{\jstar})\|_{\infty}+\E_{P,2}\|\E_{P,1}(\fhat_{\jstar})-f\|_{\infty}\\
		&\leq (C^{**}+C(B_U,B_L,\bpsi_{0,0}^0,\bpsi_{0,0}^1))\sqrt{\frac{2^{\jstar d}\jstar d}{n}}\\&+\E_{P,2}\|\Pi(f(\frac{g}{\ghat}-1)|V_{\jstar})\|_{\infty}+\|f-\Pi(f|V_{\jstar})\|_{\infty}\\ 
		& \leq (C^{**}+C(B_U,B_L,\bpsi_{0,0}^0,\bpsi_{0,0}^1))\sqrt{\frac{2^{\jstar d}\jstar d}{n}}\\&+C(M,\bpsi_{0,0}^0,\bpsi_{0,0}^1)2^{-\jstar s}+\E_{P,2}\|\Pi(f(\frac{g}{\ghat}-1)|V_{\jstar})\|_{\infty}.\\ \label{eqn:fhat_jhatlessthanjstar}
		\ee
		Now, by standard computations involving compactly supported wavelet bases and properties of $\ghat$
		\be 
		\E_{P,2}\|\Pi(f(\frac{g}{\ghat}-1)|V_{\jstar})\|_{\infty}&\leq C(\bpsi_{0,0}^0,\bpsi_{0,0}^1)\E_{P,2}\|f(\frac{g}{\ghat}-1)\|_{\infty}\\
		&\leq C(B_U,B_L,\bpsi_{0,0}^0,\bpsi_{0,0}^1)\E_{P,2}\|\ghat-g\|_{\infty}\\
		&\leq C(B_U,B_L,M,\gmax,\bpsi_{0,0}^0,\bpsi_{0,0}^1)\left(\frac{n}{\log{n}}\right)^{-\frac{\gamma}{2\gamma+d}}.\label{eqn:fhat_gestimationeffect}
		\ee
		Combining \eqref{eqn:fhat_jhatlessthanjstar}, \eqref{eqn:fhat_gestimationeffect}, definition of $\jstar$, and the fact that $\gamma>s$, we  have
		\be 
		\E_{P}\|\fhat_{\jhat}-f\|_{\infty}\mathbbm{1}(\jhat\leq \jstar)&\leq C(B_U,B_L,M,\gmax,\bpsi_{0,0}^0,\bpsi_{0,0}^1)\left(\frac{n}{\log{n}}\right)^{-\frac{s}{2s+d}}, \label{eqn:lepski_rightchocie_f}
		\ee
		provided that $C^**$ is chosen depending only the known parameters of the problem. Now using arguments similar to those leading to \eqref{eqn:control_II_density_step1} we have
		\be 
		\E_{P}\|\fhat_{\jhat}-f\|_{\infty}\mathbbm{1}(\jhat> \jstar)&\leq C(B_U,B_L,\bpsi_{0,0}^0,\bpsi_{0,0}^1)\sum_{j>\jstar}2^{jd}\P_{P}(\jhat=j).\label{eqn:lepski_wrongchocie_basic}
		\ee
		We now complete the control over II by suitably bounding $\P_{P}(\jhat=j)$. To this end, note that for any $j >\jstar$,
		\be
		\ & \P_{P}(\jhat=j) \\& \leq \sum_{j>\jstar}\mathbb{P}_{P}\left(\|\fhat_j-\fhat_{\jstar}\|_{\infty}> C^{**}\sqrt{\frac{2^{jd}jd}{n}}\right)\\
		&\leq \sum_{j>\jstar}\E_{P,2}\left\{\begin{array}{c}\mathbb{P}_{P,1}\left(\|\fhat_{\jstar}-\E_{P,1}\left(\fhat_{\jstar}\right)\|_{\infty}> \frac{C^{**}}{2}\sqrt{\frac{2^{jd}jd}{n}}-\|\E_{P,1}\left(\fhat_{\jstar}\right)-\E_{P,1}\left(\fhat_{j}\right)\|_{\infty}\right)\\+\mathbb{P}_{P,1}\left(\|\fhat_j-\E_{P,1} \left(\fhat_{j}\right)\|_{\infty}> \frac{C^{**}}{2}\sqrt{\frac{2^{jd}jd}{n}}\right)\end{array}\right\}\\
		& \leq \sum_{j>\jstar}E_{P,2}\left\{ \begin{array}{c}\mathbb{P}_{P,1}\left(\|\fhat_{\jstar}-\E_{P,1}\left(\fhat_{\jstar}\right)\|_{\infty}> \frac{C^{**}}{2}\sqrt{\frac{2^{jd}jd}{n}}-\|\Pi\left(f\frac{g}{\ghat}|V_{\jstar}\right)-\Pi\left(f\frac{g}{\ghat}|V_{j}\right)\|_{\infty}\right)\\+\mathbb{P}_{P,1}\left(\|\fhat_j-\E_{P,1} \left(\fhat_{j}\right)\|_{\infty}> \frac{C^{**}}{2}\sqrt{\frac{2^{jd}jd}{n}}\right)\end{array}\right\}.\\
		\ee
		Now, 
		\be 
		\ &\left\|\Pi\left(f\frac{g}{\ghat}|V_{\jstar}\right)-\Pi\left(f\frac{g}{\ghat}|V_{j}\right)\right\|_{\infty}\\
		&\leq C(M,\bpsi_{0,0}^0,\bpsi_{0,0}^1)2^{-\jstar s}+C(B_U,B_L,,\bpsi_{0,0}^0,\bpsi_{0,0}^1)\|\ghat-g\|_{\infty}.
		\ee	
		Using the fact that $\sqrt{\frac{2^{jd}jd}{n}}>\sqrt{\frac{2^{\jstar d}\jstar d}{n}}$ for $j>\jstar$, we have, using the definition of $\jstar$, that there exist $C,C'>0$ depending on $M,B_U,B_L,\bpsi_{0,0}^0,\bpsi_{0,0}^1$ such that
		\be 
		\ & \P_{P}(\jhat=j) \\
		&\leq \sum_{j>\jstar} \E_{P,2}\left\{ \begin{array}{c}\mathbb{P}_{P,1}\left(\|\fhat_{\jstar}-\E_{P,1}\left(\fhat_{\jstar}\right)\|_{\infty}> (\frac{C^{**}}{2}-C)\sqrt{\frac{2^{j d}j d}{n}}\right)\\+\mathbb{P}_{P,1}\left(\|\fhat_j-\E_{P,1} \left(\fhat_{j}\right)\|_{\infty}> \frac{C^{**}}{2}\sqrt{\frac{2^{jd}jd}{n}}\right)+\mathbb{P}_{P,2}\left(\|\ghat-g\|_{\infty}> C'\sqrt{\frac{2^{\jstar d}\jstar d}{n}}\right)\end{array}\right\}.\\ \label{eqn:lepski_wrong_choice_f}
		\ee	
		Now, provided $C^{**}>2C$ is chosen large enough (depending on $B_U,B_L,\bpsi_{0,0}^0,\bpsi_{0,0}^1$) we have there exists large enough $C''$ (depending on $B_U,B_L,\bpsi_{0,0}^0,\bpsi_{0,0}^1$) such that
		\be 
		\ & \mathbb{P}_{P,1}\left(\|\fhat_{\jstar}-\E_{P,1}\left(\fhat_{\jstar}\right)\|_{\infty}> (\frac{C^{**}}{2}-C)\sqrt{\frac{2^{j d}j d}{n}}\right)\\&+\mathbb{P}_{P,1}\left(\|\fhat_j-\E_{P,1} \left(\fhat_{j}\right)\|_{\infty}> \frac{C^{**}}{2}\sqrt{\frac{2^{jd}jd}{n}}\right)\leq 2e^{-C''j d}. \label{eqn:lepski_wrongchoice_fhat_conc}
		\ee
		Henceforth, whenever required, $C, C',C''$ will be chosen to be large enough depending on the known parameters of the problem, which in turn will imply that $C^{**}$ can be chosen large enough depending on the known parameters of the problem as well. First note that, the last term in the above display can be bounded rather crudely using the following lemma.
		\begin{lem}\label{lemma_ghat_prob_bound}
			Assume $\gamma_{\min}>s_{\max}$. Then for $C', C_1,C_2>0$ (chosen large enough depending on $B_U, \bpsi_{0,0}^0,\bpsi_{0,0}^1$) one has
			\be
			\sup_{P \in \Par(s, \gamma)}\P_{P,2}\Big(\|\ghat-g\|_{\infty}>C'\sqrt{\frac{2^{j^*d}\jstar d}{n}}\Big)& \leq C_1(\lmax-\lmin)e^{-C_2\lmin  d}.
			\ee
		\end{lem}
		The proof of Lemma \ref{lemma_ghat_prob_bound} can be argued as follows. Indeed, $\ghat=\psi(\gtilde)$, where $\psi(u)$ is a $C^{\infty}$ function which is identically equal to $u$ on $[B_L, B_U]$ and has universally bounded first derivative. Therefore, it is enough to prove Lemma \ref{lemma_ghat_prob_bound} for $\gtilde$ instead of $\ghat$ and thereby invoking a simple first order Taylor series argument along with the fact that $\psi(g)\equiv g$ owing to the bounds on $g$. 
		The crux of the argument for proving Lemma \ref{lemma_ghat_prob_bound}  is that by Lemma \ref{lemma_linear_projection_tail_bound}, any $\ghat_l$ for $l \in \mathcal{T}_2$ suitably concentrates around $g$ in a radius of the order of $\sqrt{\frac{2^{ld}ld}{n}}$. The proof of the lemma is therefore very similar to the proof of adaptivity of $\ghat$ (by dividing into cases where the chosen $\lhat$ is larger and smaller than $\lstar$ respectively and thereafter invoking Lemma \ref{lemma_linear_projection_tail_bound}) and therefore we omit the details.
		
		Plugging in the result of Lemma \ref{lemma_ghat_prob_bound} into \eqref{eqn:lepski_wrong_choice_f}, and thereafter using the facts that $\gmin>s_{\max}$, $\lmax,\jmax$ are both poly logarithmic in nature,  along with equations \eqref{eqn:fhat_basic_inequality}, \eqref{eqn:lepski_rightchocie_f}, \eqref{eqn:lepski_wrongchocie_basic}, and \eqref{eqn:lepski_wrongchoice_fhat_conc}, we have the existence of an estimator $\ftilde$ depending on $M,B_U,B_L,s_{\min},s_{\max},\gmax$, such that for every $(s,\gamma) \in [s_{\min},s_{\max}]\times [\gmin,\gmax]$,
		$$\sup\limits_{P \in \mathcal{P}(s,\gamma)} \E_P\|\ftilde-f\|_{\infty} \leq C\left(\frac{n}{\log{n}}\right)^{-\frac{s}{2s+d}},$$ 
		with a large enough positive constant $C$ depending on  $M,B_U,B_L,s_{\min}, \gmax$, $\bpsi_{0,0}^0,\bpsi_{0,0}^1$. The proof that $\ftilde \in H(s, C)$ with probability at least $1 - 1 / n^2$ follows largely from the same argument used to show $\gtilde \in H(\gamma, C)$, by the boundedness of $Y$ and $\ghat$ along with the fact that $\ghat \in H(\gamma, C)$ with high probability with $\gamma>s$.
		
		However this $\ftilde$ does not satisfy the desired point-wise bounds. To achieve this, as before, let $\phi$ be a $C^{\infty}$ function such that $\psi(x) |_{[B_L, B_U]} \equiv x$ while $\frac{B_L}{2} \leq \psi(x) \leq 2B_U$ for all $x$. Finally, consider the estimator $\hat{f}(\bx) = \psi(\tilde{f}(\bx))$. We note that $|f(\bx) - \hat{f}(\bx)| \leq |f(\bx) - \tilde{f}(\bx)|$--- thus $\hat{f}$ is adaptive to the smoothness of the design density. The boundedness of the constructed estimator follows from the construction. The proof of the fact that the constructed estimator belongs to the H\"{o}lder space with the same smoothness, possibly of a different radius follows once again from of Lemma \ref{lemma:function_truncation}.
		
		Finally, we obtain the tail bound \eqref{eq:tail_f} of $\Vert \fhat - f \Vert_\infty$. As in the proof of \eqref{eq:tail_g}, we consider the tail bound of $\Vert \ftilde - f \Vert_\infty$ first, and then the desired tail bound of $\Vert \fhat - f \Vert_\infty$ follows directly from the inequality $\vert f(\bx) - \fhat(\bx) \vert \leq \vert f(\bx) - \ftilde(\bx) \vert$ almost surely.
		\begin{align*}
		& \; \P_P \left( \Vert \ftilde - f \Vert_\infty \ge (\tilde{C}^{\dag})^{\frac{d}{2 s + d}} \left( \frac{n}{\log{n}} \right)^{- \frac{s}{2 s + d}} \right) \\
		= & \; \P_P \left( \Vert \ftilde - f \Vert_\infty \ge (\tilde{C}^{\dag})^{\frac{d}{2 s + d}} \left( \frac{n}{\log{n}} \right)^{- \frac{s}{2 s + d}}, \jhat \le \jstar \right) + \P_P \left( \Vert \ftilde - f \Vert_\infty \ge (\tilde{C}^{\dag})^{\frac{d}{2 s + d}} \left( \frac{n}{\log{n}} \right)^{- \frac{s}{2 s + d}}, \jhat > \jstar \right) \\
		\le & \; \P_P \left( \Vert \ftilde - f \Vert_\infty \ge (\tilde{C}^{\dag})^{\frac{d}{2 s + d}} \left( \frac{n}{\log{n}} \right)^{- \frac{s}{2 s + d}}, \jhat \le \jstar \right) + \P_P \left( \jhat > \jstar \right).
		\end{align*}
		As was shown before, for some $\eta > 3$, we have
		$$
		\P_P \left( \jhat > \jstar \right) = \sum_{j = \jstar + 1}^{\jmax} \P_P \left( \jhat = j \right) \le \sum_{j = \jstar + 1}^{\jmax} C' \exp \{ - C'' j d \} \lesssim \frac{\jmax}{n^\eta} \lesssim \frac{\log{n}}{n^3}.
		$$
		Therefore we are left to bound the first term $\P_P \left( \Vert \ftilde - f \Vert_\infty \ge (\tilde{C}^{\dag})^{\frac{d}{2 s + d}} \left( \frac{n}{\log{n}} \right)^{- \frac{s}{2 s + d}}, \jhat \le \jstar \right)$, proceeded as below. First, we write
		\begin{align*}
		& \; \P_P \left( \Vert \ftilde - f \Vert_\infty \ge (\tilde{C}^{\dag})^{\frac{d}{2 s + d}} \left( \frac{n}{\log{n}} \right)^{- \frac{s}{2 s + d}}, \jhat \le \jstar \right) \\
		= & \; \E_{P, 2} \left[ \P_{P, 1} \left( \Vert \ftilde - f \Vert_\infty \ge (\tilde{C}^{\dag})^{\frac{d}{2 s + d}} \left( \frac{n}{\log{n}} \right)^{- \frac{s}{2 s + d}}, \jhat \le \jstar \right) \right] \\
		= & \; \E_{P, 2} \left[ \P_{P, 1} \left( \Vert \ftilde - f \Vert_\infty \ge (\tilde{C}^{\dag})^{\frac{d}{2 s + d}} \left( \frac{n}{\log{n}} \right)^{- \frac{s}{2 s + d}}, \jhat \le \jstar \right) \mathbbm{1} \left\{ \Vert \ghat - g \Vert_\infty \le C' \left( \frac{n}{\log{n}} \right)^{- \frac{s}{2 s + d}} \right\} \right] \\
		& \; + \E_{P, 2} \left[ \P_{P, 1} \left( \Vert \ftilde - f \Vert_\infty \ge (\tilde{C}^{\dag})^{\frac{d}{2 s + d}} \left( \frac{n}{\log{n}} \right)^{- \frac{s}{2 s + d}}, \jhat \le \jstar \right) \mathbbm{1} \left\{ \Vert \ghat - g \Vert_\infty > C' \left( \frac{n}{\log{n}} \right)^{- \frac{s}{2 s + d}} \right\} \right] \\
		\le & \; \E_{P, 2} \left[ \P_{P, 1} \left( \Vert \ftilde - f \Vert_\infty \ge (\tilde{C}^{\dag})^{\frac{d}{2 s + d}} \left( \frac{n}{\log{n}} \right)^{- \frac{s}{2 s + d}}, \jhat \le \jstar \right) \mathbbm{1} \left\{ \Vert \ghat - g \Vert_\infty \le C' \left( \frac{n}{\log{n}} \right)^{- \frac{s}{2 s + d}} \right\} \right] \\
		& \; + \P_{P, 2} \left( \Vert \ghat - g \Vert_\infty > C' \left( \frac{n}{\log{n}} \right)^{- \frac{s}{2 s + d}} \right) \\
		\le & \; \E_{P, 2} \left[ \P_{P, 1} \left( \Vert \ftilde - f \Vert_\infty \ge (\tilde{C}^{\dag})^{\frac{d}{2 s + d}} \left( \frac{n}{\log{n}} \right)^{- \frac{s}{2 s + d}}, \jhat \le \jstar \right) \mathbbm{1} \left\{ \Vert \ghat - g \Vert_\infty \le C' \left( \frac{n}{\log{n}} \right)^{- \frac{s}{2 s + d}} \right\} \right] \\
		& \; + \frac{\log{n}}{n^3}
		\end{align*}
		where the last inequality is implied by \eqref{eq:tail_g}. Then within the event $\Vert \ghat - g \Vert_\infty \le C' \left( \frac{n}{\log{n}} \right)^{- \frac{s}{2 s + d}}$, we have
		\allowdisplaybreaks
		\begin{align*}
		& \; \P_{P, 1} \left( \Vert \ftilde - f \Vert_\infty \ge (\tilde{C}^{\dag})^{\frac{d}{2 s + d}} \left( \frac{n}{\log{n}} \right)^{- \frac{s}{2 s + d}}, \jhat \le \jstar \right) \\
		\le & \; \P_{P, 1} \left( \Vert \ftilde - \fhat_{\jstar} \Vert_\infty + \Vert \fhat_{\jstar} - f \Vert_\infty \ge (\tilde{C}^{\dag})^{\frac{d}{2 s + d}} \left( \frac{n}{\log{n}} \right)^{- \frac{s}{2 s + d}}, \jhat \le \jstar \right) \\
		\le & \; \P_{P, 1} \left( \Vert \fhat_{\jstar} - f \Vert_\infty \ge (\tilde{C}^{\dag})^{\frac{d}{2 s + d}} \left( \frac{n}{\log{n}} \right)^{- \frac{s}{2 s + d}} - C^{**} \sqrt{\frac{2^{\jstar d} \jstar d}{n}} \right) \\
		\le & \; \P_{P, 1} \left( \Vert \fhat_{\jstar} - \E_{P, 1} \left( \fhat_{\jstar} \right) \Vert_\infty + \Vert \E_{P, 1} \left( \fhat_{\jstar} \right) - f \Vert_\infty \ge (\tilde{C}^{\dag})^{\frac{d}{2 s + d}} \left( \frac{n}{\log{n}} \right)^{- \frac{s}{2 s + d}} - C^{**} \sqrt{\frac{2^{\jstar d} \jstar d}{n}} \right) \\
		\le & \; \P_{P, 1} \left( \begin{array}{c}
		\Vert \fhat_{\jstar} - \E_{P, 1} \left( \fhat_{\jstar} \right) \Vert_\infty + \Vert \Pi \left( f \left( \frac{g}{\ghat} - 1 \right) \vert V_{\jstar} \right) \Vert_\infty \geq \\
		(\tilde{C}^{\dag})^{\frac{d}{2 s + d}} \left( \frac{n}{\log{n}} \right)^{- \frac{s}{2 s + d}} - C^{**} \sqrt{\frac{2^{\jstar d} \jstar d}{n}} - C( M, \bpsi_{0, 0}^0, \bpsi_{0, 0}^1 ) 2^{- \jstar s} 
		\end{array} \right) \\
		\le & \; \P_{P, 1} \left( \begin{array}{c}
		\Vert \fhat_{\jstar} - \E_{P, 1} \left( \fhat_{\jstar} \right) \Vert_\infty + \Vert \Pi \left( f \left( \frac{g}{\ghat} - 1 \right) \vert V_{\jstar} \right) \Vert_\infty \geq \\
		(\tilde{C}^{\dag})^{\frac{d}{2 s + d}} \left( \frac{n}{\log{n}} \right)^{- \frac{s}{2 s + d}} - (C^{**} + C( M, \bpsi_{0, 0}^0, \bpsi_{0, 0}^1 )) \sqrt{\frac{2^{\jstar d} \jstar d}{n}} 
		\end{array} \right) \\
		\le & \; \P_{P, 1} \left( \begin{array}{c}
		\Vert \fhat_{\jstar} - \E_{P, 1} \left( \fhat_{\jstar} \right) \Vert_\infty + C ( B_U, B_L, \bpsi_{0, 0}^0, \bpsi_{0, 0}^1 ) \Vert \ghat - g \Vert_\infty \geq \\
		(\tilde{C}^{\dag})^{\frac{d}{2 s + d}} \left( \frac{n}{\log{n}} \right)^{- \frac{s}{2 s + d}} - (C^{**} + C( M, \bpsi_{0, 0}^0, \bpsi_{0, 0}^1 )) \sqrt{\frac{2^{\jstar d} \jstar d}{n}} 
		\end{array} \right) \\
		\le & \; \P_{P, 1} \left( \begin{array}{c}
		\Vert \fhat_{\jstar} - \E_{P, 1} \left( \fhat_{\jstar} \right) \Vert_\infty \geq \\
		\left( (\tilde{C}^{\dag})^{\frac{d}{2 s + d}} - C' ( B_U, B_L, \bpsi_{0, 0}^0, \bpsi_{0, 0}^1 ) \right) \left( \frac{n}{\log{n}} \right)^{- \frac{s}{2 s + d}} - (C^{**} + C( M, \bpsi_{0, 0}^0, \bpsi_{0, 0}^1 ) ) \sqrt{\frac{2^{\jstar d} \jstar d}{n}}
		\end{array} \right) \\
		\le & \; \P_{P, 1} \left(	\Vert \fhat_{\jstar} - \E_{P, 1} \left( \fhat_{\jstar} \right) \Vert_\infty \geq C^{***} \sqrt{\frac{2^{\jstar d} \jstar d}{n}} \right) \\
		\le & \; C_1 e^{- C_2 \jstar d} \lesssim \frac{1}{n^\eta}
		\end{align*}
		for some $\eta > 3$, where the second inequality is due to the definition of $\jhat$, the fourth inequality follows from the property of wavelet basis of Cohen-Daubechies-Vial type, the fifth inequality is implied by the definition of $\jstar$, the sixth inequality follows from the properties of $\ghat, f, g$, the seventh inequality is a result of being restricted to the event $\Vert \ghat - g \Vert_\infty \le C' \left( \frac{n}{\log{n}} \right)^{- \frac{s}{2 s + d}}$, the eighth inequality is again implied by the definition of $\jstar$, and the ninth inequality follows from the same argument as in \eqref{eqn:lepski_wrongchoice_fhat_conc}.
		
		Finally combining the above analysis, and $\jstar < \lmin$, for sufficiently large $C_1', C_2' > 0$, we have
		\begin{align*}
		\P_P \left( \Vert \ftilde - f \Vert_\infty \ge (\tilde{C}^{\dag})^{\frac{d}{2 s + d}} \left( \frac{n}{\log{n}} \right)^{- \frac{s}{2 s + d}} \right) \lesssim \frac{\log{n}}{n^3}.
		\end{align*}
	\end{proof}
	
	\begin{rem} \label{rem:mar}
	For the mean response functional in missing data models studied in Section \ref{sec:mar}, denote $W = A Y$. There we need to estimate $b (\bx) = \E_{P} [Y | A = 1, \bX = \bx]$ instead of $f (\bx) = \E_{P} [W | \bX = \bx]$. Note that $b (\bx)$ can be rewritten as:
	$$
	b (\bx) = \int y f (y | a = 1, \bx) d y = \int \frac{w}{g_{1} (\bx)} f (y, a, \bx) d y d a
	$$
	where $g_{1} (\bx) = f(\bx | A = 1) \BP(A = 1)$, where $f (\bx | A = 1)$ is the conditional density of $\bX$ given $A = 1$. We define an adaptive estimator $\hat{b}$ of $b$ in the same way as we define $\hat{f}$ as follows. Define $\Pi (b \vert V_{j}) (\bx) \coloneqq \E_{P} \left[ \frac{W}{g_{1} (\bX)} K_{V_{j}} (\bX, \bx) \right]$, and $\hat{b}_{j} (\bx) \coloneqq \frac{1}{n} \sum_{i = 1}^{n} \frac{W_{i}}{\hat{g}_{1} (\bX_{i})} K_{V_{j}} (\bX_{i}, \bx)$, where $\hat{g}_{1} (\bx)$ is estimated in the same way as $\hat{g} (\bx)$ except that $\hat{g}_{j} (\bx)$ is replaced by $\hat{g}_{1, j} (\bX = \bx) = \frac{1}{n} \sum\limits_{i = n + 1}^{2n} A_{i} K_{V_{j}} (\bX_{i}, x)$. Then the above proof goes through if we replace all the corresponding $f$ by $b$. 
	\end{rem}
	
	\subsection*{Proof of Theorem \ref{thm_treatment_effect}}
	\begin{proof}\hspace*{\fill}
		
		\textit{\textbf{(i) Proof of Upper Bound}} \\
		
		The general scheme of the proof involves identifying a non-adaptive minimax estimator of $\phi(P)$ under the knowledge of $P \in \Par_{(\a, \b, \g)}$, demonstrating suitable conditional bias and variance properties of this sequence of estimators, and thereafter invoking Theorem \ref{theorem_general_lepski} to conclude. This routine can be carried out as follows. Suppose that we observe $n$ i.i.d. samples $\{ Y_i, A_i, \bX_i \}_{i = 1}^n$. First divide the samples into two disjoint parts $\mathcal{D} = \mathcal{D}_1 \coprod \mathcal{D}_2$ with sample sizes $n_1 = n ( 1 - 1 / \log{n} )$ and $n_2 = n / \log{n}$ respectively. Further divide the subsample $\mathcal{D}_2$ into two disjoint and equal-sized parts $\mathcal{D}_2 = \mathcal{D}_{21} \coprod \mathcal{D}_{22}$ with sample sizes $n_{21} = n_{22} = n / (2 \log{n})$, where we estimate $g$ by $\ghat$ adaptively from $\mathcal{D}_{22}$ (as in Theorem \ref{thm_adaptive_density_regression}), and estimate $a$ and $b$ by $\ahat$ and $\bhat$ respectively, adaptively from $\mathcal{D}_{21}$ (as in Theorem \ref{thm_adaptive_density_regression}). Let $\E_{P, S}$ denote the expectation while samples with indices not in $S$ held fixed, for $S \subset \{1, 21, 22\}$. A first order influence function for $\phi(P)$ at $P$ is given by $(Y - b(\bX)) (A - a(\bX)) - \phi(P)$ and a resulting first order estimator for $\phi(P)$ is $\frac{1}{n} \sum_{i=1}^n (Y_i - \bhat(\bX_i)) (A_i - \ahat(\bX_i))$, computed from the subsample $\mathcal{D}_1$. This estimator has a conditional bias $\int \left( (b(\bx) - \bhat(\bx)) (a(\bx) - \ahat(\bx)) \right) g(\bx) d \bx$. Indeed for $\bav < \frac{d}{2}$, this bias turns out to be suboptimal compared to the minimax rate of convergence of $n^{- \frac{4 \alpha + 4 \beta}{2 \alpha + 2 \beta + d}}$ in mean squared loss. The most intuitive way to proceed is to estimate and correct for the bias. If there exists a ``Dirac-kernel" $K(\bx_1,\bx_2) \in L_2\left([0,1]^d\times [0,1]^d\right)$ such that $\int h(\bx_1)K(\bx_1,\bx_2)d\bx_2=h(\bx_1)$ almost surely $\bx_1$ for all $h \in L_2 ([0,1]^d)$, then one can estimate the bias term by $\frac{1}{n (n-1)} \sum\limits_{1 \leq i_1 \neq i_2 \leq n} \frac{(Y_{i_1} - \bhat(\bX_{i_1}))}{\sqrt{g(\bX_{i_1})}} K(\bX_{i_1}, \bX_{i_2}) \frac{(A_{i_2} - \ahat({\bX_{i_2}}))}{\sqrt{g(\bX_{i_2})}}$, provided the marginal density $g$ was known. Indeed there are two concerns with the above suggestion. The first one is the knowledge of $g$. This can be relatively easy to deal with by plugging in an suitable estimate $\ghat$--although there are some subtleties involved (refer to Section \ref{section_discussion} for more on this). The primary concern though is the non-existence of a ``Dirac-kernel" of the above sort as an element of $L_2 ([0,1]^d) \times L_2 ([0,1]^d)$. This necessitates the following modification where one works with projection kernels on suitable finite-dimensional linear subspace $L$ of $L_2 ([0,1]^d)$ which guarantees existence of such kernels when the domain space is restricted to $L$. In particular, we work with the linear subspace $V_k$ (defined in \ref{section_wavelets and function spaces}) where the choice of $k$ is guided by the balance between the bias and variance properties of the resulting estimator. In particular, a choice of $k$ is guided by the knowledge of the parameter space $\Par_{(\a,\b, \g)}$. For any $k \in \mathcal{K}$, this implies that our bias corrected second-order estimator of $\phi(P)$ is given by 
		\be
		\phihat_{n, k}&=\frac{1}{n}\sum_{i=1}^n {(Y_i-\bhat(\bX_i))(A_i-\ahat(\bX_i))}\\&-\frac{1}{n(n-1)}\sum\limits_{1\leq i_1\neq i_2\leq  n}S\left(\frac{(Y_{i_1}-\bhat(\bX_{i_1})))}{\sqrt{\ghat(\bX_{i_1})}}K_{V_k}(\bX_{i_1},\bX_{i_2})\frac{(A_{i_2}-\ahat({\bX_{i_2}}))}{\sqrt{\ghat(\bX_{i_2})}}\right).
		\ee
		Note that the division by $\ghat$ is permitted by the properties guaranteed by Theorem \ref{thm_adaptive_density_regression}.
		Indeed this sequence of estimators is in the form of those considered by Theorem \ref{theorem_general_lepski} with $$\Ltilde_1(O)=(Y-\bhat(\bX))(A-\ahat(\bX)),$$
		$$\Ltilde_{2l}(O)=\frac{(Y-\bhat(\bX)))}{\sqrt{\ghat(\bX)}},\quad \Ltilde_{2r}(O)=\frac{(A-\ahat(\bX)))}{\sqrt{\ghat(\bX)}},$$
		where by Theorem \ref{thm_adaptive_density_regression} $\max\{|\Ltilde_1(O)|, |\Ltilde_{2l}(O)|, |\Ltilde_{2r}(O)|\}\leq C(B_L, B_U)$. Therefore it remains to show that the sequence $\phihat_{n, k}$ satisfies the bias and variance properties (A) and (B) necessary for application of Theorem \ref{theorem_general_lepski}. \\
		
		We first verify the conditional bias property. Utilizing the representation of the first order bias as stated above, we have
		\be 
		\ & | \E_{P, 1} \left( \phihat_{n, k} - \phi(P) \right) | \\
		& = \left| \begin{array}{c} 
			\int \left( ( b(\bx) - \bhat(\bx) ) ( a (\bx) - \ahat(\bx) ) \right) g(\bx) d \bx \\
			- \E_{P, 1} \left[ S \left( \frac{(Y_{1} - \bhat(\bX_{1}))}{\sqrt{\ghat(\bX_{1})}} K_{V_k} (\bX_{1}, \bX_{2}) \frac{(A_{2} - \ahat({\bX_{2}}))}{\sqrt{\ghat(\bX_{2})}} \right) \right]
		\end{array} \right|.
		\label{eqn:treatmenteffect_second_orderbias_initial}
		\ee
		Now, using the notation $\delb(\bx) = b(\bx) - \bhat(\bx)$ and $\dela(\bx) = a(\bx) - \ahat(\bx)$, we have
		\be 
		\ & \E_{P, 1} \left[\frac{( Y_{1} - \bhat(\bX_{1}) )}{\sqrt{\ghat(\bX_{1})}} K_{V_k} (\bX_{1}, \bX_{2}) \frac{( A_{2} - \ahat({\bX_{2}}) )}{\sqrt{\ghat(\bX_{2})}} \right] \\
		& = \int \int \frac{\delb(\bx_1) g(\bx_1)}{\sqrt{\ghat(\bx_{1})}} K_{V_k} (\bx_{1}, \bx_{2}) \frac{\dela(\bx_2) g(\bx_2)}{\sqrt{\ghat(\bx_{2})}} d\bx_1 d\bx_2 \\
		& = \int \frac{\delb(\bx_1) g(\bx_1)}{\sqrt{\ghat(\bx_{1})}} \Pi \left( \frac{\dela g}{\sqrt{\ghat}} | V_k \right)(\bx_1) d \bx_1\\
		& = \int \frac{\delb(\bx_1) g(\bx_1)}{\sqrt{\ghat(\bx_{1})}} \frac{\dela(\bx_1) g(\bx_1)}{\sqrt{\ghat(\bx_1)}} d \bx_1\\
		& - \int \frac{\delb(\bx_1) g(\bx_1)}{\sqrt{\ghat(\bx_{1})}} \Pi \left( \frac{\dela g}{\sqrt{\ghat}} | V_k^{\perp} \right) (\bx_1) d \bx_1 \\
		& = \int \left( (b(\bx) - \bhat(\bx)) (a(\bx) - \ahat(\bx)) \right) g(\bx) d \bx \\
		& + \int \dela(\bx_1) \delb(\bx_1) g^2(\bx_1) \left( \frac{1}{\ghat(\bx_1)} - \frac{1}{g(\bx_1)} \right) d \bx_1 \\
		& - \int \frac{\delb(\bx_1) g(\bx_1)}{\sqrt{\ghat(\bx_{1})}} \Pi \left( \frac{\dela g}{\sqrt{\ghat}} | V_k^{\perp} \right) (\bx_1) d \bx_1. \label{eqn:treatmenteffect_second_orderbias_secondterm}
		\ee
		Plugging \eqref{eqn:treatmenteffect_second_orderbias_secondterm} into \eqref{eqn:treatmenteffect_second_orderbias_initial}, we get,
		\be 
		\ & | \E_{P, 1} \left( \phihat_{n, k} - \phi(P) \right) | \\
		& = \left| \begin{array}{c}
			\int \dela(\bx_1) \delb(\bx_1) g^2(\bx_1) \left( \frac{1}{\ghat(\bx_1)} - \frac{1}{g(\bx_1)} \right) d \bx_1 \\
			- \int \frac{\delb(\bx_1) g(\bx_1)}{\sqrt{\ghat(\bx_{1})}} \Pi \left(\frac{\dela g}{\sqrt{\ghat}} | V_k^{\perp} \right) (\bx_1) d \bx_1
		\end{array} \right|.
		\label{eqn:treatmenteffect_remainder_bias}
		\ee
		Now, by repeatedly applying Cauchy-Schwarz inequality and invoking results in Theorem \ref{thm_adaptive_density_regression}, we have\\
		\be 
		\ & \left| \int \dela(\bx_1) \delb(\bx_1) g^2(\bx_1) \left( \frac{1}{\ghat(\bx_1)} - \frac{1}{g(\bx_1)} \right) d\bx_1 \right| \\
		& \leq \left( \int \frac{g^4(\bx_1)}{g^2(\bx_1) \ghat^2(\bx_1)} \left( \ghat(\bx_1) - g(\bx_1) \right)^2 d \bx_1 \right)^{\frac{1}{2}} \\
		& \times \left( \int \left( \ahat(\bx_1) - a(\bx_1) \right)^4 d \bx_1 \right)^{\frac{1}{4}} \left( \int \left( \bhat(\bx_1) - b(\bx_1) \right)^4 d \bx_1 \right)^{\frac{1}{4}} \\
		& \leq \frac{B_U^2}{B_L^2} \left\| \ghat - g \right\|_2 \left\| \ahat - a \right\|_4 \left\| \bhat - b \right\|_4 \\
		& \leq \frac{B_U^2}{B_L^2} C^{\frac{d}{2 \gamma + d} + \frac{d}{2 \alpha + d} + \frac{d}{2 \beta + d}} \left( \frac{n_2 / 2}{\log{(n_2 / 2)}} \right)^{- \frac{\alpha}{2 \alpha + d} - \frac{\beta}{2 \beta + d} - \frac{\gamma}{2 \gamma + d}}. \label{eqn:treatmentefftect_bias_part1}
		\ee
		Moreover, when restricted to the ``good'' event 
		\begin{align*}
		\mathcal{I}_{2; b, a, g} (n_2) \coloneqq \left\{ O \in \mathcal{D}_2: \begin{array}{c}
		\Vert \hat{b} - b \Vert_\infty \leq C_b^{\frac{d}{2 \beta + d}} \left( \frac{n_{2} / 2}{\log({n_{2} / 2})} \right)^{- \frac{\beta}{2 \beta + d}}
		\cap \; \Vert \hat{a} - a \Vert_\infty \leq C_a^{\frac{d}{2 \alpha + d}} \left( \frac{n_{2} / 2}{\log{(n_{2} / 2 )}} \right)^{- \frac{\alpha}{2 \alpha + d}} \cap \\
		\Vert \hat{g} - g \Vert_\infty \leq \tilde{C}^{\frac{d}{2 \alpha + d}} \left( \frac{n_{2} / 2}{\log{(n_{2} / 2 )}} \right)^{- \frac{\gamma}{2 \gamma + d}} \cap \hat{g} \in H(\gamma, C) \cap \hat{b} \in H(\beta, C) \cap \hat{a} \in H(\alpha, C)
		\end{array} \right\}
		\end{align*}
		with probability at least $1 - \frac{C \log{n_2}}{n_2^\eta}$ for some absolute constant $C$ and $\eta > 2$ by Theorem \ref{thm_adaptive_density_regression},
		\be 
		\ & \left| \int \frac{\delb(\bx_1) g(\bx_1)}{\sqrt{\ghat(\bx_{1})}} \Pi \left( \frac{\dela g}{\sqrt{\ghat}} | V_k^{\perp} \right)(\bx_1) d\bx_1 \right| \\
		& = \left| \int \Pi \left( \frac{\delb g}{\sqrt{\ghat}} | V_k^{\perp} \right)(\bx_1) \Pi \left( \frac{\dela g}{\sqrt{\ghat}} | V_k^{\perp} \right) (\bx_1) d \bx_1 \right| \\
		& \leq \left\| \Pi \left( \frac{\delb g}{\sqrt{\ghat}} | V_k^{\perp} \right) \right\|_2 \left\| \Pi \left( \frac{\dela g}{\sqrt{\ghat}} | V_k^{\perp} \right) \right\|_2 \\
		& \leq C k^{ - \frac{\beta + \alpha}{d} }, \label{eqn:treatmentefftect_bias_part2}
		\ee
		where the last line follows for some constant $C$ (depending on $M, B_U, B_L,\gmax$) by Theorem \ref{thm_adaptive_density_regression}, definition in equation \eqref{eqn:besov_multidim}, and noting that $\|\Pi\left(h|V_j\right)\|_{\infty}\leq C(B_U)$ if $\|h\|_{\infty}\leq B_U$. Therefore, one has combining \eqref{eqn:treatmenteffect_remainder_bias}, \eqref{eqn:treatmentefftect_bias_part1}, and \eqref{eqn:treatmentefftect_bias_part2}, that for a constant $C$ (depending on $M, B_U, B_L, \gmin, \gmax$) and $\gmin(\epsilon) \coloneqq \frac{\gmin}{1 + \epsilon}$.
		\be 
		\ & \left| \E_{P, 1} \left( \phihat_{n, k} - \phi(P) \right) \right| \\
		& \leq C \left[ \left( \frac{n_2 / 2}{\log{(n_2 / 2)}} \right)^{- \frac{\alpha}{2 \alpha + d} - \frac{\beta}{2 \beta + d} - \frac{\gamma}{2 \gamma + d}} + k^{- \frac{\beta + \alpha}{d}} \right] \\
		& \leq C \left[ (n_2 / 2)^{- \frac{\alpha}{2 \alpha + d} - \frac{\beta}{2 \beta + d} - \frac{\gmin(\epsilon)}{2 \gmin(\epsilon) + d}} (n_2 / 2)^{\frac{\gmin(\epsilon)}{2 \gmin(\epsilon) + d} - \frac{\gamma}{2 \gamma + d}} \log{(n_2 / 2)}^{\frac{\alpha}{2 \alpha + d} + \frac{\beta}{2 \beta + d} + \frac{\gamma}{2 \gamma + d}} + k^{- \frac{\beta + \alpha}{d}} \right] \\
		& \leq C \left[ (n_2 / 2)^{- \frac{\alpha}{2 \alpha + d} - \frac{\beta}{2 \beta + d} - \frac{\gmin(\epsilon)}{2 \gmin(\epsilon) + d}} (n_2 / 2)^{- \frac{\gamma - \gmin(\epsilon)}{2 \gamma + d}} \log^3{(n_2 / 2)} + k^{- \frac{\beta + \alpha}{d}} \right] \\
		& \leq C \left[ (n_2 / 2)^{- \frac{\alpha}{2 \alpha + d} - \frac{\beta}{2 \beta + d} - \frac{\gmin(\epsilon)}{2 \gmin(\epsilon) + d}} (n_2 / 2)^{- \frac{\epsilon\gmin}{(1 + \epsilon) (2 \gmax + d)}} \log^3{(n_2 / 2)} + k^{- \frac{\beta + \alpha}{d}} \right] \\
		& \leq C \left[ n_2^{- \frac{\alpha}{2 \alpha + d} - \frac{\beta}{2 \beta + d} - \frac{\gmin(\epsilon)}{2 \gmin(\epsilon) + d}} n_2^{- \frac{\epsilon\gmin}{(1 + \epsilon) (2 \gmax + d)}} \log^3{n_2} + k^{- \frac{\beta + \alpha}{d}} \right] \\
		& \leq C \left[ (n / \log{n})^{- \frac{\alpha}{2 \alpha + d} - \frac{\beta}{2 \beta + d} - \frac{\gmin(\epsilon)}{2 \gmin(\epsilon) + d}} (n / \log{n})^{- \frac{\epsilon\gmin}{(1 + \epsilon) (2 \gmax + d)}} \log^3{(n / \log{n})} + k^{- \frac{\beta + \alpha}{d}} \right] \\
		& \leq C \left[ n^{- \frac{\alpha}{2 \alpha + d} - \frac{\beta}{2 \beta + d} - \frac{\gmin(\epsilon)}{2 \gmin(\epsilon) + d}} + k^{- \frac{\beta + \alpha}{d}} \right],
		\ee
		where the constant $C$ changes its value from line to line. Now, letting $\theta = (\a, \b, \g)$, $f_1(\theta) = \bav$ and $f_2(\theta) = - \frac{\alpha}{2 \alpha + d} - \frac{\beta}{2 \beta + d} - \frac{\gmin(\epsilon)}{2 \gmin(\epsilon) + d}$ we have that the conditional bias property (A) corresponding to Theorem \ref{theorem_general_lepski} holds with the given choice of $f_1$ and $f_2$ and a constant $C$ depending on $M, B_U, B_L, \gmax$ for $\{P \in \Par_{\theta}: f_1(\theta) = \bav, f_2(\theta) > \frac{2 \alpha + 2 \beta}{2 \alpha + 2 \beta + d} \}$.
		
		The proof of the validity of the conditional variance property corresponding to Theorem \ref{theorem_general_lepski} is easy to derive by standard Hoeffding decomposition of $\phihat_{n, k}$ followed by applications of moment bounds in Lemmas \ref{lemma_ustat_tail_use} and \ref{lemma_linear_tail_bound}. For calculations of similar flavor, refer to proof of Theorem 1.3 in \cite{mukherjee2018optimal}. Note that this is the step where we have used the fact that $\bav \leq \frac{d}{4}$, since otherwise the linear term dominates resulting in $O(\frac{1}{n})$ the variance. Then invoking Theorem \ref{theorem_general_lepski}, we have
		\be 
		\sup_{P \in \Par_{\theta}: \atop f_1(\theta) = \frac{\alpha + \beta}{2}, f_2(\theta) > \frac{2 \alpha + 2 \beta}{2 \alpha + 2 \beta + d}} \E_P \left(\phihat_{n, \khat} - \phi(P) \right)^2 \leq C \left( \frac{\sqrt{\log{n}}}{n} \right)^{\frac{4 \alpha + 4 \beta}{d + 2 \alpha + 2 \beta}}.
		\ee
		
		Noting that for $\theta \in \Theta$, since $\gmin>2(1+\epsilon)\max\{\alpha,\beta\}$, one has automatically, $f_2(\theta)>\frac{2\alpha+2\beta}{2\alpha+2\beta+d}$, which completes the proof of the upper bound.
		
		\textit{\textbf{(ii) Proof of Lower Bound}} \\
		To prove a lower bound matching the upper bound above, note that $\phi(P) = \E_P \left( \cov_P \left(Y, A | \bX \right) \right) = \E_P \left( A Y \right) - \E_P \left( a(\bX) b(\bX) \right)$. Indeed, $\E_P \left( AY \right)$ can be estimated at $\sqrt{n}$-rate by sample average of $A_i Y_i$. Therefore, it suffices to prove a lower bound for adaptive estimation of $\E_P \left( a(\bX) b(\bX) \right)$ (our proof continues to hold for $\E_P \left( \cov_P \left(Y, A | \bX \right) \right)$ since the perturbations created for the lower bound of $\E_P \left( A Y \right)$ is trivial and we only present the lower bound for estimation of the most important term of its estimation). Let $c(\bX) = \E_P \left( Y | A = 1, \bX \right) - \E_P \left( Y | A = 0, \bX \right)$, which implies owing to the binary nature of $A$ that $\E_P \left( Y | A, \bX \right) = c(\bX) \left( A - a(\bX) \right) + b(\bX)$. For the purpose of lower bound it is convenient to parametrize the data generating mechanism by $(a, b, c, g)$, which implies that $\phi(P) = \int a(\bx) b(\bx) g(\bx) d \bx$. With this parametrization, we show that the same lower bound holds in a smaller class of problems where $g \equiv 1$ on $[0, 1]^d$. Specifically consider
		\be 
		\Theta_{sub} = \left\{ \begin{array}{c}
			P = (a, b, c, g): \\ 
			a \in H (\alpha, M), b \in H(\beta, M),\ \frac{\alpha + \beta}{2} < \frac{d}{4}, \\ g \equiv 1, (a(\bx), b(\bx)) \in [B_L, B_U]^2 \ \forall \ \bx \in [0, 1]^d 
		\end{array} \right\}.
		\ee
		The likelihood of $O \sim P$ for $P \in \Theta_{sub}$ can then be written as
		\be 
		\begin{array}{c}
			a(\bX)^A (1 - a(\bX))^{1 - A} \\
			\times \left( c(\bX) (1 - a(\bX)) + b(\bX) \right)^{YA} \left( 1 - c(\bX) (1 - a(\bX)) - b(\bX) \right)^{(1-Y) A} \\ 
			\times \left( - c(\bX) a(\bX) + b(\bX) \right)^{Y (1 - A)} \left( 1 + c(\bX) a(\bX) - b(\bX) \right)^{(1 - Y) (1 - A)}. \label{eqn:likelihood_treatment_effect}
		\end{array}
		\ee
		Let for some $(\alpha, \beta, \gamma)$ tuple in the original problem $\Theta$, one has $$\sup_{P \in \Par_{(\alpha, \beta, \gamma)}} \E_P \left( \phihat - \phi(P) \right)^2 \leq C \left( \frac{\sqrt{\log{n}}}{n} \right)^{\frac{4 \alpha + 4 \beta}{d + 2 \alpha + 2 \beta}}.$$
		
		Now, let $H: [0, 1]^d \rightarrow \mathbb{R}$ be a $C^{\infty}$ function supported on $\left[ 0, \frac{1}{2} \right]^d$ such that $\int H (\bx) d \bx = 0$ and $\int H^2 (\bx) d \bx = 1$ and let for $k \in \mathbb{N}$ (to be decided later) $\Omega_1, \ldots, \Omega_k$ be the translations of the cube $k^{-\frac{1}{d}} \left[ 0, \frac{1}{2} \right]^d$ that are disjoint and contained in $[0, 1]^d$. Let  $\bx_1, \ldots, \bx_k$ denote the bottom left corners of these cubes.
		
		Assume first that $\alpha < \beta$. We set for $\lambda = ( \lambda_1, \ldots, \lambda_k ) \in \{ -1, +1 \}^k$ and $\alpha \leq  \beta' < \beta$,
		$$a_{\lambda} (\bx) = \frac{1}{2} + \Big(\frac{1}{k}\Big)^{\frac{\alpha}{d}} \sum\limits_{j = 1}^k \lambda_j H \Big( (\bx - \bx_j) k^{\frac{1}{d}} \Big),$$
		
		$$b_{\lambda} (\bx) = \frac{1}{2} + \Big(\frac{1}{k}\Big)^{\frac{\beta'}{d}} \sum\limits_{j = 1}^k \lambda_j H \Big( (\bx - \bx_j) k^{\frac{1}{d}} \Big),$$
		
		$$c_{\lambda} (\bx) = \frac{\frac{1}{2} - b_{\lambda} (\bx)}{1 - a_{\lambda} (\bx)}.$$
		A properly chosen $H$ guarantees $a_{\lambda} \in H(\alpha, M)$ and $b_{\lambda} \in H(\beta', M)$ for all $\lambda$. Let 
		\be 
		\Theta_0 & = \left\{ P^n: P = \left( a_{\lambda}, \frac{1}{2}, 0, 1 \right), \lambda \in \{-1,+1\}^k \right\},
		\ee
		and 
		\be 
		\Theta_1 & = \left\{ P^n: P = (a_{\lambda}, b_{\lambda}, c_{\lambda}, 1), \lambda \in \{ -1, +1 \}^k \right\}.
		\ee
		Finally let $\Theta_{test} = \Theta_0 \cup \Theta_1$. Let $\pi_0$ and $\pi_1$ be uniform priors on $\Theta_0$ and $\Theta_1$ respectively. It is easy to check that by our choice of $H$, $\phi(P) = \frac{1}{4}$ on $\Theta_0$ and $\phi(P) = \frac{1}{4} + \left(\frac{1}{k}\right)^{\frac{\alpha + \beta'}{d}}$ for $P \in \Theta_1$. Therefore, using notation from Lemma \ref{lemma:prop_cai}, $\mu_1 = \frac{1}{4}$, $\mu_2 = \frac{1}{4} + \left(\frac{1}{k}\right)^{\frac{\alpha + \beta'}{d}}$, and $\sigma_1 = \sigma_2 = 0$. Since $\Theta_0 \subseteq P(\alpha, \beta, \gamma)$, we must have that worst case error of estimation over $\Theta_0$ is bounded by $C \left( \frac{\sqrt{\log{n}}}{n} \right)^{\frac{4 \alpha + 4 \beta}{d + 2 \alpha + 2 \beta}}$. Therefore, the $\pi_0$ average bias over $\Theta_0$ is also bounded by $C \left( \frac{\sqrt{\log{n}}}{n} \right)^{\frac{2 \alpha + 2 \beta}{d + 2 \alpha + 2 \beta}}$. This implies by Lemma \ref{lemma:prop_cai}, that the $\pi_1$ average bias over $\Theta_1$ (and hence the worst case bias over $\Theta_1$) is bounded below by 
		\be 
		\left(\frac{1}{k}\right)^{\frac{\alpha + \beta'}{d}} - C \left(\frac{\sqrt{\log{n}}}{n}\right)^{\frac{2 \alpha + 2 \beta}{d + 2 \alpha + 2 \beta}} - C \left(\frac{\sqrt{\log{n}}}{n}\right)^{\frac{2 \alpha + 2 \beta}{d + 2 \alpha + 2 \beta}} \eta, \label{eqn:bias_lower_bound_treatment_effect}
		\ee
		where $\eta$ is the chi-square divergence between the probability measures $\int P^n d \pi_0(P^n)$ and $\int P^n d \pi_1(P^n)$. We now bound $\eta$ using Proposition \ref{prop_chisquare_affinity}. 
		
		To put ourselves in the notation of Proposition \ref{prop_chisquare_affinity}, let for $\lambda \in \{-1, +1\}^k$, $P_{\lambda}$ and $Q_{\lambda}$ be the probability measures identified from $\Theta_0$ and $\Theta_1$ respectively. 
		
		Therefore, with $\chi_j = \{0, 1\} \times \{0, 1\} \times \Omega_j$, we indeed have for all $j = 1, \ldots, k$, $P_{\lambda} (\chi_j) = Q_{\lambda} (\chi_j) = p_j$ where there exists a constant $c$ such that $p_j = \frac{c}{k}$. 
		
		Letting $\pi$ be the uniform prior over $\{-1, +1\}^k$ it is immediate that $\eta = \chi^2 \left( \int P_{\lambda} d ( \pi(\lambda) ),\int Q_{\lambda} d ( \pi(\lambda) ) \right)$. 
		
		It now follows by calculations similar to the proof of Theorem 4.1 in \cite{robins2009semiparametric}, that for a constant $C'>0$
		
		\be 
		\chi^2\left(\int P_{\lambda}d(\pi(\lambda)),\int Q_{\lambda}d(\pi(\lambda))\right)& \leq \exp\left(C'\frac{n^2}{k}\left(k^{-\frac{4\beta'}{d}}+k^{-4\frac{\alpha+\beta'}{2d}}\right)\right)-1.
		\ee
		Now choosing $k=\left(\frac{n}{\sqrt{c_*\log{n}}}\right)^{\frac{2d}{d+2\alpha+2\beta'}}$, we have 
		\be 
		\chi^2\left(\int P_{\lambda}d(\pi(\lambda)),\int Q_{\lambda}d(\pi(\lambda))\right)& \leq n^{2C'c_*}-1.
		\ee
		Therefore choosing $c_*$ such that $2C'c_*+\frac{2\alpha+2\beta'}{2\alpha+2\beta'+d}<\frac{2\alpha+2\beta}{2\alpha+2\beta+d}$, we have the desired result by \eqref{eqn:bias_lower_bound_treatment_effect}. The proof for $\alpha> \beta$ is similar after changing various quantities to:
		$$a_{\lambda}(\bx)=\frac{1}{2}+\Big(\frac{1}{k}\Big)^{\frac{\alpha'}{d}}\sum\limits_{j=1}^k\lambda_j H\Big((\bx-\bx_j)k^{\frac{1}{d}}\Big), \quad \beta\leq \alpha'<\alpha,$$
		
		$$b_{\lambda}(\bx)=\frac{1}{2}+\Big(\frac{1}{k}\Big)^{\frac{\beta}{d}}\sum\limits_{j=1}^k\lambda_j H\Big((\bx-\bx_j)k^{\frac{1}{d}}\Big),$$
		\be 
		c_{\lambda}(\bX)=\frac{(\frac{1}{2}-a_{\lambda}(\bX))b_{\lambda}(\bX)}{a_{\lambda}(\bX)(1-a_{\lambda}(\bX))},
		\ee
		\be 
		\Theta_0&=\left\{P^n: P=\left(\frac{1}{2},b_{\lambda},0,1\right): \lambda \in \{-1,+1\}^k\right\},
		\ee
		and 
		\be 
		\Theta_1&=\left\{P^n: P=(a_{\lambda},b_{\lambda},c_{\lambda},1): \lambda \in \{-1,+1\}^k\right\}.
		\ee
		For the case of $\alpha=\beta$, choose $\alpha'<\alpha$ and therefore, $\alpha'<\beta$ and thereafter work with
		$$a_{\lambda}(\bx)=\frac{1}{2}+\Big(\frac{1}{k}\Big)^{\frac{\alpha'}{d}}\sum\limits_{j=1}^k\lambda_j H\Big((\bx-\bx_j)k^{\frac{1}{d}}\Big),$$
		
		$$b_{\lambda}(\bx)=\frac{1}{2}+\Big(\frac{1}{k}\Big)^{\frac{\beta}{d}}\sum\limits_{j=1}^k\lambda_j H\Big((\bx-\bx_j)k^{\frac{1}{d}}\Big),$$
		
		$$c_{\lambda}(\bx)=\frac{\frac{1}{2}-b_{\lambda}(\bx)}{1-a_{\lambda}(\bx)}.$$
		
		\be 
		\Theta_0&=\left\{P^n: P=\left(a_{\lambda},\frac{1}{2},0,1\right): \lambda \in \{-1,+1\}^k\right\},
		\ee
		and 
		\be 
		\Theta_1&=\left\{P^n: P=(a_{\lambda},b_{\lambda},c_{\lambda},1): \lambda \in \{-1,+1\}^k\right\}.
		\ee
		This completes the proof of the lower bound.
	\end{proof}
	
	\subsection*{Proof of Theorem \ref{thm_missing_data}}
	\begin{proof}
		\textit{\textbf{(i) Proof of Upper Bound}} \\
		The general scheme of the proof is the same as that of Theorem \ref{thm_treatment_effect} and involves identifying a non-adaptive minimax estimator of $\phi(P)$ under the knowledge of $P \in \Par_{(\a,\b,\g)}$, demonstrating suitable bias and variance properties of this sequence of estimators, and thereafter invoking Theorem \ref{theorem_general_lepski} to conclude. This routine can be carried out as follows. We observe $n$ i.i.d. samples $\{Y_iA_i,A_i,\bX_i\}_{i=1}^n$. First divide the samples into two disjoint parts $\mathcal{D} = \mathcal{D}_1 \coprod \mathcal{D}_2$ with sample sizes $n_1 = n ( 1 - 1 / \log{n} )$ and $n_2 = n / \log{n}$ respectively. Further divide the subsample $\mathcal{D}_2$ into two disjoint and equal-sized parts $\mathcal{D}_2 = \mathcal{D}_{21} \coprod \mathcal{D}_{22}$ with sample sizes $n_{21} = n_{22} = n / (2 \log{n})$. We estimate $\E\left(A|\bx\right)$ and $b$ by $\widehat{\E} \left(A|\bx\right)$ and $\bhat(\bx) \coloneqq \widehat{\E} \left(Y|A=1,\bx\right)$ respectively, adaptively from $\mathcal{D}_{21}$ (as in Theorem \ref{thm_adaptive_density_regression}). Note that $g(\bX)=f(\bX|A=1)P(A=1)$. Therefore, also estimate $\P_P ( A = 1)$ by $\hat{\pi} \coloneqq \frac{1}{n_{22}}\sum_{i \in \mathcal{D}_{22}} A_i$ i.e. the sample average of $A$'s from $\mathcal{D}_{22}$ and $\fhat_1$ is estimated as an estimator of $f(X|A=1)$ also from $\mathcal{D}_{22}$ using density estimation technique among observations with $A=1$. Finally, our estimates of $a$ and $g$ are $\ahat(\bx)=\frac{1}{\widehat{\E}\left(A|\bx\right)}$ and $\ghat = \fhat_1 \hat{\pi}$ respectively. In the following, we will freely use Theorem \ref{thm_adaptive_density_regression}, for desired properties of $\ahat,\bhat$, and $\ghat$. In particular, following the proof of Theorem \ref{thm_adaptive_density_regression}, we can actually assume that our choice of $\ghat$ also satisfies the necessary conditions of boundedness away from $0$ and $\infty$, as well as membership in $H(\gamma, C)$ with high probability for a large enough $C>0$.  A first order influence function for $\phi(P)$ at $P$ is given by $A a(\bX) (Y - b(\bX)) + b(\bX) - \phi(P)$ and a resulting first order estimator for $\phi(P)$ is $\frac{1}{n_1} \sum_{i \in \mathcal{D}_1} {A_i \ahat(\bX_i) (Y_i - \bhat(\bX_i)) + \bhat(\bX_i)}$. This estimator has a bias $- \E_{P, 2} \int\left((b(\bx)-\bhat(\bx))(a(\bx)-\ahat(\bx))\right)g(\bx)d\bx$. Indeed for $\bav<\frac{d}{2}$, this bias turns out to be suboptimal compared to the minimax rate of convergence of $n^{-\frac{4\alpha+4\beta}{2\alpha+2\beta+d}}$ in mean squared loss. Similar to the proof of Theorem \ref{thm_treatment_effect} we use a second-order bias corrected estimator as follows. 
		
		Once again we work with the linear subspace $V_k$ (defined in \ref{section_wavelets and function spaces}) where the choice of $k$ is guided by the balance between the conditional bias and variance properties of the resulting estimator. In particular, a choice of $k$ is guided by the knowledge of the parameter space $\Par_{(\a,\b, \g)}$. For any $k \in \mathcal{K}$, our bias corrected second-order estimator of $\phi(P)$ is given by 
		\be
		\phihat_{n, k} & = \frac{1}{n_1} \sum_{i = 1}^{n_1} {A_i \ahat(\bX_i) (Y_i - \bhat(\bX_i)) + \bhat(\bX_i)} \\
		& + \frac{1}{n_1(n_1 - 1)} \sum\limits_{1\leq i_1\neq i_2\leq  n_1}S\left(\frac{A_{i_1}(Y_{i_1}-\bhat(\bX_{i_1})))}{\sqrt{\ghat(\bX_{i_1})}}K_{V_k}(\bX_{i_1},\bX_{i_2})\frac{(A_{i_2}\ahat({\bX_{i_2}})-1)}{\sqrt{\ghat(\bX_{i_2})}}\right)
		\ee
		Note that division by $\ghat$ is permitted by the properties guaranteed by Theorem \ref{thm_adaptive_density_regression}.
		Indeed this sequence of estimators is in the form of those considered by Theorem \ref{theorem_general_lepski} with $$\Ltilde_1(O)=A\ahat(\bX)(Y-\bhat(\bX))+\bhat(\bX),$$
		$$\Ltilde_{2l}(O)=-\frac{A(Y-\bhat(\bX)))}{\sqrt{\ghat(\bX)}}, \quad \Ltilde_{2r}(O)=\frac{(A\ahat(\bX)-1)}{\sqrt{\ghat(\bX)}},$$
		where by Theorem \ref{thm_adaptive_density_regression} $\max\{|\Ltilde_1(O)|, |\Ltilde_{2l}(O)|, |\Ltilde_{2r}(O)|\}\leq C(B_L, B_U)$. Therefore it remains to show that the sequence $\phihat_{n, k}$ satisfies the conditional bias and variance properties (A) and (B) necessary for application of Theorem \ref{theorem_general_lepski}. Using the conditional independence of $Y$ and $A$ given $\bX$, one has the following calculations exactly parallel to that in proof of Theorem \ref{thm_treatment_effect}, that for a constant $C$ (depending on $M,B_U,B_L,\gmin,\gmax$), within the ``good'' event 
		\begin{align*}
		\mathcal{I}_{2; b, a, g} (n_2) \coloneqq \left\{ O \in \mathcal{D}_2: \begin{array}{c}
		\Vert \hat{b} - b \Vert_\infty \leq C_b^{\frac{d}{2 \beta_b + d}} \left( \frac{n_{2} / 2}{\log({n_{2} / 2})} \right)^{- \frac{\beta_b}{2 \beta_b + d}}
		\cap \; \Vert \hat{a} - a \Vert_\infty \leq C_a^{\frac{d}{2 \alpha + d}} \left( \frac{n_{2} / 2}{\log{(n_{2} / 2 )}} \right)^{- \frac{\alpha}{2 \alpha + d}} \cap \\
		\Vert \hat{g} - g \Vert_\infty \leq \tilde{C}^{\frac{d}{2 \alpha + d}} \left( \frac{n_{2} / 2}{\log{(n_{2} / 2 )}} \right)^{- \frac{\gamma}{2 \gamma + d}} \cap \hat{g} \in H(\gamma, C) \cap \hat{b} \in H(\beta, C) \cap \hat{a} \in H(\alpha, C)
		\end{array} \right\}
		\end{align*}
		with probability at least $1 - \frac{C \log{n_2}}{n_2^\eta}$ for some absolute constant $C$ and $\eta > 2$ by Theorem \ref{thm_adaptive_density_regression},
		\be 
		 \left\vert \E_{P, 1} \left(\phihat_{n, k}-\phi(P)\right) \right\vert \leq C \left[n^{-\frac{\alpha}{2\alpha+d}-\frac{\beta}{2\beta+d}-\frac{\gmin(\epsilon)}{2\gmin(\epsilon)+d}} + k^{-\frac{\alpha+\beta}{d}}\right],
		\ee
		where $\gmin(\epsilon) \coloneqq \frac{\gmin}{1 + \epsilon}$. Now, letting $\theta = (\a, \b, \g)$, $f_1(\theta) = \bav$ and $f_2(\theta)=-\frac{\alpha}{2\alpha+d}-\frac{\beta}{2\beta+d}-\frac{\gmin(\epsilon)}{2\gmin(\epsilon)+d}$ we have the bias property corresponding to Theorem \ref{theorem_general_lepski} holds with the given choice of $f_1$ and $f_2$ and a constant $C$ depending on $M, B_U, B_L,\gmax$ for $\{P\in \Par_{\theta}: f_1(\theta)=\bav,f_2(\theta)>\frac{2\alpha+2\beta}{2\alpha + 2\beta +d}\}$.  
		The proof of the validity of the conditional variance property corresponding to Theorem \ref{theorem_general_lepski} is once again easy to derive by standard Hoeffding decomposition of $\phihat_{n, k}$ followed by applications of moment bounds in Lemmas \ref{lemma_ustat_tail_use} and \ref{lemma_linear_tail_bound}.
		
		\be 
		\sup_{P \in \Par_{\theta}:\atop f_1(\theta)= \frac{\alpha + \beta}{2}, f_2(\theta)>\frac{2 \alpha + 2 \beta}{2 \alpha + 2 \beta + d}}\E_P\left(\phihat_{n, \khat} - \phi(P)\right)^2 \leq 8C \left(\frac{\sqrt{\log{n}}}{n}\right)^{\frac{4 \alpha + 4 \beta}{d + 2 \alpha + 2 \beta}}.
		\ee
		
		Noting that for $\theta \in \Theta$, since $\gmin>2(1+\epsilon)\max\{\alpha,\beta\}$, one has automatically, $f_2(\theta)>\frac{2\alpha+2\beta}{2\alpha+2\beta+d}$, completes the proof of the upper bound.\\

		\textit{\textbf{(ii) Proof of Lower Bound}} \\
		First note that we can parametrize our distributions by the tuple of functions $(a,b,g)$. We show that the same lower bound holds in a smaller class of problems where $g \equiv 1/2$ on $[0,1]^d$. Specifically consider
		\be 
		\Theta_{sub}=\left\{\begin{array}{c}P=(a,b,g):\\ a \in H(\alpha,M), b\in H(\beta, M),\  \frac{\alpha+\beta}{2}<\frac{d}{4},\\ g \equiv 1/2, (a(\bx),b(\bx))\in [B_L,B_U]^2 \ \forall \ \bx \in [0,1]^d \end{array}\right\}.
		\ee
		The observed data likelihood of $O \sim P$ for $P\in \Theta_{sub}$ can then be written as
		\be 
		(a(\bX)-1)^{1-A}\left(b^Y(\bX)(1-b(\bX))^{1-Y}\right)^A.
		\label{eqn:likelihood_missing_data}
		\ee
		Let for some $(\alpha,\beta,\gamma)$ tuple in the original problem $\Theta$, one has $$\sup_{P\in \Par_{(\alpha,\beta,\gamma)}}\E_P\left(\phihat-\phi(P)\right)^2 \leq C\left(\frac{\sqrt{\log{n}}}{n}\right)^{\frac{4\alpha+4\beta}{d+2\alpha+2\beta}}.$$
		
		Now, let $H:[0,1]^d\rightarrow \mathbb{R}$ be a $C^{\infty}$ function supported on $\left[0,\frac{1}{2}\right]^d$ such that $\int H(\bx)d\bx=0$ and $\int H^2(\bx)d\bx=1$ and let for $k \in \mathbb{N}$ (to be decided later) $\Omega_1,\ldots,\Omega_k$ be the translations of the cube $k^{-\frac{1}{d}}\left[0,\frac{1}{2}\right]^d$ that are disjoint and contained in $[0,1]^d$. Let  $\bx_1,\ldots,\bx_k$ denote the bottom left corners of these cubes.
		
		Assume first that $\alpha < \beta$. We set for $\lambda=(\lambda_1,\ldots,\lambda_k)\in \{-1,+1\}^k$ and $\alpha\leq \beta'<\beta$,
		$$a_{\lambda}(\bx)=2+\Big(\frac{1}{k}\Big)^{\frac{\alpha}{d}}\sum\limits_{j=1}^k\lambda_j H\Big((\bx-\bx_j)k^{\frac{1}{d}}\Big),$$
		
		$$b_{\lambda}(\bx)=\frac{1}{2}+\Big(\frac{1}{k}\Big)^{\frac{\beta'}{d}}\sum\limits_{j=1}^k\lambda_j H\Big((\bx-\bx_j)k^{\frac{1}{d}}\Big).$$
		
		A properly chosen $H$ guarantees $a_{\lambda}\in H(\alpha,M)$ and $b_{\lambda}\in H(\beta',M)$ for all $\lambda$. Let 
		\be 
		\Theta_0&=\left\{P^n: P=\left(a_{\lambda},1/2,1/2\right): \lambda \in \{-1,+1\}^k\right\},
		\ee
		and 
		\be 
		\Theta_1&=\left\{P^n: P=(a_{\lambda},b_{\lambda},1/2): \lambda \in \{-1,+1\}^k\right\}.
		\ee
		Finally let $\Theta_{test}=\Theta_0 \cup \Theta_1$. Let $\pi_0$ and $\pi_1$ be uniform priors on $\Theta_0$ and $\Theta_1$ respectively. It is easy to check that by our choice of $H$, $\phi(P)=\frac{1}{2}$ on $\Theta_0$ and $\phi(P)=\frac{1}{2}+\frac{1}{2}\left(\frac{1}{k}\right)^{\frac{\alpha+\beta'}{d}}$ for $P \in \Theta_1$. Therefore, using notation from Lemma \ref{lemma:prop_cai}, $\mu_1=\frac{1}{2}$, $\mu_2=\frac{1}{2}+\frac{1}{2}\left(\frac{1}{k}\right)^{\frac{\alpha+\beta'}{d}}$, and $\sigma_1=\sigma_2=0$. Since $\Theta_0 \subseteq  P(\alpha,\beta,\gamma)$, we must have that worst case error of estimation over $\Theta_0$ is bounded by $C\left(\frac{\sqrt{\log{n}}}{n}\right)^{\frac{4\alpha+4\beta}{d+2\alpha+2\beta}}$.  Therefore, the $\pi_0$ average bias over $\Theta_0$ is also bounded by $C\left(\frac{\sqrt{\log{n}}}{n}\right)^{\frac{2\alpha+2\beta}{d+2\alpha+2\beta}}$. This implies by Lemma \ref{lemma:prop_cai}, that the $\pi_1$ average bias over $\Theta_1$ (and hence the worst case bias over $\Theta_1$) is bounded below by a constant multiple of 
		\be 
		\left(\frac{1}{k}\right)^{\frac{\alpha+\beta'}{d}}-\left(\frac{\sqrt{\log{n}}}{n}\right)^{\frac{2\alpha+2\beta}{d+2\alpha+2\beta}}-\left(\frac{\sqrt{\log{n}}}{n}\right)^{\frac{2\alpha+2\beta}{d+2\alpha+2\beta}}\eta, \label{eqn:bias_lower_bound_missing_data}
		\ee
		where $\eta$ is the chi-square divergence between the probability measures $\int P^nd\pi_0(P^n)$ and $\int P^nd\pi_1(P^n)$. We now bound $\eta$ using Proposition \ref{prop_chisquare_affinity}.

		To put ourselves in the notation of Proposition \ref{prop_chisquare_affinity}, let for $\lambda \in \{-1,+1\}^k$, $P_{\lambda}$ and $Q_{\lambda}$ be the probability measures identified from $\Theta_0$ and $\Theta_1$ respectively.

		Therefore, with $\chi_j=\{0,1\}\times \{0,1\} \times \Omega_j$, we indeed have for all $j=1,\ldots,k$, $P_{\lambda}(\chi_j)=Q_{\lambda}(\chi_j)=p_j$ where there exists a constant $c$ such that $p_j=\frac{c}{k}$.

		Letting $\pi$ be the uniform prior over $\{-1,+1\}^k$ it is immediate that $\eta=\chi^2\left(\int P_{\lambda}d(\pi(\lambda)),\int Q_{\lambda}d(\pi(\lambda))\right)$. 
		
		It now follows by calculations similar to the proof of Theorem 4.1 in \cite{robins2009semiparametric}, that for a constant $C'>0$
		
		\be 
		\chi^2\left(\int P_{\lambda}d(\pi(\lambda)),\int Q_{\lambda}d(\pi(\lambda))\right)& \leq \exp\left(C'\frac{n^2}{k}\left(k^{-\frac{4\beta'}{d}}+k^{-4\frac{\alpha+\beta'}{2d}}\right)\right)-1.
		\ee
		Now choosing $k=\left(\frac{n}{\sqrt{c_*\log{n}}}\right)^{\frac{2d}{d+2\alpha+2\beta'}}$, we have 
		\be 
		\chi^2\left(\int P_{\lambda}d(\pi(\lambda)),\int Q_{\lambda}d(\pi(\lambda))\right)& \leq n^{2C'c_*}-1.
		\ee
		Therefore choosing $c_*$ such that $2C'c_*+\frac{2\alpha+2\beta'}{2\alpha+2\beta'+d}<\frac{2\alpha+2\beta}{2\alpha+2\beta+d}$, we have the desired result by \eqref{eqn:bias_lower_bound_treatment_effect}. The proof for $\alpha> \beta$ is similar after changing various quantities to:
		$$a_{\lambda}(\bx)=2+\Big(\frac{1}{k}\Big)^{\frac{\alpha'}{d}}\sum\limits_{j=1}^k\lambda_j H\Big((\bx-\bx_j)k^{\frac{1}{d}}\Big), \quad \beta\leq \alpha'<\alpha,$$
		
		$$b_{\lambda}(\bx)=\frac{1}{2}+\Big(\frac{1}{k}\Big)^{\frac{\beta}{d}}\sum\limits_{j=1}^k\lambda_j H\Big((\bx-\bx_j)k^{\frac{1}{d}}\Big).$$
		\be 
		\Theta_0&=\left\{P^n: P=\left(2,b_{\lambda},1/2\right): \lambda \in \{-1,+1\}^k\right\},
		\ee
		and 
		\be 
		\Theta_1&=\left\{P^n: P=(a_{\lambda},b_{\lambda},1/2): \lambda \in \{-1,+1\}^k\right\}.
		\ee
		For the case of $\alpha=\beta$, choose $\alpha'<\alpha$ and therefore, $\alpha'<\beta$ and thereafter work with
		$$a_{\lambda}(\bx)=2+\Big(\frac{1}{k}\Big)^{\frac{\alpha'}{d}}\sum\limits_{j=1}^k\lambda_j H\Big((\bx-\bx_j)k^{\frac{1}{d}}\Big),$$
		
		$$b_{\lambda}(\bx)=\frac{1}{2}+\Big(\frac{1}{k}\Big)^{\frac{\beta}{d}}\sum\limits_{j=1}^k\lambda_j H\Big((\bx-\bx_j)k^{\frac{1}{d}}\Big)$$

		\be 
		\Theta_0&=\left\{P^n: P=\left(a_{\lambda},1/2,1/2\right): \lambda \in \{-1,+1\}^k\right\},
		\ee
		and 
		\be 
		\Theta_1&=\left\{P^n: P=(a_{\lambda},b_{\lambda},1/2): \lambda \in \{-1,+1\}^k\right\}.
		\ee
		This completes the proof of the lower bound.
		
	\end{proof}
	
	\subsection*{Proof of Theorem \ref{thm_quadratic_functional_regression}}
	\begin{proof}
		\textit{\textbf{(i) Proof of Upper Bound}} \\
		We observe $n$ i.i.d. samples $\{Y_i, \bX_i\}_{i=1}^n$. First divide the samples into two disjoint parts $\mathcal{D} = \mathcal{D}_1 \coprod \mathcal{D}_2$ with sample sizes $n_1 = n ( 1 - 1 / \log{n} )$ and $n_2 = n / \log{n}$ respectively. Further divide the subsample $\mathcal{D}_2$ into two disjoint and equal-sized parts $\mathcal{D}_2 = \mathcal{D}_{21} \coprod \mathcal{D}_{22}$ with sample sizes $n_{21} = n_{22} = n / (2 \log{n})$. We estimate $g$ by $\ghat$ adaptively from $\mathcal{D}_{22}$ (as in Theorem \ref{thm_adaptive_density_regression}), and estimate $b$ by  $\bhat$, adaptively from $\mathcal{D}_{21}$ (as in Theorem \ref{thm_adaptive_density_regression}). For any $k \in \mathcal{K}$, consider
		\be
		\phihat_{n, k} &= \frac{1}{n}\sum_{i=1}^n {(2Y_i-\bhat(\bX_i))\bhat(\bX_i)}\\&+\frac{1}{n(n-1)}\sum\limits_{1\leq i_1\neq i_2\leq  n}\frac{(Y_{i_1}-\bhat(\bX_{i_1})))}{\sqrt{\ghat(\bX_{i_1})}}K_{V_k}(\bX_{i_1},\bX_{i_2})\frac{(Y_{i_2}-\bhat({\bX_{i_2}}))}{\sqrt{\ghat(\bX_{i_2})}}
		\ee
		
		Indeed this sequence of estimators is in the form of those considered by Theorem \ref{theorem_general_lepski} with $$\Ltilde_1(O)=(2Y-\bhat(\bX))\bhat(\bX),$$
		$$\Ltilde_{2l}(O)=-\frac{(Y-\bhat(\bX)))}{\sqrt{\ghat(\bX)}},\quad \Ltilde_{2r}(O)=\frac{(Y-\bhat(\bX)))}{\sqrt{\ghat(\bX)}},$$
		where by Theorem \ref{thm_adaptive_density_regression} $\max\{|\Ltilde_1(O)|, |\Ltilde_{2l}(O)|, |\Ltilde_{2r}(O)|\}\leq C(B_L, B_U)$. Therefore it remains to show that the sequence $\phihat_{n, k}$ satisfies the bias and variance properties (A) and (B) necessary for application of Theorem \ref{theorem_general_lepski}.
		We first verify the bias property. Utilizing the representation of the first order bias as stated above, we have
		\be 
		\ & |\E_{P, 1} \left(\phihat_{n, k}-\phi(P)\right)|\\
		& = \left|\begin{array}{c}
			\int\left((b(\bx)-\bhat(\bx)\right)^2g(\bx)d\bx \\ 
			- \E_{P, 1} \left[S\left(\frac{(Y_{1}-\bhat(\bX_{1})))}{\sqrt{\ghat(\bX_{1})}}K_{V_j}(\bX_{1},\bX_{2})\frac{(Y_{2}-\bhat({\bX_{2}}))}{\sqrt{\ghat(\bX_{2})}}\right)\right]\end{array}\right|\\ 
		\label{eqn:quadratic_second_orderbias_initial}
		\ee
		Now,
		by calculations similar to the proof of Theorem \ref{thm_treatment_effect}, one can show
		that for a constant $C$ (depending on $M,B_U,B_L,\gmin,\gmax$), within the ``good'' event 
		\begin{align*}
		\mathcal{I}_{2; b, g} (n_2) \coloneqq \left\{ O \in \mathcal{D}_2: \begin{array}{c}
		\Vert \hat{b} - b \Vert_\infty \leq C_b^{\frac{d}{2 \beta_b + d}} \left( \frac{n_{2} / 2}{\log({n_{2} / 2})} \right)^{- \frac{\beta_b}{2 \beta_b + d}}
		\cap \; \Vert \hat{g} - g \Vert_\infty \leq \tilde{C}^{\frac{d}{2 \gamma + d}} \left( \frac{n_{2} / 2}{\log{(n_{2} / 2 )}} \right)^{- \frac{\gamma}{2 \gamma + d}} \\
		\cap \; \hat{g} \in H(\gamma, C) \cap \hat{b} \in H(\beta, C)
		\end{array} \right\}.
		\end{align*}
		with probability at least $1 - \frac{C \log{n_2}}{n_2^\eta}$ for some absolute constant $C$ and $\eta > 2$ by Theorem \ref{thm_adaptive_density_regression},
		\be 
		\ & \left\vert \E_{P, 1} \left(\phihat_{n, k}-\phi(P)\right) \right\vert \\
		& \leq C \left[n^{-\frac{2\beta}{2\beta+d}-\frac{\gmin(\epsilon)}{2\gmin(\epsilon)+d}} +k^{-\frac{\alpha+\beta}{d}}\right], \quad \gmin(\epsilon) \coloneqq \frac{\gmin}{1+\epsilon}
		\ee
		Now, letting $\theta=(\b,\g)$, $f_1(\theta)=\beta$ and $f_2(\theta)=-\frac{2\beta}{2\beta+d}-\frac{\gmin(\epsilon)}{2\gmin(\epsilon)+d}$, the rest of the proof follows along the lines of the proof of Theorem \ref{thm_treatment_effect}.
		\\
		\textit{\textbf{(ii) Proof of Lower Bound}} \\
		The proof of the lower bound is very similar to that of the lower bound proof in Theorem \ref{thm_treatment_effect}, after identifying $Y=A$ almost surely $\P$, and hence is omitted.
	\end{proof}
	\section{Technical Lemmas}\label{section_appendix_technical_lemmas}
	\subsection{Constrained Risk Inequality}
	
	A main tool for producing adaptive lower bound arguments is a general version of constrained risk inequality due to \cite{cai2011testing}, obtained as an extension of \cite{brown1996constrained}. For the sake of completeness, we begin with a summary of these results. Suppose $Z$ has distribution $\mathbb{P}_{\theta}$ where $\theta$ belongs to some parameter space $\Theta$. Let $\hat{Q}=\hat{Q}(Z)$ be an estimator of a function $Q(\theta)$ based on $Z$ with bias $B(\theta) \coloneqq \E_{\theta}(\hat{Q})-Q(\theta)$. Now suppose that $\Theta_0$ and $\Theta_1$ form a disjoint partition of $\Theta$ with priors $\pi_0$ and $\pi_1$ supported on them respectively. Also, let $\mu_i=\int Q(\theta)d\pi_i$ and $\sigma_i^2=\int (Q(\theta)-\mu_i)^2d\pi_i$, $i=0,1$ be the mean and variance of $Q(\theta)$ under the two priors $\pi_0$ and $\pi_1$. Letting $\gamma_i$ be the marginal density w.r.t. some common dominating measure of $Z$ under $\pi_i$, $i=0,1$, let us denote by $\E_{\gamma_0}(g(Z))$ the expectation of $g(Z)$ w.r.t. the marginal density of $Z$ under prior $\pi_0$ and distinguish it from $\E_{\theta}(g(Z))$, which is the expectation under $\mathbb{P}_{\theta}$. Lastly, denote the chi-square divergence between $\gamma_0$ and $\gamma_1$ by $\chi=\left\{\E_{\gamma_0}\left(\frac{\gamma_1}{\gamma_0}-1\right)^2\right\}^{\frac{1}{2}}$. Then we have the following result.
	\begin{lem}[\cite{cai2011testing}]\label{lemma:prop_cai}
		If $\int \E_{\theta}\left(\hat{Q}(Z)-Q(\theta)\right)^2d\pi_0(\theta)\leq \epsilon^2$, then
		$$\left\vert\int B(\theta)d\pi_1(\theta)-\int B(\theta)d\pi_0(\theta)\right\vert\geq |\mu_1-\mu_0|-(\epsilon+\sigma_0)\chi.$$
	\end{lem}
	Since the maximum risk is always at least as large as the average risk, this immediately yields a lower bound on the minimax risk.

	\subsection{Tail and Moment Bounds}
	The U-statistics appearing in this paper are mostly based on projection kernels sandwiched between arbitrary bounded functions. This necessitates generalizing the U-statistics bounds obtained in \cite{bull2013adaptive} as in \cite{mukherjee2018optimal} . 
	
	\begin{lem}\label{lemma_ustat_tail_use}
	\phantomsection
		$\bW_1,\ldots,\bW_n \sim \P$ are i.i.d. random vectors of observations such that $\bX_i\in [0,1]^d$ is a sub-vector of $\bW_i$ for each $i$. There exists a constant $C \coloneqq C(B,B_U,J_0)>0$ such that the following hold
		\begin{enumerate}[(i)]
			\item \begin{align}
			\ & \mathbb{P}\Big(\Big|\frac{1}{n(n-1)}\sum_{i_1 \neq i_2}R\left(\bW_{i_1},\bW_{i_2}\right)-\E\left(R\left(\bW_1,\bW_2\right)\right) \Big|\geq t\Big) \\& \leq e^{-C nt^2}+e^{-\frac{C t^2}{a_1^2}}+e^{-\frac{C t}{a_2}}+e^{-\frac{C \sqrt{t}}{\sqrt{a_3}}}, \nonumber
			\end{align}
			\item \be 
			\E\left(|\frac{1}{n(n-1)}\sum_{i_1 \neq i_2}R\left(\bW_{i_1},\bW_{i_2}\right)-\E\left(R\left(\bW_1,\bW_2\right)\right) |^{2q}\right)&\leq \left(C\frac{2^{jd}}{n^2}\right)^q,
			\ee
		\end{enumerate}
		where $a_1= \frac{1}{n-1}2^{\frac{jd}{2}}$, $a_2=\frac{1}{n-1}\left(\sqrt{\frac{2^{jd}}{n}}+1\right)$, $a_3=\frac{1}{n-1}\left(\sqrt{\frac{2^{jd}}{n}}+\frac{2^{jd}}{n}\right)$,\\ $R(\bW_1,\bW_2)=S\left(\Ltilde_{2l}(O_{1})K_{V_j} \left(\bX_{1},\bX_{2}\right)\Ltilde_{2r}(O_{2}) \right)$ with
		$\max\{|\Ltilde_{2l}(O)|, |\Ltilde_{2r}(O)|\}\leq B$ , almost surely $\bW$, and $\bX_i \in [0,1]^d$ are i.i.d. with density $g$ such that $g(\bx)\leq B_U$ for all $\bx \in [0,1]^d$.

	\end{lem}
	\begin{proof}
		The proof of part (i) can be found in \cite{mukherjee2018optimal}. However, for the sake of completeness we provide the proof here again. We do the proof for the special case where $\Ltilde_{2l}=\Ltilde_{2r}=L$. However, the details of the argument show that the proof continues to hold to symmetrized U-statistics as defined here.
		
		The proof hinges on the following tail bound for second-order degenerate U-statistics \citep{gine2016mathematical} which is due to \cite{gine2000exponential} with constants by \cite{houdre2003exponential} and is crucial for our calculations.

		\begin{lem}\label{lemma_ustat_tail}
		\phantomsection
			Let $U_n$ be a degenerate U-statistic of order 2 with kernel $R$ based on an i.i.d. sample $W_1,\ldots,W_n$. Then there exists a constant $C$ independent of $n$, such that 
			\begin{align*}
			P[ | \sum_{i \neq j } R(W_1, W_2) | \geq C( \Lambda_1 \sqrt{u} + \Lambda_2  u + \Lambda_3 u^{3/2} + \Lambda_4 u^2 )] \leq 6 \exp(- u), 
			\end{align*}
			where, we have,
			\begin{align*}
			\Lambda_1^2 &= \frac{n(n-1)}{2} E[ R^2(W_1 , W_2)], \nonumber\\
			\Lambda_2 &= n \sup \{E[R(W_1, W_2) \zeta(W_1) \xi(W_2)] : E[\zeta^2(W_1)]\leq 1, E[\xi^2(W_1)] \leq 1 \}, \nonumber\\
			\Lambda_3 &= \| n E[R^2 (W_1 , \cdot) \| _{\infty} ^{\frac{1}{2}}, \nonumber \\ 
			\Lambda_4 &= \| R \| _{\infty}. \nonumber
			\end{align*}
		\end{lem}

		We use this lemma to establish Lemma \ref{lemma_ustat_tail_use}.
		By Hoeffding's decomposition one has
		\be
		\ &\frac{1}{n(n-1)}\sum_{i_1 \neq i_2}R\left(\bW_{i_1},\bW_{i_2}\right)-\E\left(R\left(\bW_1,\bW_2\right)\right)
		\\&=\frac{2}{n}\sum_{i_1=1}^n \left[\E_{\bW_{i_1}}R\left(\bW_{i_1},\bW_{i_2}\right)-\E R\left(\bW_{i_1},\bW_{i_2}\right)\right]\\
		&+\frac{1}{n(n-1)}\sum_{i_1 \neq i_2}\left[\begin{array}{c}R\left(\bW_{i_1},\bW_{i_2}\right)-\E_{\bW_{i_1}}R\left(\bW_{i_1},\bW_{i_2}\right)\\-\E_{\bW_{i_2}}R\left(\bW_{i_1},\bW_{i_2}\right)+\E R\left(
			\bW_{i_1},\bW_{i_2}\right)\end{array}\right]\\
		& \coloneqq T_1+T_2
		\ee
		\subsubsection{Analysis of $T_1$}
		Noting that $T_1=\frac{2}{n}\sum_{i_1=1}^n H(\bW_{i_1})$ where $H(\bW_{i_1})=\E\left(R\left(\bW_{i_1},\bW_{i_2}|\bW_{i_1}\right)\right)-\E R\left(\bW_{i_1},\bW_{i_2}\right)$ we control $T_1$ by standard Hoeffding's inequality. First note that, 
		\be
		\ &|H(\bW_{i_1})|\\
		&=|\sum_{k \in \Z_j}\sum_{v \in \{0,1\}^d}\left[L\left(\bW_{i_1}\right)\psi_{jk}^v\left(\bX_{i_1}\right)\E\left(\psi_{jk}^v\left(\bX_{i_2}\right)L\left(\bW_{i_2}\right)\right)-\left(\E\left(\psi_{jk}^v\left(\bX_{i_2}\right)L\left(\bW_{i_2}\right)\right)\right)^2\right]|\\
		& \leq \sum_{k \in \Z_j}\sum_{v \in \{0,1\}^d}|L\left(\bW_{i_1}\right)\psi_{jk}^v\left(\bX_{i_1}\right)\E\left(\psi_{jk}^v\left(\bX_{i_2}\right)L\left(\bW_{i_2}\right)\right)|\\&+\sum_{k \in \Z_j}\sum_{v \in \{0,1\}^d}\left(\E\left(\psi_{jk}^v\left(\bX_{i_2}\right)L\left(\bW_{i_2}\right)\right)\right)^2
		\ee
		First, by standard compactness argument for the wavelet bases,
		\be
		|\E\left(\psijkv\left(\bX\right)L(\bW)\right)|&\leq \int|\E\left(L(\bW)|\bX=\bx\right)\Big(2^{\frac{jd}{2}}\prod\limits_{l=1}^d\psi_{00}^{v_l}(2^j x_l-k_l )\Big)||g(\bx)|d\bx \\
		&\leq \constant 2^{-\frac{jd}{2}}. \label{eqn:innerproductbound}
		\ee

		Therefore,
		\be
		\sum_{k \in \Z_j}\sum_{v \in \{0,1\}^d}\left(\E\left(\psi_{jk}^v\left(\bX_{i_2}\right)L\left(\bW_{i_2}\right)\right)\right)^2 & \leq \constant. \label{eqn:expected_kernel}
		\ee
		Also, using the fact that for each fixed $\bx \in [0,1]^d$, the number of indices $k \in \Z_j$ such that $\bx$ belongs to support of at least one of $\psi_{jk}^v$ is bounded by a constant depending only on $\psi^0_{00}$ and $\psi_{00}^1$. Therefore combining \eqref{eqn:innerproductbound} and \eqref{eqn:expected_kernel},
		\be
		\ &\sum_{k \in \Z_j}\sum_{v \in \{0,1\}^d}|L\left(\bW_{i_1}\right)\psi_{jk}^v\left(\bX_{i_1}\right)\E\left(\psi_{jk}^v\left(\bX_{i_2}\right)L\left(\bW_{i_2}\right)\right)|\\
		&  \leq \constant 2^{-\frac{jd}{2}}2^{\frac{jd}{2}}=\constant. \\ \label{eqn:linear_term_uniform_bound}
		\ee
		Therefore, by \eqref{eqn:linear_term_uniform_bound} and Hoeffding's inequality,
		\be
		\P\left(|T_1|\geq t\right)\leq 2e^{-\constant nt^2}. \label{ustat_linear_term_exponential_bound}
		\ee
		\subsubsection{Analysis of $T_2$}
		Since $T_2$ is a degenerate U-statistic, its analysis is based on Lemma \ref{lemma_ustat_tail}. In particular,
		\be
		T_2&= \frac{1}{n(n-1)}\sum_{i_1 \neq i_2}R^{*}\left(\bW_{i_1},\bW_{i_2}\right)
		\ee
		where 
		\be 
		\ &R^{*}\left(\bW_{i_1},\bW_{i_2}\right)\\
		&=\sum_{k \in \Z_j}\sum_{v \in \{0,1\}^d}\left\{\begin{array}{c}\left(L(\bW_{i_1})\psi_{jk}^v\left(\bX_{i_1}\right)-\E\left(\psi_{jk}^v\left(\bX_{i_1}\right)\E\left(L(\bW_{i_1})|\bX_{i_1}\right)\right)\right)\\ \times \left(L(\bW_{i_2})\psi_{jk}^v\left(\bX_{i_2}\right)-\E\left(\psi_{jk}^v\left(\bX_{i_2}\right)\E\left(L(\bW_{i_2})|\bX_{i_2}\right)\right)\right)
		\end{array}\right\}.
		\ee
		Letting $\Lambda_i$, $i=1,\ldots,4$ be the relevant quantities as in Lemma \ref{lemma_ustat_tail}, we have the following lemma.
		\begin{lem}\label{lemma_lamda_estimates}
			There exists a constant $C=\constant$ such that
			$$\Lambda_1^2 \leq C \frac{n(n-1)}{2}2^{jd},\ \Lambda_2 \leq C n,\ \Lambda_3^2 \leq C n 2^{jd},\ \Lambda_4 \leq C 2^{\frac{jd}{2}}.$$
		\end{lem}	
		\begin{proof}
			First we control $\Lambda_1$. To this end, note that by simple calculations, using bounds on $L,g$, and orthonormality of $\psijkv$'s we have,
			\be
			\  \Lambda_1^2 &=\frac{n(n-1)}{2}\E \left(\left\{R^{*}\left(\bW_1,\bW_2\right)\right\}^2\right)\leq 3n(n-1)\E\left(R^2\left(\bW_1,\bW_2\right)\right)\\
			&=3n(n-1)\E \left(L^2\left(\bW_1\right)K_{V_j}^2\left(\bX_1,\bX_2\right)L^2\left(\bW_2\right)\right)\\&\leq 3n(n-1)B^4 \int \int \Big[\sumk \sumv \psijkv\left(\bx_1\right)\psijkv\left(\bx_2\right)\Big]^2 g(\bx_1)g(\bx_2)d\bx_1 d\bx_2\\
			& \leq 3n(n-1)B^4 B_U^2 \int \int \Big[\sumk \sumv \psijkv\left(\bx_1\right)\psijkv\left(\bx_2\right)\Big]^2 d\bx_1 d\bx_2\\&= 3n(n-1)B^4 B_U^2\sumk \sumv\int  \left(\psijkv\left(\bx_1\right)\right)^2d\bx_1\int  \left(\psijkv\left(\bx_2\right)\right)^2d\bx_2\\
			& \leq \constant n(n-1)2^{jd}.
			\ee
			Next we control
			$$\Lambda_2=n\sup\left\{\E \left(R^*\left(\bW_1,\bW_2\right)\zeta\left(\bW_1\right)\xi\left(\bW_2\right)\right): \E \left(\zeta^2(\bW_1)\right)\leq 1, \E \left(\xi^2(\bW_2)\right)\leq 1\right\}.$$
			
			To this end, we first control
			\be 
			\ & |\E\left(L(\bW_1)K_{V_j}\left(\bX_1,\bX_2\right)L(\bW_2)\zeta(\bW_1)\xi(\bW_2)\right)|\\
			&=|\int \int \E(L(\bW_1)\zeta(\bW_1)|\bX_1=\bx_1)K_{V_j}\left(\bx_1,\bx_2\right)\E(L(\bW_2)\xi(\bW_2)|\bX_2=\bx_2)g(\bx_1)g(\bx_2)d\bx_1 d\bx_2|\\
			&=|\int \E(L(\bW)\zeta(\bW)|\bX=\bx) \Pi\left(\E(L(\bW)\xi(\bW)|\bX=\bx)g(\bx)|V_j\right)g(\bx)d\bx|\\
			& \leq \left(\int \E^2(L(\bW)\zeta(\bW)|\bX=\bx)g^2(\bx)d\bx\right)^{\frac{1}{2}}\left(\int \Pi^2\left(\E(L(\bW)\xi(\bW)|\bX=\bx)g(\bx)|V_j\right)d\bx\right)^{\frac{1}{2}} \\
			& \leq \left(\int \E(L^2(\bW)\zeta^2(\bW)|\bX=\bx)g^2(\bx)d\bx\right)^{\frac{1}{2}}\left(\int \E(L^2(\bW)\xi^2(\bW)|\bX=\bx)g^2(\bx)d\bx \right)^{\frac{1}{2}}\\
			& \leq B^2 B_U \sqrt{\E(\zeta^2(\bW_1))\E (\xi^2(\bW_2))}\leq B^2 B_U.
			\ee
			Above we have used Cauchy-Schwarz inequality, Jensen's inequality, and the fact that projections contract norms. Also,
			\be
			\ & |\E \left(\E \left(L(\bW_1)K_{V_j}\left(\bX_1,\bX_2\right)L(\bW_2)|\bW_1\right)\zeta(\bW_1)\xi(\bW_2)\right)|\\
			& = |\E \left[L(\bW_1)\Pi\left(\E \left(L(\bW_1)g(\bX_1)|\bX_1\right)|V_j\right)\zeta(\bW_1)\xi(\bW_2)\right]|\\
			& = |\E \left[L(\bW_1)\Pi\left(\E \left(L(\bW_1)g(\bX_1)|\bX_1\right)|V_j\right)\zeta(\bW_1)\right]||\E(\xi(\bW_2))|\\
			& \leq |\int \Pi(\E (L(\bW)\zeta(\bW)|\bX=\bx)g(\bx)|V_j)\Pi(\E (L(\bW)|\bX=\bx)g(\bx)|V_j)d\bx| \leq B^2 B_U,
			\ee
			where the last step once again uses contraction property of projection, Jensen's inequality, and bounds on $L$ and $g$. Finally, by Cauchy-Schwarz inequality and \eqref{eqn:expected_kernel},
			\be
			\ &\E \left[\E \left(L(\bW_1)K_{V_j}\left(\bX_1,\bX_2\right)L(\bW_2)\right)\zeta(\bW_1)\xi(\bW_2)\right]\\
			& \leq \sumk \sumv \E^2 \left(L(\bW)\psijkv(\bX)\right)\leq \constant.
			\ee
			This completes the proof of $\Lambda_2 \leq \constant n$.
			Turning to $\Lambda_3=n \|\E\left[\left(R^*(\bW_1,\cdot)\right)^2\right]\|^{\frac{1}{2}}_{\infty}$ we have that
			\be 
			\ &\left(R^*(\bW_1,\bw_2)\right)^2 \\& \leq 2 \left[R(\bW_1,\bw_2)-\E (R(\bW_1,\bW_2)|\bW_1)\right]^2+2\left[\E (R(\bW_1,\bW_2)|\bW_2=\bw_2)-\E \left(R(\bW_1,\bW_2)\right)\right]^2.
			\ee
			Now,
			\be 
			\ & \E \left[R(\bW_1,\bw_2)-\E (R(\bW_1,\bW_2)|\bW_1)\right]^2 \\& \leq 2 \E \left(L^2(\bW_1)K^2_{V_j}\left(\bX_1,\bx_2\right)L^2(\bw_2)\right)
			\\&+2\E \Big(\sumk\sumv L(\bW_1)\psijkv(\bX_1)\E \left(\psijkv(\bX_2)L(\bW_2)\right)\Big)^2\\
			& \leq 2B^4 B_U^2 \sumk\sumv \left(\psijkv(\bx_2)\right)^2+2\E (H^2(\bW_2))\leq \constant 2^{jd}
			\ee
			where the last inequality follows from arguments along the line of \eqref{eqn:linear_term_uniform_bound}. Also, using inequalities \eqref{eqn:expected_kernel} and \eqref{eqn:linear_term_uniform_bound}
			\be
			\ & \left[\E (R(\bW_1,\bW_2)|\bW_2=\bw_2)-\E \left(R(\bW_1,\bW_2)\right)\right]^2\\&=\Big[\sumk \sumv \E \left(L(\bW_1)\psijkv(\bX_1)\right)\left(\E\left(L(\bW_1)\psijkv(\bX_1)\right)-\psijkv(\bx_2)L(\bw_2)\right)\Big]^2\\&\leq \constant.
			\ee
			This completes the proof of controlling $\Lambda_3$. Finally, using compactness of the wavelet basis,
			\be
			\|R(\cdot,\cdot)\|_{\infty}& \leq B^2 \sup_{\bx_1,\bx_2}\sumk \sumv |\psijkv(\bx_1)||\psijkv(\bx_2)|\leq \constant 2^{jd}.
			\ee
			Combining this with arguments similar to those leading to \eqref{eqn:linear_term_uniform_bound}, we have $\Lambda_4\leq \constant 2^{jd}$.
		\end{proof}
		Therefore, using Lemma \ref{lemma_ustat_tail} and Lemma \ref{lemma_lamda_estimates} we have 
		\be 
		\P\Big(|T_2|\geq \frac{\constant}{n-1}\Big(\sqrt{2^{jd}t}+t+\sqrt{\frac{2^{jd}}{n}}t^{\frac{3}{2}}+\frac{2^{jd}}{n}t^2\Big)\Big)& \leq 6 e^{-t}.
		\ee
		Finally using $2t^{\frac{3}{2}}\leq t+t^2$ we have,
		\be
		\Pf\left[|T_2|> a_1\sqrt{t}+a_2t+a_3 t^2\right]& \leq 6 e^{-t} \label{ustat_quadratic_term_bound}
		\ee
		where $a_1= \frac{\constant}{n-1}2^{\frac{jd}{2}}$, $a_2=\frac{\constant}{n-1}\left(\sqrt{\frac{2^{jd}}{n}}+1\right)$, and $a_3=\frac{\constant}{n-1}\left(\sqrt{\frac{2^{jd}}{n}}+\frac{2^{jd}}{n}\right)$. Now if $h(t)$ is such that $a_1\sqrt{h(t)}+a_2 h(t)+a_3h^2(t)\leq t$, then one has by \eqref{ustat_quadratic_term_bound},
		\be
		\mathbb{P}\left[|T_2|\geq t\right]& \leq \mathbb{P}\left[|T_2|\geq a_1\sqrt{h(t)}+a_2 h(t)+a_3h^2(t)\right]\leq 6 e^{-6h(t)}.
		\ee
		Indeed, there exists such an $h(t)$ such that $h(t)=b_1 t^2 \wedge b_2 t \wedge b_3 \sqrt{t}$ where $b_1=\frac{\constant}{a_1^2}$, $b_2=\frac{\constant}{a_2}$, and $b_3=\frac{\constant}{\sqrt{a_3}}$. Therefore, there exists a $C=\constant$ such that 
		\be 
		\mathbb{P}\left[|T_2|\geq t\right]\leq e^{-\frac{C t^2}{a_1^2}}+e^{-\frac{C t}{a_2}}+e^{-\frac{C \sqrt{t}}{\sqrt{a_3}}}. \label{ustat_quadratic_term_bound_use}
		\ee
		
		\subsubsection{Combining Bounds on $T_1$ and $T_2$} Applying union bound along with \ref{ustat_linear_term_exponential_bound} and \ref{ustat_quadratic_term_bound_use} completes the proof of Lemma \ref{lemma_ustat_tail_use} part (i).
		
		For the proof of part (ii) note that with the notation of the proof of part (i) we have by Hoeffding decomposition
		\be 
		\E\left(|\frac{1}{n(n-1)}\sum_{i_1 \neq i_2}R\left(\bW_{i_1},\bW_{i_2}\right)-\E\left(R\left(\bW_1,\bW_2\right)\right) |^{2q}\right)&\leq 2(\E|T_1|^{2q}+\E|T_2|^{2q}).
		\ee
		The proof will be completed by individual control of the two moments above.
		\subsubsection{Analysis of $\E|T_1|^{2q}$}
		Recall that $T_1=\frac{2}{n}\sum_{i_1=1}^n H(\bW_{i_1})$ where $H(\bW_{i_1})=\E\left(R\left(\bW_{i_1},\bW_{i_2}|\bW_{i_1}\right)\right)-\E R\left(\bW_{i_1},\bW_{i_2}\right)$ and $|H(\bW)|\leq \constant$ almost surely. Therefore by Rosenthal's inequality \eqref{lemma_linear_tail_bound} we have
		\be 
		\E|T_1|^{2q}&\leq \left(\frac{2}{n}\right)^{2q}\left[\sum_{i=1}^n \E|H(\bW_i)|^{2q}+\left\{\sum_{i=1}^{n}\E|H(\bW_i)|^{2}\right\}^{q}\right]\\
		&\leq \left(\frac{2\constant}{n}\right)^{2q}\left(n+n^{q}\right)\leq C(B,B_U,J_0)^qn^{-q}.
		\ee
		\subsubsection{Analysis of $\E|T_2|^{2q}$}
		Recall that
		\be
		\mathbb{P}\left[|T_2|\geq t\right]& \leq \mathbb{P}\left[|T_2|\geq a_1\sqrt{h(t)}+a_2 h(t)+a_3h^2(t)\right]\leq 6 e^{-6h(t)}
		\ee
		where $h(t)=b_1 t^2 \wedge b_2 t \wedge b_3 \sqrt{t}$ with $b_1=\frac{\constant}{a_1^2}$, $b_2=\frac{\constant}{a_2}$, and $b_3=\frac{\constant}{\sqrt{a_3}}$. Therefore
		\be
		\ &\E (|T_2|^{2q})\\&=2q\int_0^{\infty}x^{2q-1}\Pf(|T_2|\geq x)dx\\
		&\leq 2q \int_0^{\infty}x^{2q-1}\Pf(|T_2|\geq a_1\sqrt{h(x)}+a_2 h(x) +a_3 h^2(x))dx\\
		&\leq 12q \int_0^{\infty}x^{2q-1} e^{-h(x)}dx\\
		& = 12q \int_0^{\infty}x^{2q-1} e^{-\left\{b_1 x^2 \wedge b_2 x \wedge b_3 \sqrt{x}\right\}}dx\\
		&\leq 12q \left[\int_0^{\infty}x^{2q-1} e^{-b_1 x^2}dx+\int_0^{\infty}x^{2q-1} e^{- b_2 x}dx+\int_0^{\infty}x^{2q-1} e^{-b_3 \sqrt{x}}dx\right]\\
		&=12q\left(\frac{\Gamma(q)}{2b_1^{q}}+\frac{\Gamma(2q)}{b_2^{2q}}+\frac{2\Gamma(4q)}{b_3^{4q}}\right)\leq  \left(C\frac{2^{jd}}{n^2}\right)^{q}
		\ee
		for a constant $C=\constant$, by our choices of $b_1,b_2,b_3$.
	\end{proof}

	Since the estimators arising in this paper also have a linear term, we will need the following standard Bernstein and Rosenthal type tail and moment bounds \citep{petrov1995limit}.
	\begin{lem}\label{lemma_linear_tail_bound}
	\phantomsection
		If $\bW_1,\ldots,\bW_n \sim \P$ are i.i.d. random vectors such that $|L(\bW)|\leq B$ almost surely $\P$, then for $q\geq 2$ one has for large enough constants $C(B)$ and $C(B,q)$
		\be 
		\P(|\frac{1}{n}\sum_{i=1}^n\left(L(\bW_i)-\E(L(\bW_i))\right)|\geq t)\leq 2e^{-nt^2/C(B)},
		\ee
		and 
		\be 
		\ & \E(|\sum_{i=1}^n\left(L(\bW_i)-\E(L(\bW_i))\right)|^q)\\&\leq \left[\sum_{i=1}^n \E\left(|L(\bW_i)-\E(L(\bW_i))|^q\right)+\left[\sum_{i=1}^n \E\left(|L(\bW_i)-\E(L(\bW_i))|^2\right)\right]^{q/2}\right]\\&\leq C(B,q)n^{\frac{q}{2}}.
		\ee
	\end{lem}
	We will also need the following concentration inequality for linear estimators based on wavelet projection kernels, the proof of which can be done along the lines of proofs of Theorem 5.1.5 and Theorem 5.1.13 of \cite{gine2016mathematical}.
	\begin{lem}\label{lemma_linear_projection_tail_bound}
		Consider i.i.d. observations $\bW_i=(Y,\bX)_i$, $i=1,\ldots,n$ where $\bX_i\in [0,1]^d$ with marginal density $g$. Let $\hat{m}(\bx)=\frac{1}{n}\sum_{i=1}^n L(\bW_i)K_{V_l}\left(\bX_i,\bx\right)$, such that $\max\{\|g\|_{\infty},\|L\|_{\infty}\}\leq B_U$. If $\frac{2^{ld}ld}{n}\leq 1$,  there exist $C,C_1,C_2>0$, depending on $B_U$ and scaling functions $\psi_{0,0}^0,\psi_{0,0}^1$ respectively, such that 
		\be 
		\E(\|\hat{m}-\E(\hat{m})\|_{\infty})\leq C\sqrt{\frac{2^{ld}ld}{n}},
		\ee
		and for any $t>0$
		\be 
		\P\left(n\|\hat{m}-\E(\hat{m})\|_{\infty}>\frac{3}{2}n\E(\|\hat{m}-\E(\hat{m})\|_{\infty})+\sqrt{C_1 n 2^{ld}t}+C_2 2^{ld}t\right)\leq e^{-t}.
		\ee
	\end{lem}
\end{document}